\documentclass[sigconf,authorversion=true]{acmart}

\usepackage{xcolor}
\usepackage{adjustbox}
\usepackage{colortbl}
\usepackage{tabularx}
\usepackage{booktabs} 
\usepackage{graphicx}
\usepackage{subfig}
\usepackage[linesnumbered, vlined, ruled]{algorithm2e}
\SetArgSty{textup}

\usepackage{color, colortbl}
\usepackage{braket}
\definecolor{MyHiLiRow}{gray}{0.9}

\newcommand{\argmin}{\mathop{\rm argmin}\limits}

\def\vector#1{\mbox{\boldmath $#1$}}

\AtBeginDocument{%
  \providecommand\BibTeX{{%
    \normalfont B\kern-0.5em{\scshape i\kern-0.25em b}\kern-0.8em\TeX}}}

\setcopyright{acmcopyright}

\copyrightyear{2022}
\acmYear{2022}
\setcopyright{acmcopyright}\acmConference[GECCO '22]{Genetic and Evolutionary Computation Conference}{July 9--13, 2022}{Boston, MA, USA}
\acmBooktitle{Genetic and Evolutionary Computation Conference (GECCO '22), July 9--13, 2022, Boston, MA, USA}
\acmPrice{15.00}
\acmDOI{10.1145/3512290.3528778}
\acmISBN{978-1-4503-9237-2/22/07}

\begin{document}

\title[A Two-phase Framework with a B\'{e}zier Simplex-based Interpolation Method]{A Two-phase Framework with a B\'{e}zier Simplex-based Interpolation Method for Computationally Expensive Multi-objective Optimization}





\author{Ryoji Tanabe}
\affiliation{%
  \institution{Yokohama National University \& RIKEN AIP}
  \city{Yokohama}
  \state{Kanagawa}
  \country{Japan}}
  \email{rt.ryoji.tanabe@gmail.com}

\author{Youhei Akimoto}
\orcid{0000-0003-2760-8123}
\affiliation{%
  \institution{University of Tsukuba \& RIKEN AIP}
  \city{Tsukuba}
  \state{Ibaraki}
  \country{Japan}}
  \email{akimoto@cs.tsukuba.ac.jp}
  
\author{Ken Kobayashi}
\affiliation{%
  \institution{Fujitsu Limited \& RIKEN AIP}
  \city{Kawasaki}
  \state{Kanagawa}
  \country{Japan}}
  \email{ken-kobayashi@fujitsu.com}

\author{Hiroshi Umeki}
\affiliation{%
  \institution{Yokohama National University \& RIKEN AIP}
  \city{Yokohama}
  \state{Kanagawa}
  \country{Japan}}
  \email{drumehiron@gmail.com}
 
\author{Shinichi Shirakawa}
\affiliation{%
  \institution{Yokohama National University \& RIKEN AIP}
  \city{Yokohama}
  \state{Kanagawa}
  \country{Japan}}
  \email{shirakawa-shinichi-bg@ynu.ac.jp}

\author{Naoki Hamada}
\affiliation{%
  \institution{KLab Inc. \& RIKEN AIP}
  \city{Minato}
  \state{Tokyo}
  \country{Japan}}
  \email{hamada-n@klab.com}






\renewcommand{\shortauthors}{R. Tanabe, Y. Akimoto, K. Kobayashi, H. Umeki, S. Shirakawa, and N. Hamada}


\begin{abstract}

This paper proposes a two-phase framework with a B\'{e}zier simplex-based interpolation method (TPB) for computationally expensive multi-objective optimization.
The first phase in TPB aims to approximate a few Pareto optimal solutions by optimizing a sequence of single-objective scalar problems.
The first phase in TPB can fully exploit a state-of-the-art single-objective derivative-free optimizer.
The second phase in TPB utilizes a B\'{e}zier simplex model to interpolate the solutions obtained in the first phase.
The second phase in TPB fully exploits the fact that a B\'{e}zier simplex model can approximate the Pareto optimal solution set by exploiting its simplex structure when a given problem is simplicial.
We investigate the performance of TPB on the 55 bi-objective BBOB problems.
The results show that TPB performs significantly better than HMO-CMA-ES and some state-of-the-art meta-model-based optimizers.

\end{abstract}




\maketitle

\section{Introduction}
\label{sec:introduction}

\textit{General context.}
This paper considers computationally expensive multi-objective black-box numerical optimization.
Some real-world optimization problems require computationally expensive simulation to evaluate the solution (e.g., \cite{DanielsRETF18,YangPESB19}).
In this case, only a limited budget of function evaluations is available for multi-objective optimization.
%
Instead of general evolutionary multi-objective optimization (EMO) algorithms (e.g., NSGA-II \cite{DebAPM02} and MOEA/D \cite{ZhangL07}), meta-model-based approaches \cite{TabatabaeiHHMS15,ChughSHM19} have been generally used for computationally expensive multi-objective optimization.
%




Some mathematical derivative-free optimizers (e.g., NEWUOA \cite{Powell08}, BOBYQA \cite{Powell09}, and SLSQP \cite{Kraft88}) have shown their effectiveness for computationally expensive \textit{single-objective} black-box numerical optimization.
For example, Hansen et al. \cite{HansenARFP10} investigated the performance of 31 optimizers on the noiseless BBOB function set \cite{HansenFRA09bbob}.
Their results showed that NEWUOA achieves the best performance in the 31 optimizers for a small number of function evaluations.
The results in \cite{PosikH12,RiosS13} also reported the excellent convergence performance of NEWUOA.
The results in \cite{Hansen19} demonstrated that SLSQP can quickly find the optimal solution on some unimodal functions.
In \cite{BajerPRH19}, Bajer et al. showed that BOBYQA outperforms some meta-model-based optimizers including SMAC \cite{HutterHL11} and lmm-CMA \cite{BouzarkounaAD11}.
\textit{Motivation.}
Let $g_{\vector{w}}: \mathbb{R}^M \rightarrow \mathbb{R}$ be a scalarizing function that maps an $M$-dimensional objective vector to a scalar value.
Let also $\vector{W} = \{\vector{w}_k\}^K_{k=1}$ be a set of $K$ uniformly distributed weight vectors.
Under certain conditions, the optimal solution of a single-objective scalar optimization problem can be a weakly Pareto optimal solution (see Chapter 3.5 in \cite{Miettinen98}).
Therefore, $K$ weakly Pareto optimal solutions can potentially be obtained by solving a sequence of $K$ single-objective scalar optimization problems $\{g_{\vector{w}_k}\}_{k=1}^{K}$.
Any single-objective optimizer can be applied to the $K$ scalar optimization problems in principle.
When the number of function evaluations is limited, a mathematical derivative-free optimizer is likely to be suitable for this purpose based on the above review.



Actually, the first warm start phase in HMO-CMA-ES \cite{LoshchilovG16} adopts this idea.
HMO-CMA-ES was designed to achieve good anytime performance for bi-objective optimization in terms of the hypervolume indicator \cite{ZitzlerT98}.
HMO-CMA-ES is a hybrid multi-objective optimizer that consists of four phases.
The first out of the four phases in HMO-CMA-ES applies BOBYQA to a sequence of $K$ scalar optimization problems $\{g_{\vector{w}_k}\}_{k=1}^{K}$ for only the first $10 \times N$ function evaluations, where $N$ is the number of variables.
Let $\vector{X}$ be the set of all solutions found so far by BOBYQA.
At the end of the first phase, HMO-CMA-ES selects five solutions from $\vector{X}$ by applying environmental selection in SMS-EMOA \cite{BeumeNE07}.
Then, the second phase in HMO-CMA-ES performs a steady-state MO-CMA-ES \cite{IgelSH06} with the initial population of the five solutions.
%
Brockhoff et al. \cite{BrockhoffPAH21} showed that HMO-CMA-ES performs significantly better than some multi-objective optimizers for the first $10 \times N$ function evaluations, including NSGA-II \cite{DebAPM02}, COMO-CMA-ES \cite{ToureHAB19}, and DMS \cite{CustodioMVV11}.
Thus, their results indicate the effectiveness of mathematical derivative-free approaches to solving a scalar problem for computationally expensive multi-objective optimization.




One drawback of the above-discussed scalar optimization approach is that it can achieve only $K$ solutions that are sparsely distributed in the objective space, even in the best case.
Since only a limited number of function evaluations are available for computationally expensive optimization, $K$ needs to be as small as possible.
Due to the small value of $K$, the above-discussed scalar optimization approach cannot obtain a set of non-dominated solutions that cover the entire Pareto front in the objective space.

However, we believe that the issue of the above-discussed scalar optimization approach can be addressed by using a solution interpolation method.
Let $\vector{X}=\{\vector{x}_k\}_{k=1}^{K}$ be a set of $K$ solutions obtained by optimizing a sequence of $K$ single-objective scalar optimization problems $\{g_{\vector{w}_k}\}_{k=1}^{K}$.
Densely distributed solutions in the objective space can potentially be obtained by interpolating the $K$ sparsely distributed solutions in $\vector{X}$.
Some solution interpolation methods have been proposed in the literature (see Section \ref{sec:related_work}).
Unfortunately, existing methods were not designed for interpolating only a few (say $K \in \mathcal{O}(M)$) solutions.
In addition, we are particularly interested in optimization with a small budget of function evaluations.




The B\'{e}zier simplex is an extended version of the B\'{e}zier curve \cite{Farin02} to higher dimensions.
For a certain class of problems, the B\'{e}zier simplex has a capability to interpolate $K$ solutions, approximating the entire set of Pareto optimal solutions.
More precisely, Hamada et al. \cite{HamadaHIKT2020} showed that the set of Pareto optimal solutions is homeomorphic to an $(M-1)$-dimensional simplex under certain conditions. 
In such a case, Kobayashi et al. \cite{KobayashiHSTBS19} proved that a B\'{e}zier simplex model can approximate the Pareto optimal solution set.
They also proposed an algorithm for fitting a B\'{e}zier simplex by extending the B\'ezier curve fitting \cite{BorgesP2002}. 
Their results in \cite{KobayashiHSTBS19} demonstrated that it achieved an accurate approximation with a small number of solutions.
Thus, we expect that the B\'{e}zier simplex model can effectively interpolate the $K$ sparsely distributed solutions.


%

\textit{Contribution.}
Motivated by the above discussion, this paper proposes a two-phase framework with a B\'{e}zier simplex-based interpolation method (TPB) for computationally expensive multi-objective black-box optimization.
%
The first phase performs a mathematical derivative-free optimizer on a sequence of $K$ single-objective scalar optimization problems $\{g_{\vector{w}_k}\}_{k=1}^{K}$.
%
The second phase fits a B\'{e}zier simplex model to the $K$ solutions obtained in the first phase.
Then, TPB samples interpolated solutions from the B\'{e}zier simplex model.
We investigate the performance of TPB on the bi-objective BBOB function set \cite{BrockhoffAHT22}.
We also compare TPB with HMO-CMA-ES and state-of-the-art meta-model-based multi-objective optimizers.




\textit{Outline.}
%
Section \ref{sec:preliminaries} provides some preliminaries.
Section \ref{sec:related_work} reviews related work.
Section \ref{sec:proposed_method} introduces TPB.
Section \ref{sec:setting} describes our experimental setting.
Section \ref{sec:results} shows analysis results.
Section \ref{sec:conclusion} concludes this paper.



\textit{Code availability.}
The code of TPB is available at \url{https://github.com/ryojitanabe/tpb}.



\section{Preliminaries}
\label{sec:preliminaries}

\subsection{Multi-objective optimization}
\label{sec:def_MOPs}

We tackle a multi-objective minimization of a vector-valued objective function $\vector{f}: \mathbb{X} \to \mathbb{R}^{M}$, where $\mathbb{X} \subseteq \mathbb{R}^{N}$ is the search space.
Note that $M$ is the dimension of the objective space, and $N$ is the dimension of the search space.
Let $\vector{f} = (f_1, \dots, f_M)$, where $f_m : \mathbb{X} \to \mathbb{R}$ is called the $m$-th objective function.
The image of $\vector{f}$, $\mathbb{R}^M$ in our case, is called the objective space. 
Throughout of this paper, we consider a box constrained search space, i.e., $\mathbb{X} = [\texttt{LB}_1, \texttt{UB}_1] \times \cdots \times [\texttt{LB}_N, \texttt{UB}_N]$, where $\texttt{LB}_n$ and $\texttt{UB}_n$ are the lower and upper bounds of the $n$-th coordinate of the search space.

Our objective is to find a finite set of solutions $\vector{B}$ that approximates the Pareto front $\mathcal{P}(\vector{f})$, which is defined as follows:
\begin{equation}
\mathcal{P}(\vector{f}) = \{\vector{f}(\vector{x}) \mid \vector{x} \in \mathbb{X},\ \not\exists \vector{y}\in\mathbb{X}\ \text{s.t.}\ \vector{f}(\vector{y}) \prec \vector{f}(\vector{x})\},
\end{equation}
where $\vector{a} \prec \vector{b}$ (for $\vector{a}, \vector{b} \in \mathbb{R}^{M}$) represents the Pareto dominance relation ($1$ if $a_i \leq b_i$ holds for all $i = 1,\dots,M$ and $a_j < b_j$ holds for some $j$, and $0$ otherwise).
A solution $\vector{x}^\star$ is said to a Pareto optimal solution if no solution in $\mathbb{X}$ can dominate $\vector{x}^\star$.
The Pareto optimal solution set $\vector{X}^\star$ is the set of all $\vector{x}^\star$.
The objective is informally stated as to find a set of approximate Pareto optimal solutions that are well-distributed on $\mathcal{P}(\vector{f})$. The quality of $\vector{B}$ is often measured by a quality indicator  such as the hypervolume indicator \cite{ZitzlerT98}. 

In this paper, we suppose that we can access the objective function only through an expensive  black-box query $\vector{f}: \vector{x} \mapsto \vector{f}(\vector{x})$. Its indication is summarized below.
(1) The Jacobian and higher order information of $\vector{f}$ is unavailable (\emph{derivative-free optimization}). 
(2) The characteristic constants of $\vector{f}$ such as the Lipschitz constant are unavailable (\emph{black-box optimization}). 
(3) Evaluation of $\vector{f}(\vector{x})$ is computationally expensive (\emph{expensive optimization}). 
(4) Each objective function value $f_m(\vector{x})$ cannot be obtained with a lower computational cost. 
Therefore, the cost of the optimization process is measured by the number of $\vector{f}$-calls. 
We assume that it is limited up to $20 \times N, \dots, 40 \times N$.

\subsection{Simplicial problem}
\label{sec:simplicial_prob}

Kobayashi \textit{et al.}~\cite{KobayashiHSTBS19} defined a class of multi-objective optimization problems whose Pareto optimal solution set and Pareto front can be seen topologically as a simplex.
Let $M$ be a positive integer. 
The \emph{standard ($M-1$)-simplex} is denoted by
\[
    \Delta^{M-1} 
    = 
    \Set{
        (t_1, \dots, t_M) \in \mathbb{R}^M 
        | 
        {\displaystyle \sum_{m=1}^M t_m = 1,~t_m\geq 0}
    }
    .
\]
Let $I:= \{1,\dots,M\}$ be the index set on the objective functions.
For each non-empty subset $J\subseteq I$, we define 
\begin{equation*}
    \Delta^{J}:= \{(t_1,\dots,t_M)\in \Delta^{M-1}\mid t_m= 0~(m\notin J)\}
\end{equation*}
and 
\begin{equation*}
 \vector{f}_J:= (f_j)_{j\in J}: \mathbb{X}\to \mathbb{R}^{|J|}.
\end{equation*}

\begin{definition}
    For a given objective function $\vector{f}: \mathbb{X} \to \mathbb{R}^{M}$, the multi-objective optimization problem of minimizing $\vector{f}$ is \textit{simplicial} if there exists a map $\vector{\phi}:\Delta^{M-1}\to \mathbb{X}$ such that for each non-empty subset $J\subseteq I$, its restriction $\left.\vector{\phi}\right|_{\Delta^J}:\Delta^J\to \mathbb{X}$ gives the following homeomorphisms:
\begin{align}
    \left.\vector{\phi}\right|_{\Delta^J} &: \Delta^J\to \vector{X}^\star (\vector{f}_J),\\
    \left.\vector{f}\circ\vector{\phi}\right|_{\Delta^J} &: \Delta^J\to \mathcal{P}(\vector{f}_J).
\end{align}
\end{definition}

\subsection{B\'{e}zier simplex fitting}
\label{sec:bez}

We denote the set of non-negative integers (including zero) by $\mathbb{N}$.
Let $D$ be an arbitrary integer in $\mathbb{N}$, and 
\[
\mathbb{N}_D^M := \Set{(d_1,\dots, d_M) \in \mathbb{N}^M | {\displaystyle \sum_{m = 1}^M d_m = D}}.
\]
An \emph{$(M - 1)$-B\'ezier simplex of degree $D$} is a mapping $ \vector{b}: \Delta^{M - 1} \to \mathbb{R}^N$ determined by \emph{control points} $\vector{p}_{\vector{d}} \in \mathbb{R}^N$ $(\vector{d} \in \mathbb{N}_D^M)$ as follows:
\begin{equation}\label{eqn:bezier-simplex}
    \vector{b}(\vector{t}) := \sum_{\vector{d} \in \mathbb{N}_D^M} \binom{D}{\vector{d}} \vector{t}^{\vector{d}} \vector{p}_{\vector{d}},
\end{equation}
where $\binom{D}{\vector{d}} := \frac{D!}{d_1! \dots d_M!} (\in\mathbb{R})$ is a multinomial coefficient, and $\vector{t}^{\vector{d}} := t^{d_1}_1 \dots t^{d_M}_M (\in\mathbb{R})$ is a monomial for each $\vector{t} := (t_1, \dots, t_M) \in \Delta^{M - 1}$ and $\vector{d} := (d_1, \dots, d_M) \in \mathbb{N}^M_D$. 
The following theorem ensures that the Pareto optimal solution set and Pareto front of any simplicial problem can be approximated with arbitrary accuracy by a B\'ezier simplex of an appropriate degree:
\begin{theorem}[{Kobayashi \textit{et al.} \cite[Theorem~1]{KobayashiHSTBS19}}]
    Let $\vector{\phi}: \Delta^{M-1}\to\mathbb{R}^N$ be a continuous map.
    There is an infinite sequence of B\'ezier simplices $\vector{b}^{(i)}: \Delta^{M-1}\to\mathbb{R}^N$ such that
    \[
        \lim_{i \to \infty}\sup_{\vector{t}\in\Delta^{M-1}} \left\|\vector{\phi}(\vector{t}) - \vector{b}^{(i)}(\vector{t})\right\| = 0.
    \]
\end{theorem}

With this result, Kobayashi \textit{et al.} \cite{KobayashiHSTBS19} proposed the B\'ezier simplex fitting method to describe the Pareto optimal solution set of a simplicial problem. 
Suppose that we have a set of approximate Pareto optimal solutions $\{(\vector{x}_l, \vector{t}_l) \in \mathbb{X}\times \Delta^{M - 1}~\mid l = 1, \dots, L\}$, where $\vector{x}_l$ and $\vector{t}_l$ are the $l$-th approximate Pareto optimal solution and its corresponding parameter, respectively.  
The B\'ezier simplex fitting method adjusts the control points by minimizing the ordinary least squares (OLS) loss function: $\frac{1}{L}\sum_{l=1}^{L} \|\vector{x}_l - \vector{b}(\vector{t}_l)\|^2$. 
Since the OLS loss function is a convex quadratic function with respect to $\vector{p}_{\vector{d}}$, its minimization problem can be solved efficiently, for example, by solving a normal equation.


\section{Related work}
\label{sec:related_work}

%
Two-phase approaches have been well studied in the context of multi-objective optimization (e.g., \cite{HamadaSKO08,HiranoY13,HuYL17,Regis21}).
TPLS$+$PLS \cite{PaqueteS03,Dubois-LacosteLS13} is one of the most representative two-phase approaches for combinatorial optimization.
Roughly speaking, the first phase in multi-objective two-phase approaches aims to find well-converged solutions to the Pareto front.
Then, the second phase aims to generate a set of well-diversified solutions based on the solutions obtained in the first phase.
Generally, two-phase approaches can produce only a poor-quality solution set when it stops before the maximum budget of function evaluations \cite{Dubois-LacosteLS11amai}.
Thus, the anytime performance of most two-phase approaches is poor.
Here, we say that the anytime performance of an optimizer is good if it can obtain a well-approximated solution set at any time during the search process.
%
%
The substantial difference between TPB and existing two-phase approaches is that the second phase in TPB incorporates solutions by utilizing a B\'{e}zier simplex model, which fully exploits the theoretical property of the Pareto optimal solution set.
In addition, unlike TPB, all two-phase approaches but \cite{Regis21} were designed for non-expensive optimization.
Here, the study \cite{Regis21} proposed a surrogate model-based approach for constrained bi-objective optimization.



Some methods for interpolating objective vectors (not solutions) obtained by an EMO algorithm have been proposed in the literature \cite{HartikainenMW11,HartikainenMW12,BhattacharjeeSR17J}.
A decision-maker can determine her/his preference by visually examining interpolated objective vectors.
One of the most representative approaches is the PAINT method \cite{HartikainenMW12}, which interpolates an objective vector set using the Delaunay triangulation.
Note that these interpolation methods cannot provide an inverse mapping from the objective space to the search space.
In contrast, the second phase in TPB aims to interpolate solutions (not objective vectors) to approximate the Pareto front.

The Pareto estimation method \cite{GiagkiozisF14} aims to increase the number of non-dominated solutions obtained by an EMO algorithm.
The Pareto estimation method uses a neural network model to find an inverse mapping from the objective space to the search space.
GAN-LMEF \cite{WangHYZJT22} interpolates randomly generated solutions on the manifold by using dimensionality reduction, clustering, and GAN \cite{GoodfellowPMXWOCB14}.
These two methods aim to interpolate a sufficiently large number of solutions. 
In contrast, the second phase in TPB aims to interpolate only $K$ solutions (i.e., $K = M+1 = 2+1 =3$ in this study) by utilizing a B\'{e}zier simplex model.

Some EMO algorithms (e.g., RM-MEDA \cite{ZhangZJ08}) exploit the simplex structure of the Pareto optimal solution set.
BezEA \cite{MareeAB20} evolves a control point set for a B{\'{e}}zier curve to generate a high-quality solution set in terms of the ``smoothness'' measure, which was proposed in \cite{MareeAB20}.
Unlike these EMO algorithms, TPB exploits the property of the Pareto optimal solution set by using the theoretically well-founded B\'{e}zier simplex.
No previous study also proposed an EMO algorithm based on the simplex structure of the Pareto optimal solution set for computationally expensive optimization.


\section{Proposed framework}
\label{sec:proposed_method}

This section describes the proposed TPB, which consists of the first phase (Section \ref{sec:opt_phase}) and the second phase (Section \ref{sec:int_phase}).
Let $\vector{W} = \{\vector{w}_k\}^K_{k=1}$ be a set of $K$ weight vectors.
We assume that $K=M+1$, which is the minimum value of $K$.

In the first phase (Section \ref{sec:opt_phase}), TPB aims to approximate $K$ Pareto optimal solutions by applying a single-objective optimizer to $K$ scalar optimization problems $\{g_{\vector{w}_k}\}_{k=1}^{K}$.
Let $\vector{B}^*$ be a set of the best $K$ solutions for the $K$ scalar problems obtained in the first phase.
Here, the $k$-th solution in $\vector{B}^*$ should correspond to the $k$-th weight vector in $\vector{W}$.
Ideally, the first phase should find $\vector{B}^*$ such that $\vector{x}_k$ in $\vector{B}^*$ minimizes its corresponding scalar problem $g_{\vector{w}_k}$.
Let \texttt{budget} be the maximum budget of function evaluations for the whole process of TPB.
The first phase in TPB can use  $\lfloor $\texttt{budget} $ \times r^{\mathrm{1st}} \rfloor$ function evaluations in the maximum case, where $r^{\mathrm{1st}} \in [0, 1]$ is a control parameter of TPB.
For example, when \texttt{budget}$=40$ and $r^{\mathrm{1st}}=0.9$, $36$ function evaluations can be used in the first phase in the maximum case.
Note that some optimizers have their own stopping criteria in addition to the maximum number of function evaluations.
For example, BOBYQA stops when reaching its minimum trust region radius.
Thus, it is possible that the first phase in TPB does not use all $\lfloor$\texttt{budget} $ \times r^{\mathrm{1st}} \rfloor$ function evaluations.

The second phase in TPB (Section \ref{sec:int_phase}) aims to interpolate the $K$ solutions in $\vector{B}^*$ by using a B\'{e}zier simplex-based interpolation method \cite{KobayashiHSTBS19}.
The B\'{e}zier simplex model can approximate the Pareto optimal solution set (see Section  \ref{sec:bez}).
In addition, the B\'{e}zier simplex-based interpolation can be done by minimizing the OLS function, which is a convex quadratic function.

Below, Sections \ref{sec:opt_phase} and \ref{sec:int_phase} describe the first and second phases in TPB, respectively.
Section \ref{sec:tpb_discussion} discusses the property of TPB. 

\subsection{First phase}
\label{sec:opt_phase}


Algorithm \ref{alg:opt_phase} shows the first phase in TPB.
In line 1 in Algorithm \ref{alg:opt_phase}, \texttt{budget}$^{\mathrm{opt}}$ is the maximum budget of function evaluations used in an optimizer on each scalar problem.
%
In line 2 in Algorithm \ref{alg:opt_phase}, $\vector{X}$ is an archive that maintains all solutions found so far.


As in D-TPLS \cite{PaqueteS03}, the first phase in TPB first performs single-objective optimization of each objective function $f \in \{f_1, \dots, f_M\}$ (lines 3--6 in Algorithm \ref{alg:opt_phase}).
This aims to approximate $M$ Pareto optimal solutions that minimize the $M$ objective functions, respectively.
Unlike D-TPLS, the $M$ solutions are mainly used for the normalization procedure in the next step (lines 7--15 in Algorithm \ref{alg:opt_phase}).
TPB sets the initial solution $\vector{x}_{\mathrm{init}}$ to the center of the search space $\vector{x}_{\mathrm{center}}$ (line 3 in Algorithm \ref{alg:opt_phase}), where the $n$-th element in $\vector{x}_{\mathrm{center}}$ is $(\texttt{LB}_n + \texttt{UB}_n) / 2$.
Then, TPB applies a pre-defined single-objective optimizer (\texttt{optimizer}) to each objective function (line 5 in Algorithm \ref{alg:opt_phase}).
Here, $\vector{Y}$ is a set of all solutions found by \texttt{optimizer}.


Next, the first phase in TPB aims to solve the remaining $K-M$ scalar problem(s).
Since TPB has solved the $M$ objective functions, TPB here does not consider the $M$ extreme weight vectors $\{\vector{e}_m\}_{m=1}^{M}$ (line 7 in Algorithm \ref{alg:opt_phase}).
TPB sets the approximated ideal point $\vector{z}^{\mathrm{ideal}} \in \mathbb{R}^M$ and the approximated nadir point $\vector{z}^{\mathrm{nadir}} \in \mathbb{R}^M$ based on $\vector{X}$ (line 8 in in Algorithm \ref{alg:opt_phase}).
Note that this step always normalizes the objective vector $\vector{f}(\vector{x})$ as follows: $\vector{f}(\vector{x}) = (\vector{f}(\vector{x}) - \vector{z}^{\mathrm{ideal}}) / (\vector{z}^{\mathrm{nadir}} - \vector{z}^{\mathrm{ideal}})$.
The initial solution $\vector{x}_{\mathrm{init}}$ is set to the best solution in $\vector{X}$ in terms of a given scalarizing function (line 9 in Algorithm \ref{alg:opt_phase}).

Finally, we set $\vector{B}^* = \{\vector{x}_{\mathrm{best}, k}\}^K_{k=1}$, where $\vector{x}_{\mathrm{best}, k}$ is the best-so-far solution of the $k$-th scalar problem (lines 12--15 in Algorithm \ref{alg:opt_phase}).
The second phase in TPB interpolates the $K$ solutions in $\vector{B}^*$.

\IncMargin{0.5em}
\begin{algorithm}[t]
\SetSideCommentRight
\SetKwInOut{Input}{input}
\SetKwInOut{Output}{output}
%
%
\texttt{budget}$^{\mathrm{opt}}$ $\leftarrow \lfloor ($\texttt{budget} $ \times r^{\mathrm{1st}}) / K \rfloor$\;
$\vector{X} \leftarrow \emptyset$\; 
$\vector{x}_{\mathrm{init}} \leftarrow \vector{x}_{\mathrm{center}}$\;
\For{$f \in \{f_1, \dots, f_M\}$}{
    $\vector{Y}$ $\leftarrow$ Run \texttt{optimizer} on $\vector{f}$ with $\vector{x}_{\mathrm{init}}$ and \texttt{budget}$^{\mathrm{opt}}$\;
    $\vector{X} \leftarrow \vector{X} \cup \vector{Y}$\; 
}
%
\For{$\vector{w} \in \vector{W} \setminus \{\vector{e}_m\}_{m=1}^{M}$}{
    Set $\vector{z}^{\mathrm{ideal}}$ and  $\vector{z}^{\mathrm{nadir}}$ based on $\vector{X}$\;
  $\vector{x}_{\mathrm{init}} \leftarrow \argmin_{\vector{x} \in \vector{X}} \{g_{\vector{w}}(\vector{x} |  \vector{z}^{\mathrm{ideal}}, \vector{z}^{\mathrm{nadir}})\}$\;  
    $\vector{Y}$ $\leftarrow$ Run \texttt{optimizer} on $g_{\vector{w}}$ with $\vector{x}_{\mathrm{init}}$ and \texttt{budget}$^{\mathrm{opt}}$\;
    $\vector{X} \leftarrow \vector{X} \cup \vector{Y}$\; 
}
%
$\vector{B}^* \leftarrow \emptyset$\; 
\For{$\vector{w} \in \vector{W}$}{
  $\vector{x}_{\mathrm{best}} \leftarrow \argmin_{\vector{x} \in \vector{X}} \{g_{\vector{w}}(\vector{x} | \vector{z}^{\mathrm{ideal}}, \vector{z}^{\mathrm{nadir}})\}$\;
  $\vector{B}^* \leftarrow \vector{B}^* \cup \{\vector{x}_{\mathrm{best}}\}$;
}
\caption{The first phase in TPB}
\label{alg:opt_phase}
\end{algorithm}
\DecMargin{0.5em}

\subsection{Second phase}
\label{sec:int_phase}

Let \texttt{budget}$^{\mathrm{1st}}$ be the number of function evaluations used in the first phase, where the maximum \texttt{budget}$^{\mathrm{1st}}$ is $\lfloor $\texttt{budget} $ \times $ $r^{\mathrm{1st}} \rfloor$.
The second phase in TPB uses the remaining \texttt{budget}$^{\mathrm{2nd}}$ $ = $ \texttt{budget} $-$ \texttt{budget}$^{\mathrm{1st}}$ function evaluations. 


Let $\vector{T}^{\mathrm{fit}}=\{\vector{t}^\mathrm{fit}_k\}^K_{k=1}$ be a set of $K$ parameter vectors, where $\vector{t}^\mathrm{fit}_k \in \Delta^{M-1}$.
TPB treats the $k$-th weight vector $\vector{w}_k$ in $\vector{W}$ as the $k$-th parameter $\vector{t}^\mathrm{fit}_k$ in $\vector{T}^{\mathrm{fit}}$.
Thus, $\vector{T}^{\mathrm{fit}}$ is identical to $\vector{W}$. 
With  $\vector{B}^*$ and $\vector{T}^{\mathrm{fit}}$, we next train a B\'ezier simplex model $\vector{b}:\Delta^{M-1}\to \mathbb{R}^N$ that takes a parameter  $\vector{t} \in \Delta^{M-1}$ as an input and outputs a minimizer of the corresponding scalarizing function.
Specifically, TPB fits a B\'{e}zier simplex model $\vector{b}$ to $\vector{B}^*$ with $\vector{T}^{\mathrm{fit}}$ by solving the OLS loss minimization problem:
\begin{equation}
\label{eqn:fit}
    \text{minimize}\quad \sum_{k=1}^{K}\|\vector{x}_k-\vector{b}(\vector{t}^{\mathrm{fit}}_k)\|^2,
\end{equation}
where $\vector{x}_k$ is the $k$-th solution in $\vector{B}^*$. 

Let $\vector{T}^{\mathrm{int}}$ be a set of \texttt{budget}$^{\mathrm{2nd}}$ parameter vectors.
After fitting the B\'{e}zier simplex model $\vector{b}$ in \eqref{eqn:fit},
TPB generates \texttt{budget}$^{\mathrm{2nd}}$ solutions by using $\vector{b}$ and $\vector{T}^{\mathrm{int}}$.
It is expected that the \texttt{budget}$^{\mathrm{2nd}}$ solutions complement the $K$ solutions in $\vector{B}^*$.
Any method can be used to generate $\vector{T}^{\mathrm{int}}$, e.g., uniform random generation.
The decision maker's preference can also be incorporated into $\vector{T}^{\mathrm{int}}$. 
In this study, we generate \texttt{budget}$^{\mathrm{2nd}}$ parameters in $\vector{T}^{\mathrm{int}}$ so that they are equally spaced.
First, we equally generate $($\texttt{budget}$^{\mathrm{2nd}}$ $+M)$ parameters on $\Delta^{M-1}$.
Then, we removed the $M$ extreme parameters $(1, 0)$ and $(0, 1)$ from $\vector{T}^{\mathrm{int}}$.
Since the first phase has found the $M$ extreme solutions, we do not need to re-generate them.
For example, when \texttt{budget}$^{\mathrm{2nd}} = 4$, we can obtain the following parameters: $\vector{t}^{\mathrm{int}}_1 = (0.2, 0.8), \vector{t}^{\mathrm{int}}_2 = (0.4, 0.6), \vector{t}^{\mathrm{int}}_3 = (0.6, 0.4),$ and $\vector{t}^{\mathrm{int}}_4 = (0.8, 0.2)$.

\subsection{Discussion}
\label{sec:tpb_discussion}


\subsubsection{Control parameters for TPB}
\label{sec:param_tpb}

The numerical control parameters for TPB include the number of weight vectors $K$, the degree in a B\'{e}zier simplex model $D$, and the budget ratio $r^{\mathrm{1st}}$.
Clearly, the best setting of $K$ and $D$ depends on the shape of the Pareto optimal solution set.
We believe that $K$ must be more than or equal to $M+1$ so that a resulting B\'{e}zier simplex model can characterize the shape of the Pareto optimal solution set.
This is because a B\'{e}zier simplex model fitting needs at least one non-extreme solution to handle the nonlinear Pareto optimal solution set.
Similarly, $D$ must be more than or equal to $2$ to handle the nonlinearity of the Pareto optimal solution set.
%
The best setting of $r^{\mathrm{1st}}$ depends on the difficulty in solving $K$ scalar problems.
If $K$ scalar problems are easy, $r^{\mathrm{1st}}$ should be a small value.
Otherwise, the first phase in TPB can waste computational resources.
However, as described at the beginning of Section \ref{sec:proposed_method}, some modern optimizers (e.g., BOBYQA) automatically terminate the search.
Thus, we believe that $r^{\mathrm{1st}}$ can be set to a relatively high value (e.g., $r^{\mathrm{1st}}=0.9$).
%

The categorical control parameters for TPB include the scalarizing function $g$ and the single-objective optimizer \texttt{optimizer}.
Although TPB can use any $g$ (e.g., the weighted Tchebycheff function), we set $g$ to the weighted sum function $g^{\mathrm{ws}}_{\vector{w}} = \sum_{m=1}^M w_m f_m (\vector{x})$ in this study.
Since $g^{\mathrm{ws}}$ is the simplest scalarizing function, $g^{\mathrm{ws}}$ is a reasonable first choice.
A mathematical derivative-free optimizer is suitable for \texttt{optimizer} for the reason discussed in Section \ref{sec:introduction}.
We set \texttt{optimizer} to BOBYQA, which is a state-of-the-art mathematical derivative-free optimizer for box-constrained optimization.
The first phase in HMO-CMA-ES also adopts $g^{\mathrm{ws}}$ and BOBYQA.

\subsubsection{Advantages and disadvantages of TPB}
\label{sec:pros_cons_tpb}

One advantage of TPB is that it can use a state-of-the-art single-objective optimizer without any change.
In contrast to meta-model-based optimizers, TPB does not require computationally expensive operations if a single-objective optimizer is computationally cheap.
TPB can also exploit the structure of the Pareto solution set by using the theoretically well-understood B\'{e}zier simplex.

As described in Section \ref{sec:related_work}, the anytime performance of two-phase approaches is generally poor.
TPB has the same disadvantage.
%
The second phase in TPB cannot interpolate $K$ solutions when a given problem is not simplicial (see Section \ref{sec:simplicial_prob}).
This is because a B\'{e}zier simplex model can represent only a standard $(M-1)$-simplex.
Fortunately, for a lot of practical real-world problems, scatter plots of approximate Pareto optimal solutions imply those problems are simplicial (e.g., \cite{ShovalSSHRNDKA2012,MastroddiG2013,VrugtGBBS2003,TanabeI2020}).




\section{Experimental setup}
\label{sec:setting}


We investigated the performance of the proposed TPB using COCO \cite{HansenARMTB21}, which is the standard benchmarking platform in the GECCO community.
We used the 55 bi-objective BBOB problems ($f_1, \dots, f_{55}$) \cite{BrockhoffAHT22} provided by COCO.
The first and second objective functions in a bi-objective BBOB problem are selected from the 24 single-objective noiseless BBOB functions \cite{HansenFRA09bbob}.
Although the DTLZ \cite{DebTLZ05} and WFG \cite{HubandHBW06} problems are the most commonly-used test problems, many previous studies (e.g., \cite{BrockhoffTH15,IshibuchiSMN16,ChenIS20}) pointed out that they have some serious issues, including the regularity of the Pareto front and the existence of distance and position variables.
In contrast, the bi-objective BBOB problems address all these issues.
Each bi-objective BBOB problem consists of 15  instances in COCO.
A single run of a multi-objective optimizer was performed on each problem instance.
In other words, 15 runs were performed for each problem.
We set the number of variables $N$ to 2, 3, 5, 10, and 20.

We used an automatic performance indicator ($I_{\mathrm{COCO}}$) \cite{BrockhoffTTWHA16} provided by COCO.
COCO uses an unbounded external archive to maintain all non-dominated solutions found so far.
When there exists at least a single solution in the archive that dominates a reference point $(1, 1)^{\top}$ in the normalized objective space $[0,1]^M$, the performance of optimizers is measured by a referenced version of the hypervolume indicator \cite{ZitzlerT98} using the archive. 
Otherwise, the performance of optimizers is measured by the smallest distance to the region of interest, which is bounded by the nadir point.


We compare TPB with HMO-CMA-ES \cite{LoshchilovG16}, ParEGO \cite{Knowles06}, MOTPE \cite{OzakiTWO20}, K-RVEA \cite{ChughJMHS18}, KTA2 \cite{SongWHJ22}, and EDN-ARMOEA \cite{GuoWGJDC22}.
We demonstrate the effectiveness of the second phase in TPB by comparing with the warm start phase in HMO-CMA-ES, which is based on a sophisticated scalarizing approach.
We are also interested in the performance of TPB compared to state-of-the-art meta-model-based optimizers. 
We used the optuna \cite{AkibaSYOK19} implementation of MOTPE and the PlatEMO \cite{TianCZJ17} implementation of the surrogate-assisted EMO algorithms.
We used the results of HMO-CMA-ES provided by the COCO data archive (\url{https://numbbo.github.io/data-archive}).



We set the control parameters for TPB based on the discussion in Section \ref{sec:param_tpb}, i.e., $K=M+1=3$ and $D=2$.
We used BOBYQA and $g^{\mathrm{ws}}$ as \texttt{optimizer} and $g$, respectively.
Here, we evaluated the performance of TPB with $r^{\mathrm{1st}} \in \{0.7, 0.8, 0.9\}$ on the first BBOB problem $f_1$ with $N=2$ in our preliminary study.
We set $r^{\mathrm{1st}}$ to $0.9$ based on the rough hand-tuning results.
We used the Py-BOBYQA \cite{CartisFMR19} implementation of BOBYQA.
Unlike HMO-CMA-ES, we used the default parameter setting of BOBYQA.
For the B\'{e}zier simplex model fitting method, we used the code provided by the authors of \cite{KobayashiHSTBS19} (\url{https://gitlab.com/hmkz/pytorch-bsf}).
We used a workstation with an Intel(R) 48-Core Xeon Platinum 8260 (24-Core$\times$2) 2.4GHz and 384GB RAM using Ubuntu 18.04.


We set the maximum budget of function evaluations (\texttt{budget}) to $20\times N$, $30\times N$, and $40\times N$.
As discussed in Section~\ref{sec:pros_cons_tpb}, TPB is not an anytime algorithm.
Thus, the behavior of TPB depends on the termination condition, i.e., \texttt{budget} in this study.
The performance of some state-of-the-art surrogate-assisted EMO algorithms (e.g., K-RVEA and KTA2) also depends on \texttt{budget}.
This is because they are not anytime algorithms similar to TPB.
For example, K-RVEA has a temperature-like parameter that determines the magnitude of the penalty value.
Generally, the best parameter setting for EMO algorithms depends on \texttt{budget} \cite{DymondKH13,BezerraLS18}.
In addition, \texttt{budget} has not been standardized in the field of computationally expensive multi-objective optimization.
For the above-discussed reasons, we used the three budget settings.

\section{Results}
\label{sec:results}


This section describes our analysis results.
Through experiments, Sections \ref{sec:analysis_tpb}--\ref{sec:para_study}  address the following research questions.\\







 \noindent \textbf{$\bullet$ Section \ref{sec:analysis_tpb}} How does TPB interpolate solutions?
 
 \noindent \textbf{$\bullet$ Section \ref{sec:vs_sota}} Is TPB competitive with state-of-the-art optimizers?

\noindent \textbf{$\bullet$ Section \ref{sec:inv_tp}} How important is the two-phase mechanism in TPB?

\noindent \textbf{$\bullet$ Section \ref{sec:para_study}} How does the choice of $K$ and $r^{\mathrm{1st}}$ influence the performance of TPB?

\subsection{On the solution interpolation in TPB}
\label{sec:analysis_tpb}

\begin{figure}[t]
  \centering
\subfloat[Search space]{  
\includegraphics[width=0.172\textwidth]{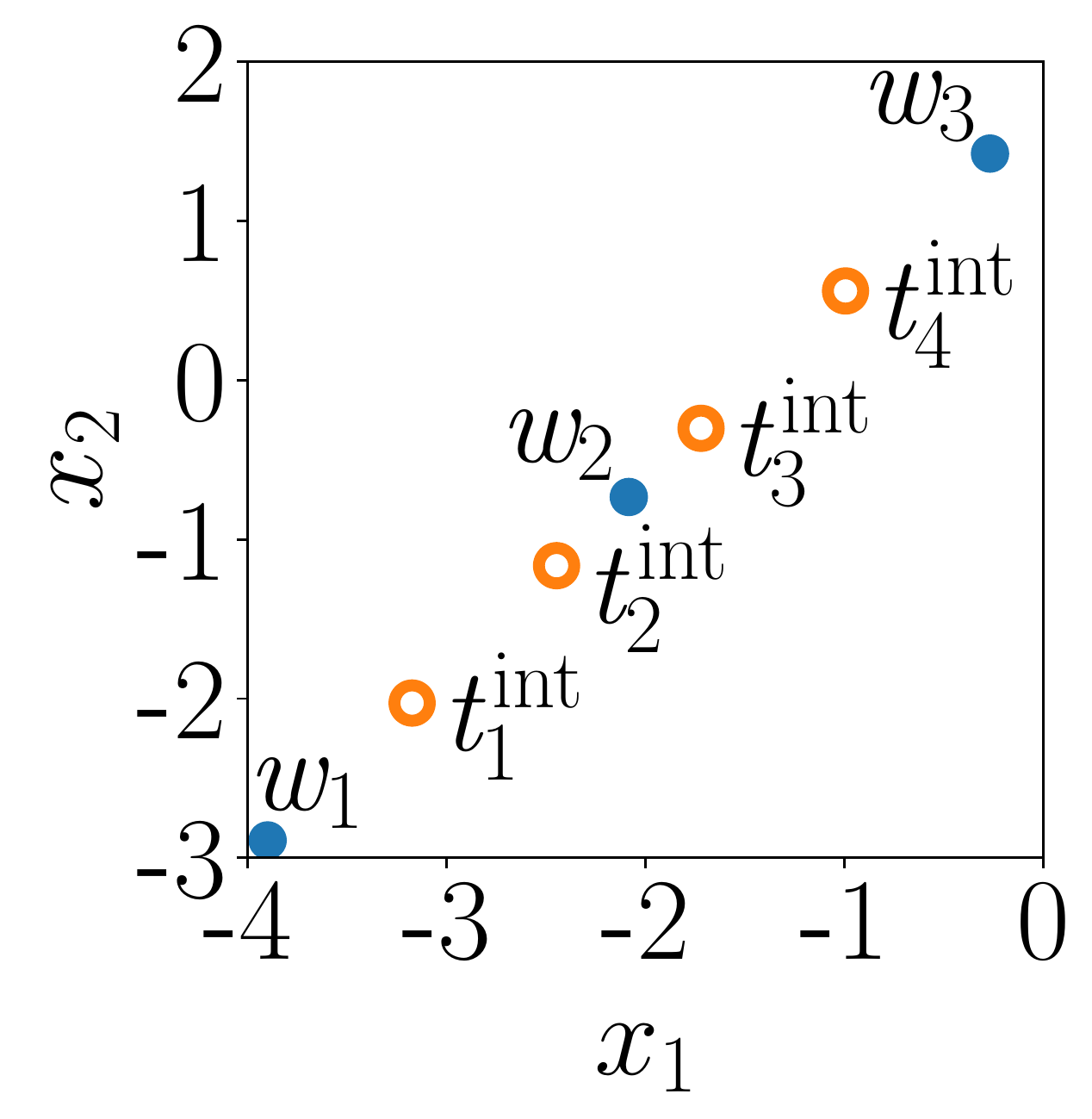}
}
\subfloat[Objective space]{
\includegraphics[width=0.2\textwidth]{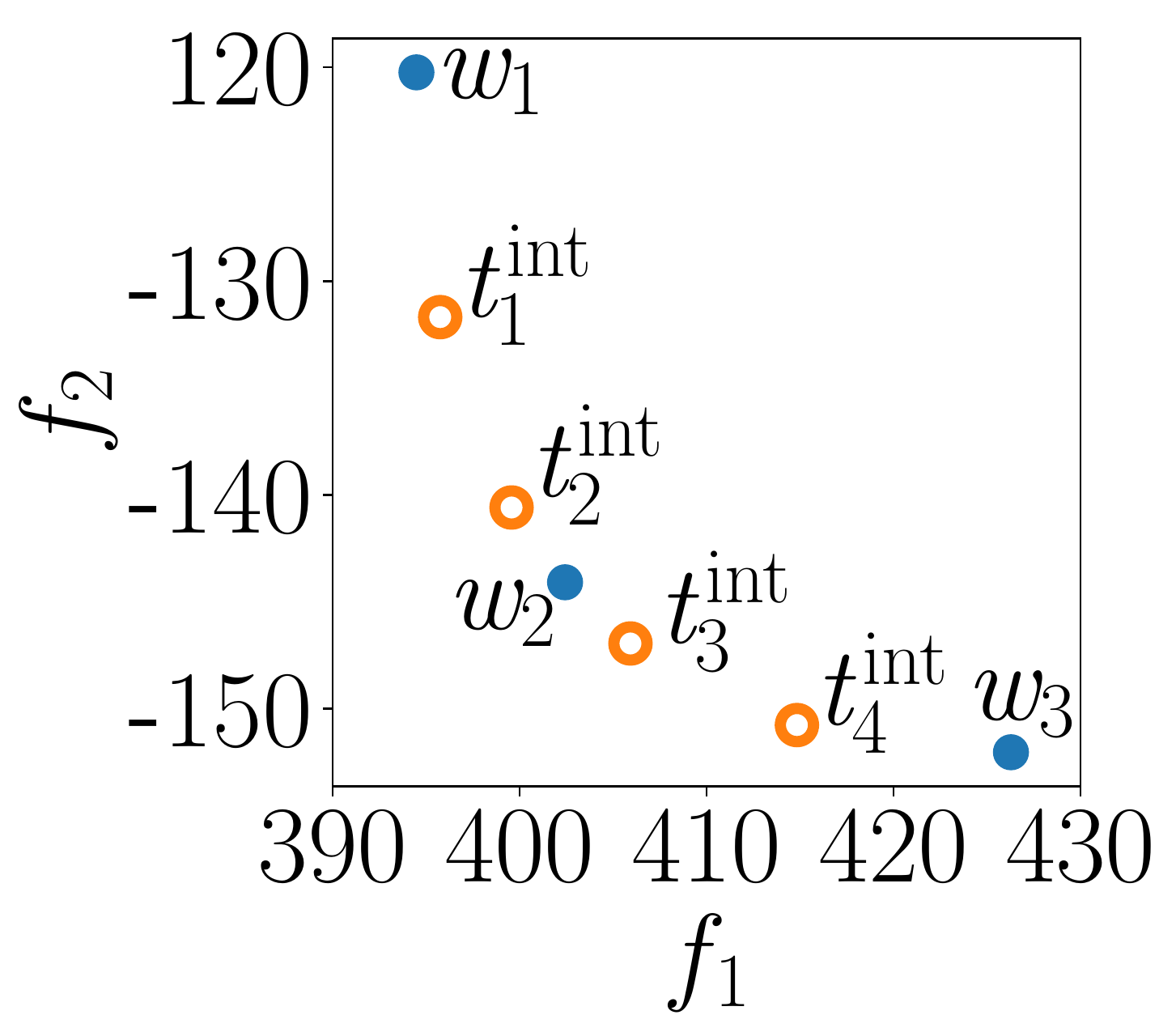}
}
\caption{Distribution of solutions obtained by TPB for $f_1$ for $N=2$ in the search and objective spaces.}
  \label{fig:dist_solutions}
\end{figure}

Figure \ref{fig:dist_solutions} shows the distribution of solutions generated by TPB for \texttt{budget} $=20 \times N$. 
Figure \ref{fig:dist_solutions} shows the results on the first instance of $f_1$ with $N=2$.
Note that the solution interpolation in TPB is performed in the search space (Figure \ref{fig:dist_solutions}(a)), not the objective space (Figure \ref{fig:dist_solutions}(b)).
%
In Figure \ref{fig:dist_solutions}, we confirmed that \texttt{budget}$^{\mathrm{1st}}$ $=0.9 \times 20 \times 2 = 36$ and \texttt{budget}$^{\mathrm{2nd}}$  $=40 - 36 = 4$.

In Figure \ref{fig:dist_solutions}, the three blue filled circles represent the three solutions  in $\vector{B}^*$ found in the first phase with the three weight vectors $\vector{w}_1 = (1, 0)$, $\vector{w}_2 = (0.5, 0.5)$, and $\vector{w}_3 = (0, 1)$.
In contrast, the four orange unfilled circles represent the four solutions generated in the second phase with the four parameters $\vector{t}^{\mathrm{int}}_1, \dots, \vector{t}^{\mathrm{int}}_4$ that are the same as in the example in Section \ref{sec:int_phase}.


As shown in Figure \ref{fig:dist_solutions}, the three solutions obtained in the first phase are well-converged to the Pareto front, but they are sparsely distributed.
The second phase makes up for this shortcoming.
As seen from Figure \ref{fig:dist_solutions}, the four solutions generated in the second phase incorporate the three solutions.
The four solutions are distributed as if they were obtained by a scalar optimization approach with the four weight vectors $(0.2, 0.8), (0.4, 0.6), (0.6, 0.4),$ and $(0.8, 0.2)$.
As a result, TPB could obtain the seven well-converged and well-distributed solutions.
As demonstrated here, the first and second phases in TPB are complementary to each other.

\subsection{Comparison with state-of-the-art optimizers}
\label{sec:vs_sota}

\begin{figure*}[t]
   \centering
\subfloat[$N=2$ and \texttt{budget} $=20 \times N$]{  
\includegraphics[width=0.24\textwidth]{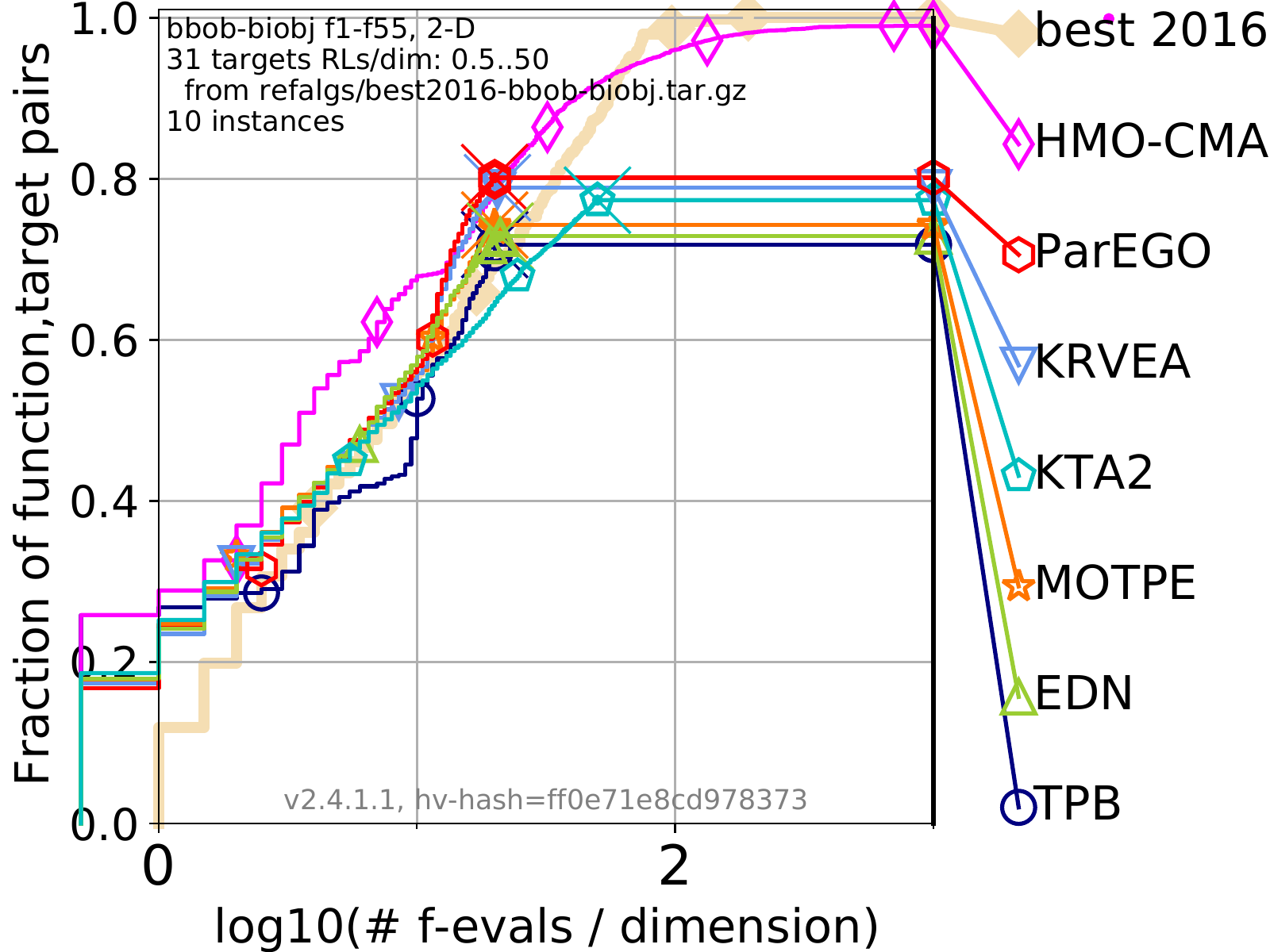}
}
\subfloat[$N=5$ and \texttt{budget} $=20 \times N$]{  
\includegraphics[width=0.24\textwidth]{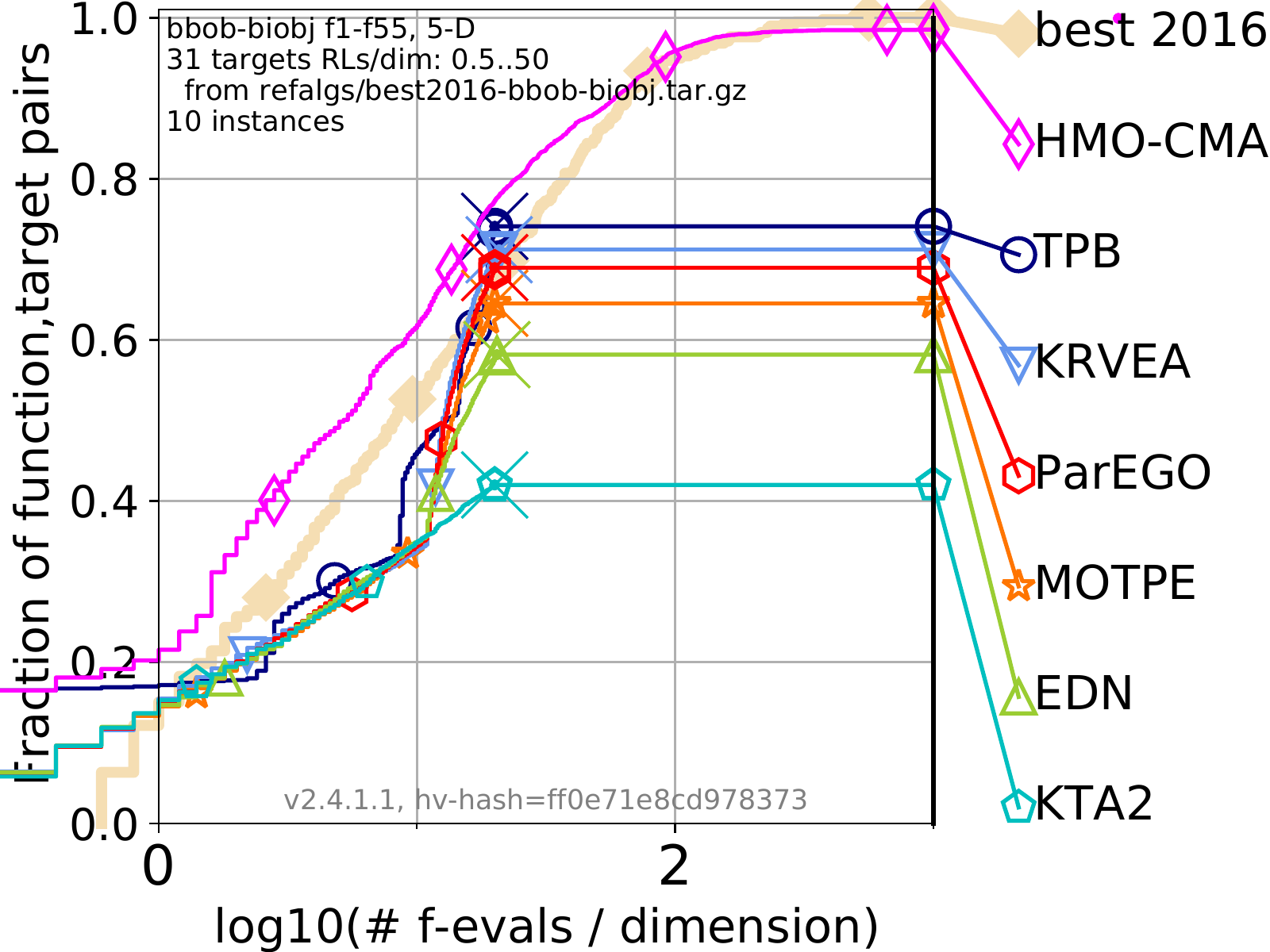}
}
\subfloat[$N=10$ and \texttt{budget} $=20 \times N$]{   
\includegraphics[width=0.24\textwidth]{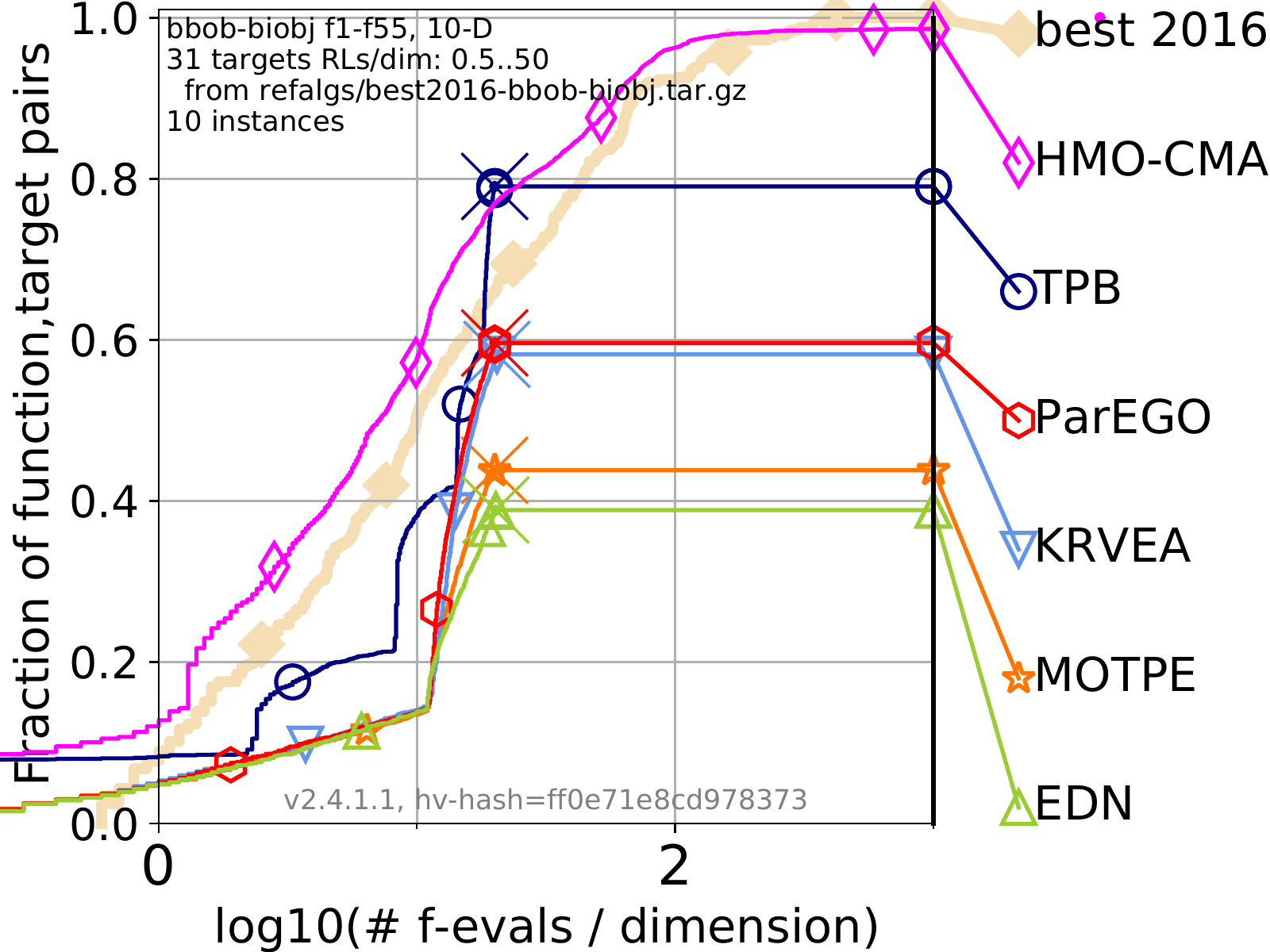}
}
\subfloat[$N=20$ and \texttt{budget} $=20 \times N$]{  
\includegraphics[width=0.24\textwidth]{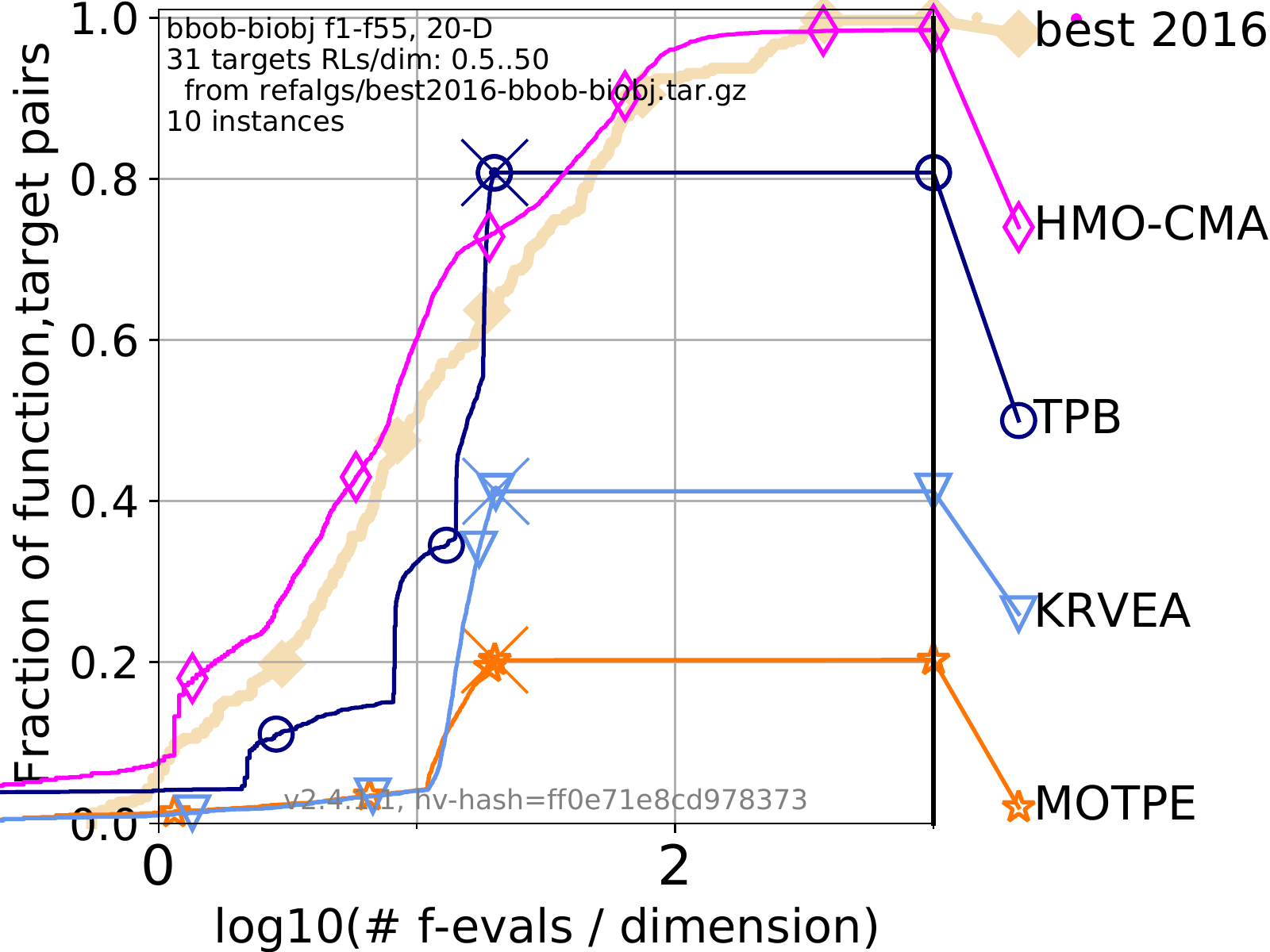}
}
\\
\subfloat[$N=2$ and \texttt{budget} $=30 \times N$]{  
\includegraphics[width=0.24\textwidth]{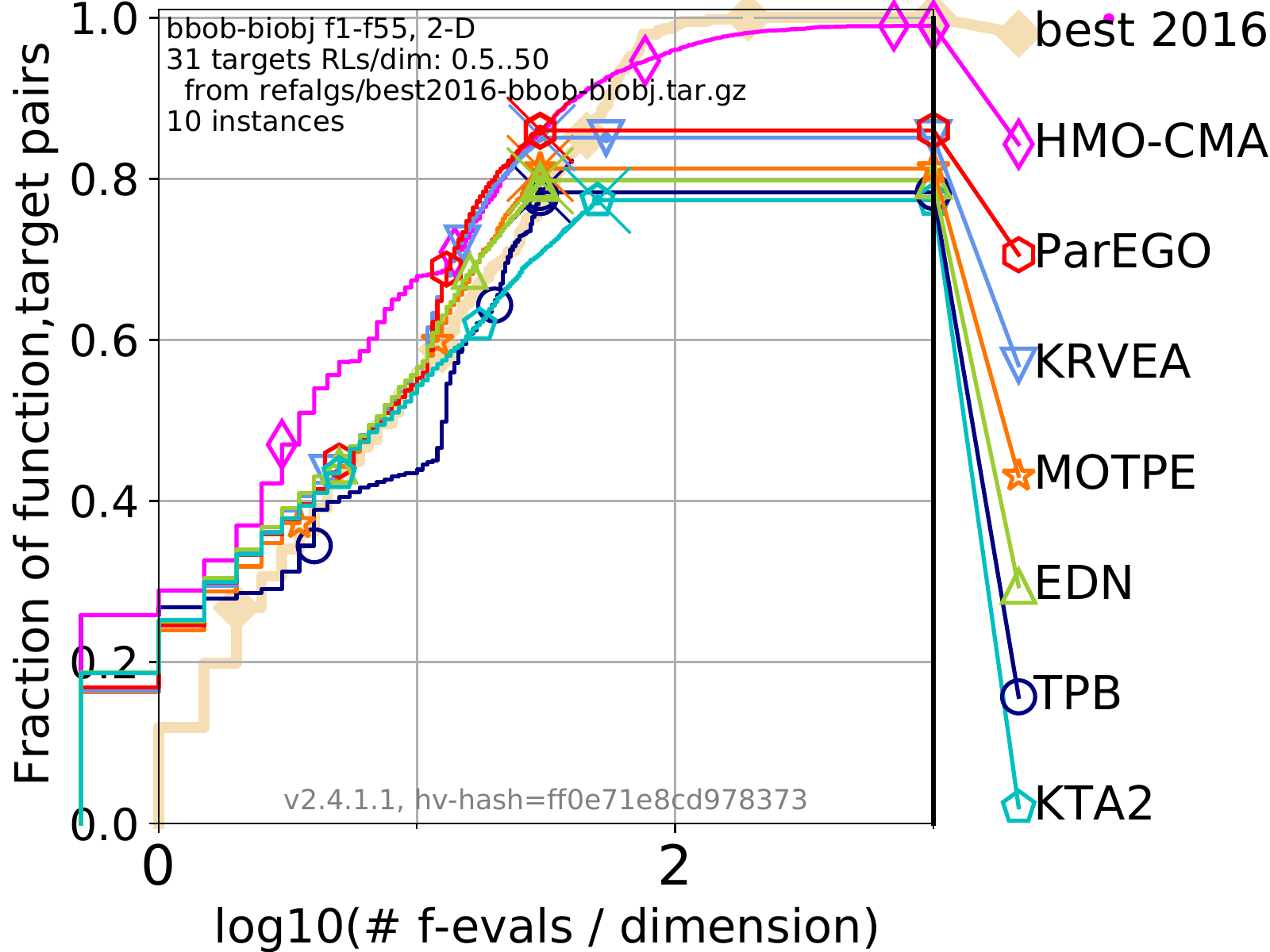}
}
\subfloat[$N=5$ and \texttt{budget} $=30 \times N$]{  
\includegraphics[width=0.24\textwidth]{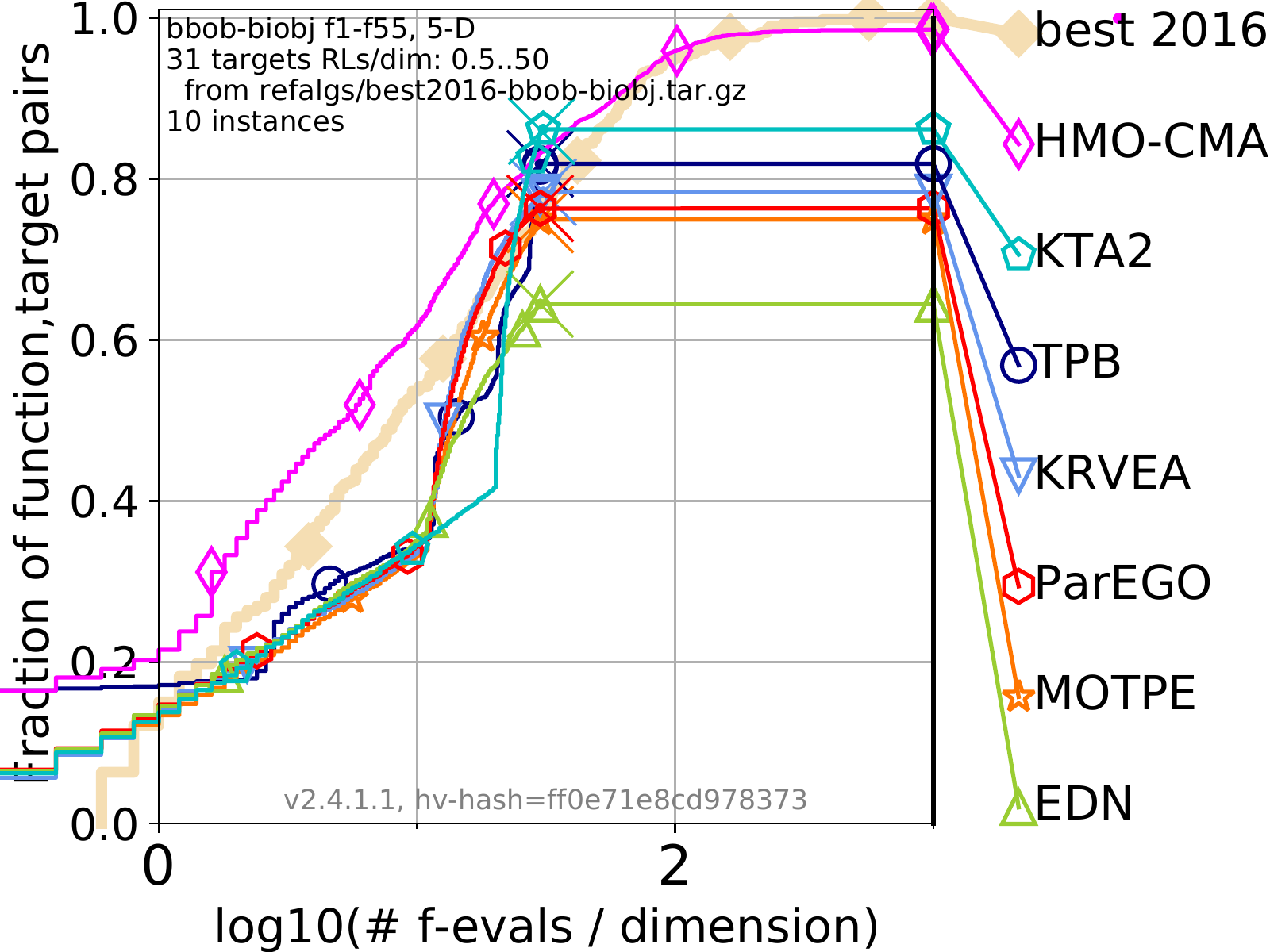}
}
\subfloat[$N=10$ and \texttt{budget} $=30 \times N$]{  
\includegraphics[width=0.24\textwidth]{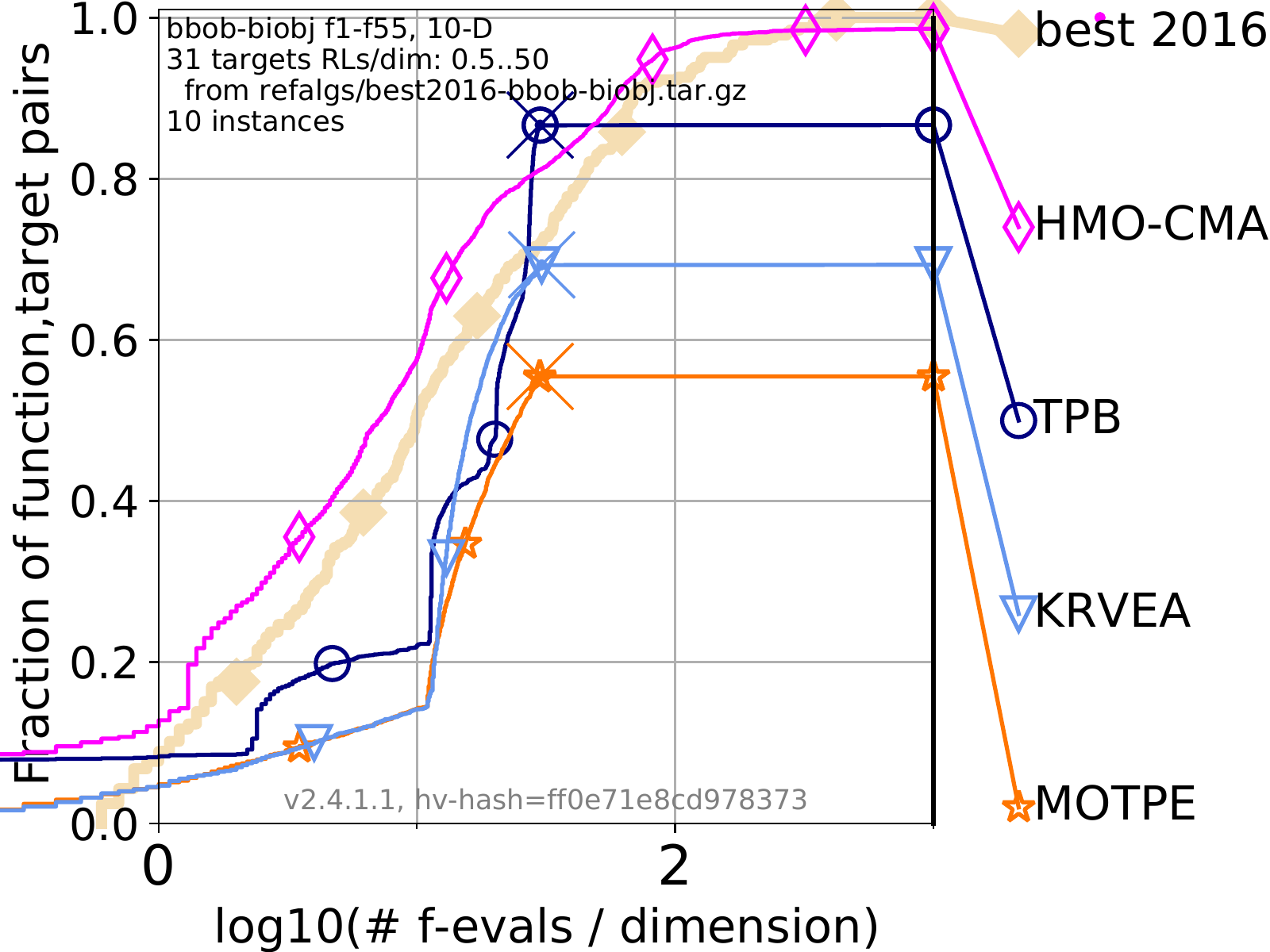}
}
\subfloat[$N=20$ and \texttt{budget} $=30 \times N$]{  
\includegraphics[width=0.24\textwidth]{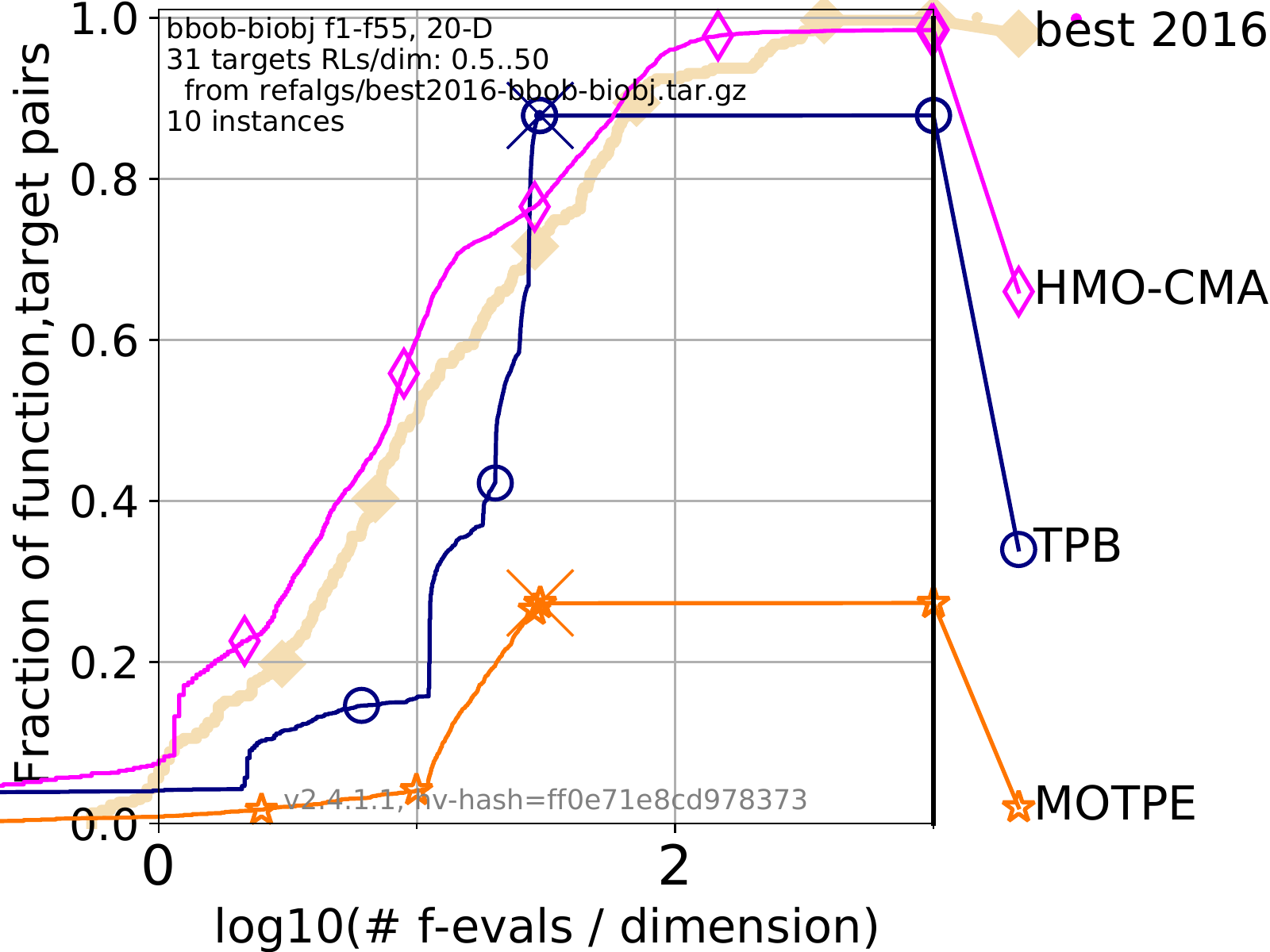}
}
\\
\subfloat[$N=2$ and \texttt{budget} $=40 \times N$]{  
\includegraphics[width=0.24\textwidth]{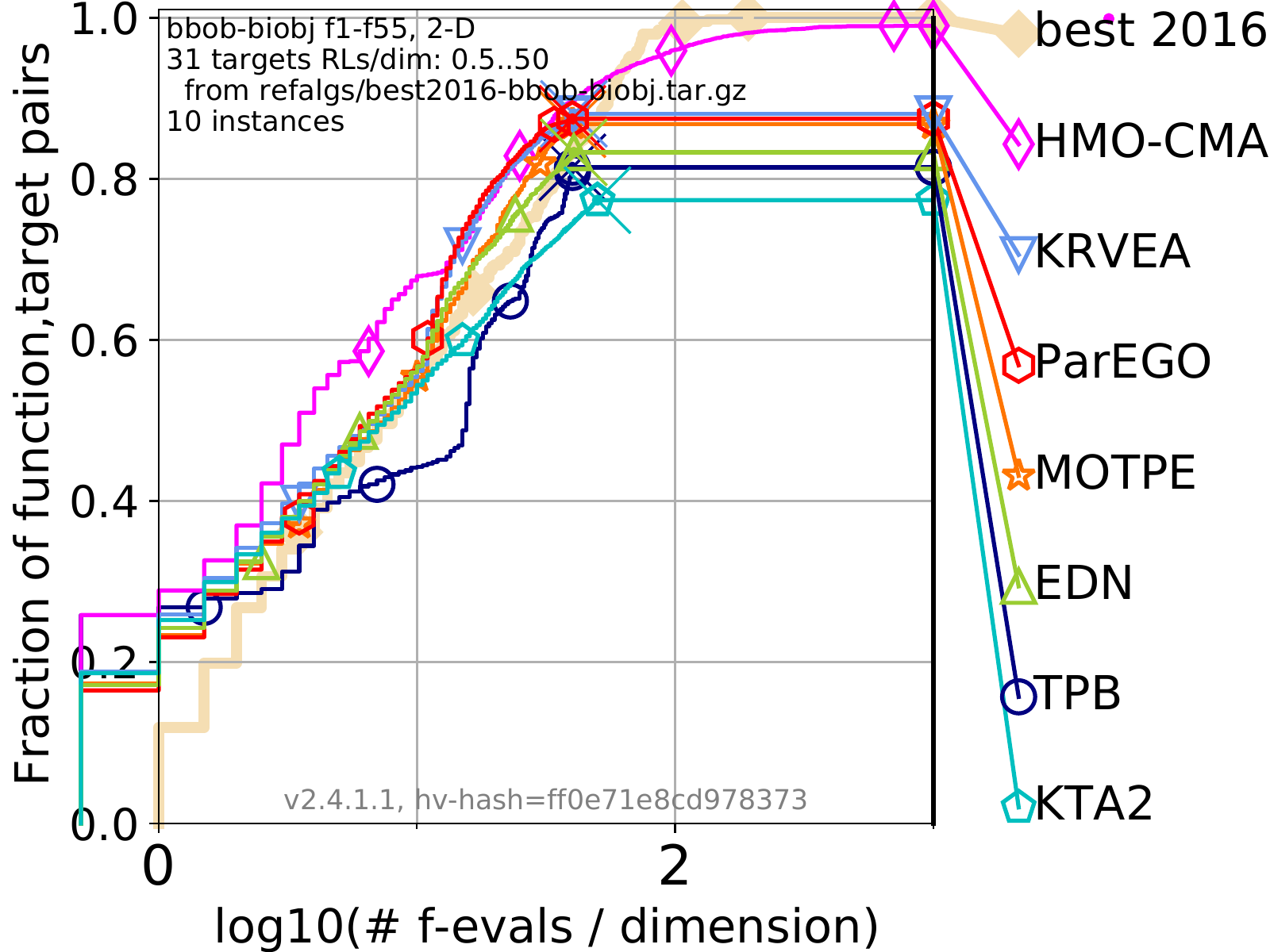}
}
\subfloat[$N=5$ and \texttt{budget} $=40 \times N$]{  
\includegraphics[width=0.24\textwidth]{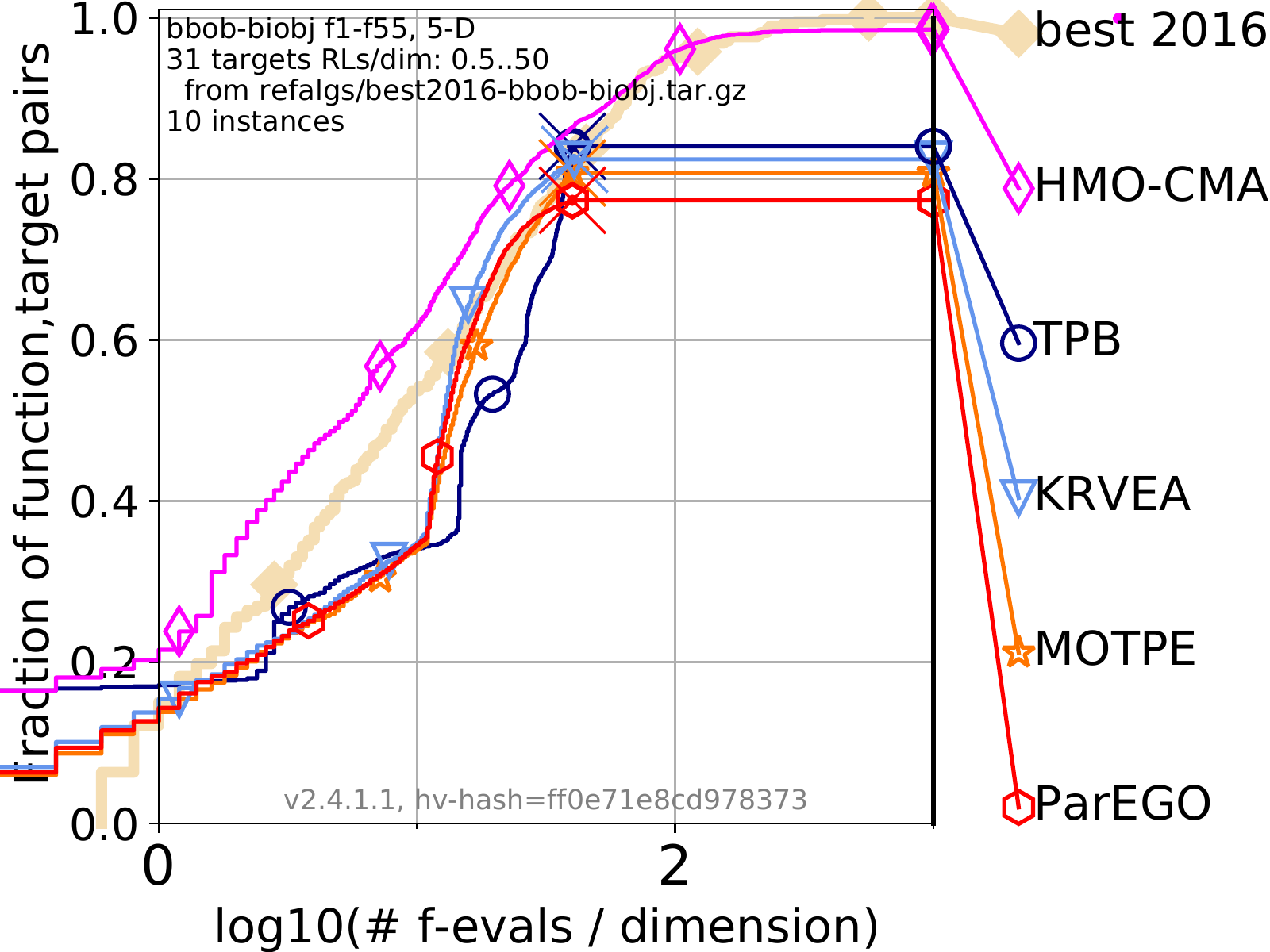}
}
\subfloat[$N=10$ and \texttt{budget} $=40 \times N$]{  
\includegraphics[width=0.24\textwidth]{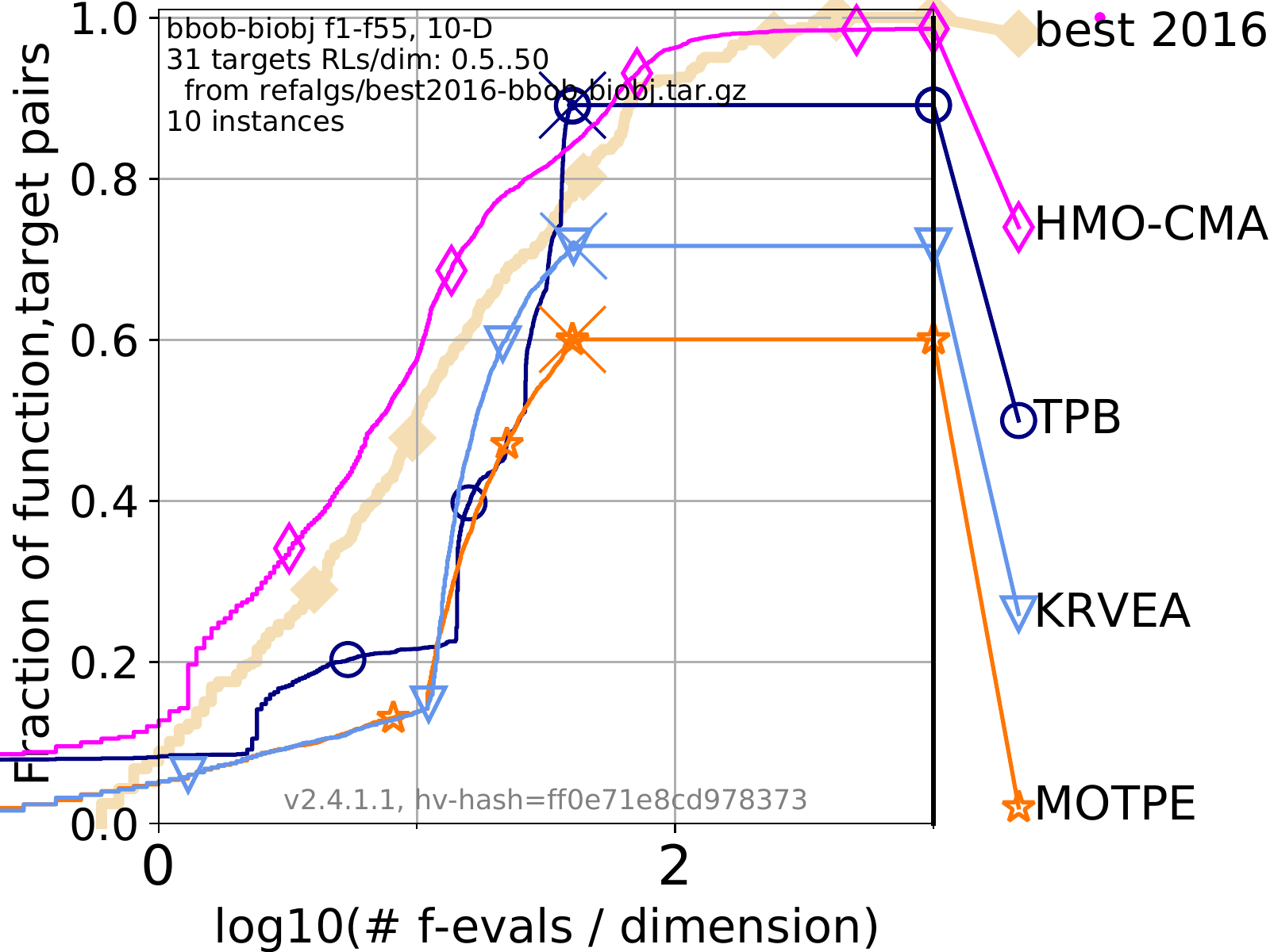}
}
\subfloat[$N=20$ and \texttt{budget} $=40 \times N$]{  
\includegraphics[width=0.24\textwidth]{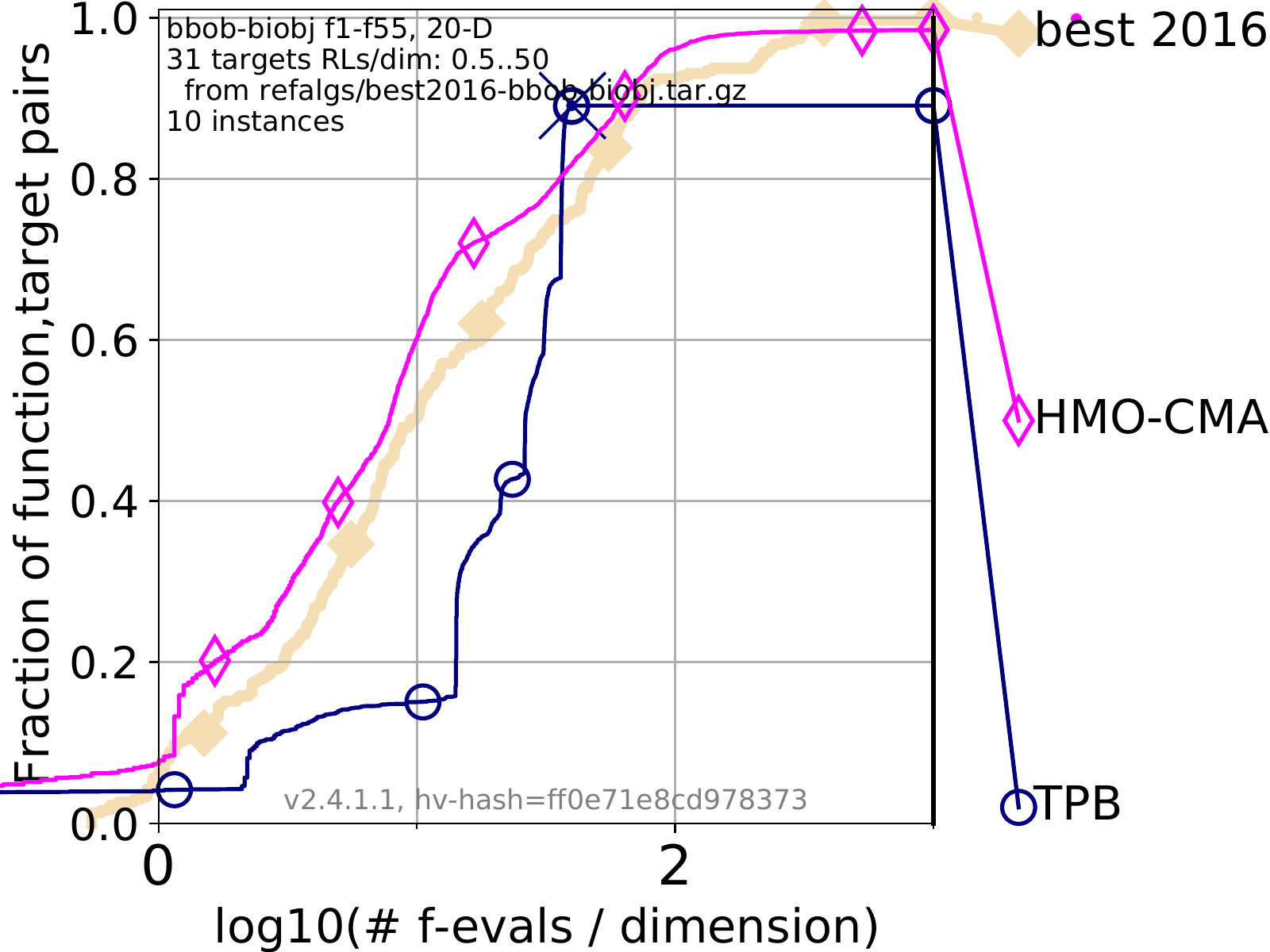}
}
\caption{Comparison with state-of-the-art optimizers. ``HMO-CMA'' and ``EDN'' stand for HMO-CMA-ES and EDN-ARMOEA, respectively.}
   \label{fig:vs_sota}
\end{figure*}

Figure \ref{fig:vs_sota} shows the results of TPB and the six optimizers on the 55 bi-objective BBOB problems with $N \in \{2, 5, 10, 20\}$ and \texttt{budget} $\in \{20 \times N, 30 \times N, 40 \times N\}$.
Recall that \texttt{budget} is the maximum budget of function evaluations.
We do not show the results for $N=3$, but they are similar to the results for $N=2$.
Most meta-model-based optimizers require extremely high computational cost, especially for higher dimensions and larger budgets.
Experiments on the 825 ($=55 \times 15$) BBOB instances for each dimension is also time-consuming.
For these reasons, we stopped an optimizer when it did not finish within a week.
The missing results in Figure \ref{fig:vs_sota} indicate that the corresponding optimizer was stopped before reaching \texttt{budget}, e.g., the results of KTA2 for $N=10$ in Figure \ref{fig:vs_sota}(c).
In Figures \ref{fig:vs_sota}, ``best 2016'' shows the performance of a \textit{virtual} best solver constructed based on the results of 15 optimizers participating in the GECCO BBOB 2016 workshop.
Thus, ``best 2016'' does not mean the best actual optimizer.
The cross in Figure \ref{fig:vs_sota} shows the number of function evaluations used in each optimizer.
Since ParEGO, K-RVEA, KTA2, and EDN-ARMOEA cannot stop exactly at a pre-defined \texttt{budget}, their crosses exceed \texttt{budget} in some cases, e.g., the results of KTA2 for $N=2$ in Figure \ref{fig:vs_sota}(a).

Figure \ref{fig:vs_sota} shows the bootstrapped empirical cumulative distribution (ECDF) \cite{BrockhoffTH15,HansenABTT16} based on the results on all 55 bi-objective BBOB problems.
We used the COCO postprocessing tool \texttt{cocopp} with the expensive option \texttt{--expensive} to generate all ECDF figures in this paper.
For each problem instance, let $\vector{I}^{\mathrm{ref}}$ be the $I_{\mathrm{COCO}}$ indicator value of the Pareto optimal solution set.
Let also $\vector{I}^{\mathrm{target}} = \vector{I}^{\mathrm{ref}} + \epsilon$ be a target value to reach, where $\epsilon$ is any one of 31 precision levels $\{0.5, ..., 50\}$ in the expensive setting.
Thus, 31 $\vector{I}^{\mathrm{target}}$ values are available for each problem instance.
The vertical axis in the ECDF figure represents the proportion of $I_{\rm COCO}$ values reached by the corresponding optimizer within specified function evaluations.
Here, the horizontal axis represents the number of function evaluations.
For example, Figure~\ref{fig:vs_sota}(d) indicates that HMO-CMA-ES solved about 60 \% of the 31 $\vector{I}^{\mathrm{target}}$ values within $10 \times N$ evaluations for $N=20$.
%

\textbf{Statistical significance} is tested with the rank-sum test  for a given $\vector{I}^{\mathrm{target}}$ value by using COCO.
Due to space limitation, we show the results in the supplementary material.
Note that the statistical test results are generally consistent with the results in Figure~\ref{fig:vs_sota}.

As shown in Figure~\ref{fig:vs_sota}, HMO-CMA-ES is the clear winner within $10 \times N$ function evaluations for any $N$.
The five meta-model-based optimizers perform almost the same until $11 \times N - 1$ function evaluations.
This is because they generate the initial solution set of size $11 \times N - 1$ by Latin hypercube sampling.
These results suggest that scalarization-based approaches with BOBYQA as in HMO-CMA-ES perform the best when only a very small number of function evaluations (i.e., $10 \times N$ evaluations) are available.

Some meta-model-based optimizers (e.g., ParEGO and K-RVEA) perform better than HMO-CMA-ES for more than $11 \times N - 1$ evaluations, especially for larger budgets.
We observed that the ranks of some meta-model-based optimizers depend on the maximum budget.
For example, for $N=5$, as shown in Figure \ref{fig:vs_sota}(b), KTA2 performs the worst when \texttt{budget} $=20 \times N$.
In contrast, as shown in Figure \ref{fig:vs_sota}(f), KTA2 performs the best at the end of the run when \texttt{budget} $=30 \times N$.
These observations indicate that the performance of some meta-model-based optimizers is sensitive to \texttt{budget}.
One may wonder about the high performance of ParEGO.
We believe that this is due to the performance evaluation based on the unbounded external archive.
Although an analysis of the performance of meta-model-based optimizers is beyond the scope of this paper, it is an interesting research direction.


As seen from Figures~\ref{fig:vs_sota}(a), (e), and (i), TPB performs poorly compared to the state-of-the-art optimizers for $N=2$.
In contrast, TPB achieves a good performance at the end of the run for $N=5$.
As shown in Figures~\ref{fig:vs_sota}(c), (d), (g), (h), (k), and (l), TPB is the best performer at the end of the run for $N=10$ and $N=20$.
These results indicate the effectiveness of TPB for \texttt{budget}$=20 \times N$, $30 \times N$, and $40 \times N$ for $N \geq 10$.

\begin{figure}[t]
\centering
\includegraphics[width=0.3\textwidth]{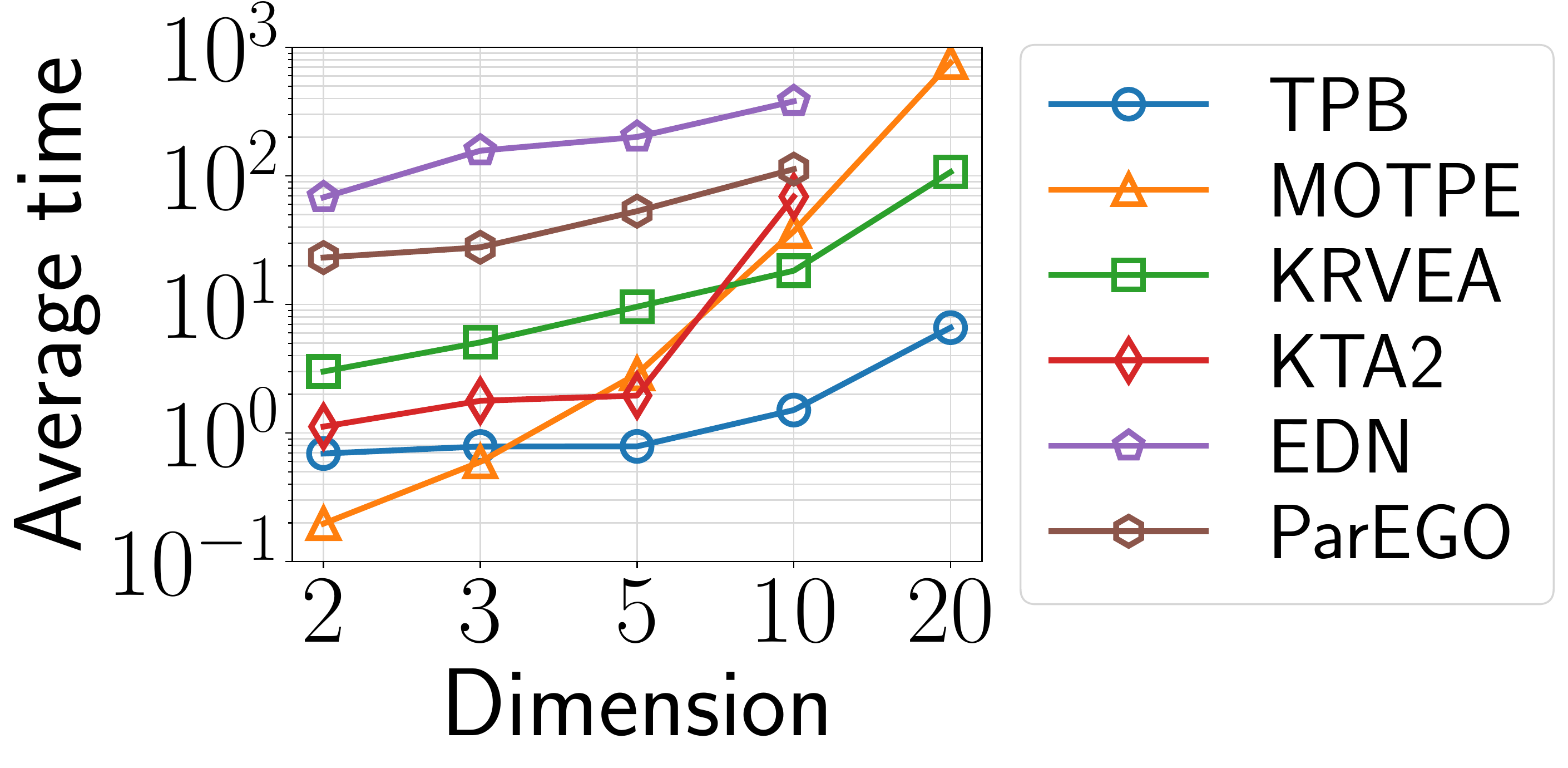}
\caption{Average computation time (sec) of each optimizer over the 15 instances of $f_1$ for \texttt{budget} $=20 \times N$.}
\label{fig:time}
\end{figure}

\begin{figure*}[t]
   \centering
\subfloat[$f_{1}$]{  
\includegraphics[width=0.24\textwidth]{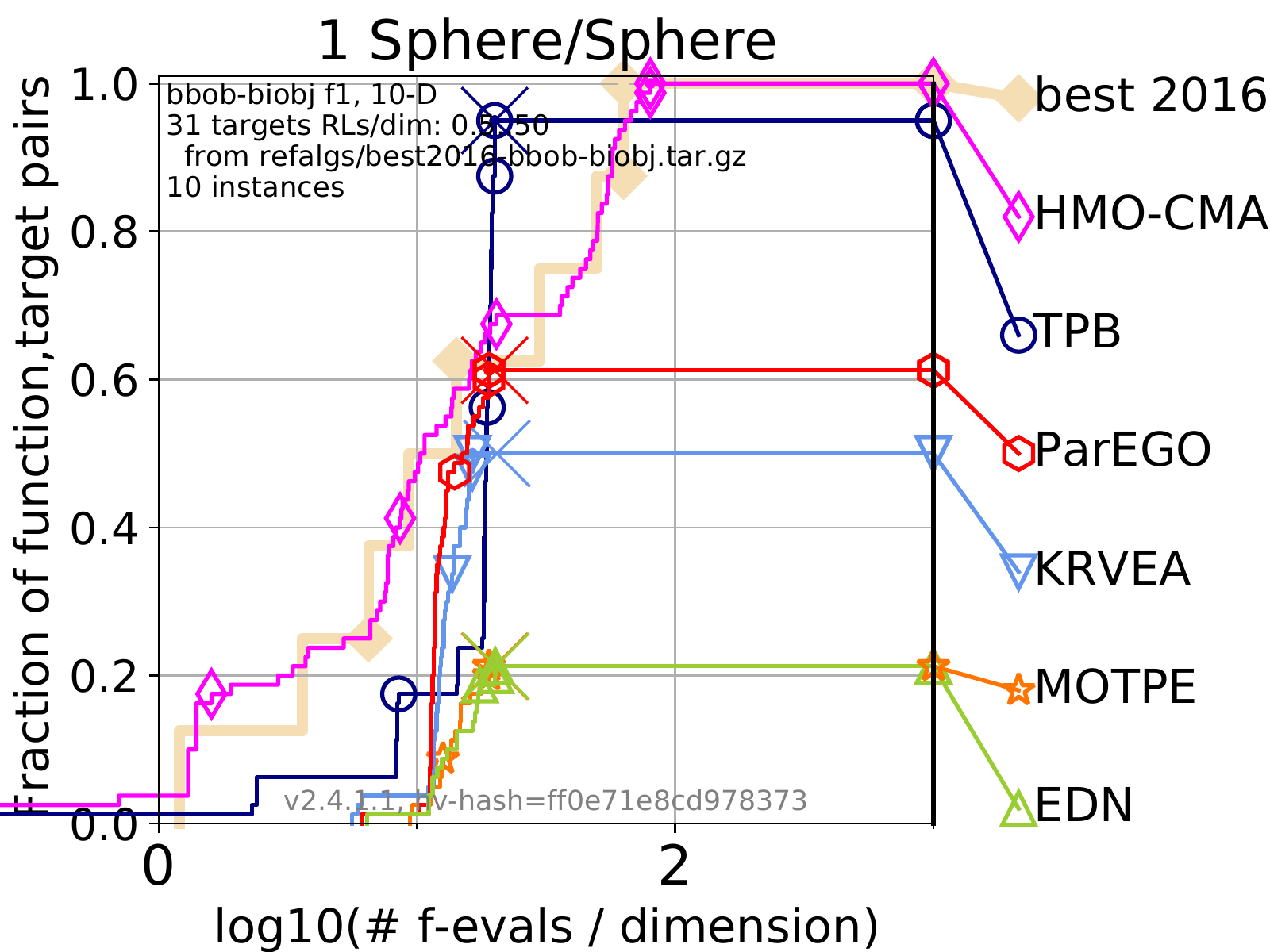}
}
\subfloat[$f_{28}$]{  
\includegraphics[width=0.24\textwidth]{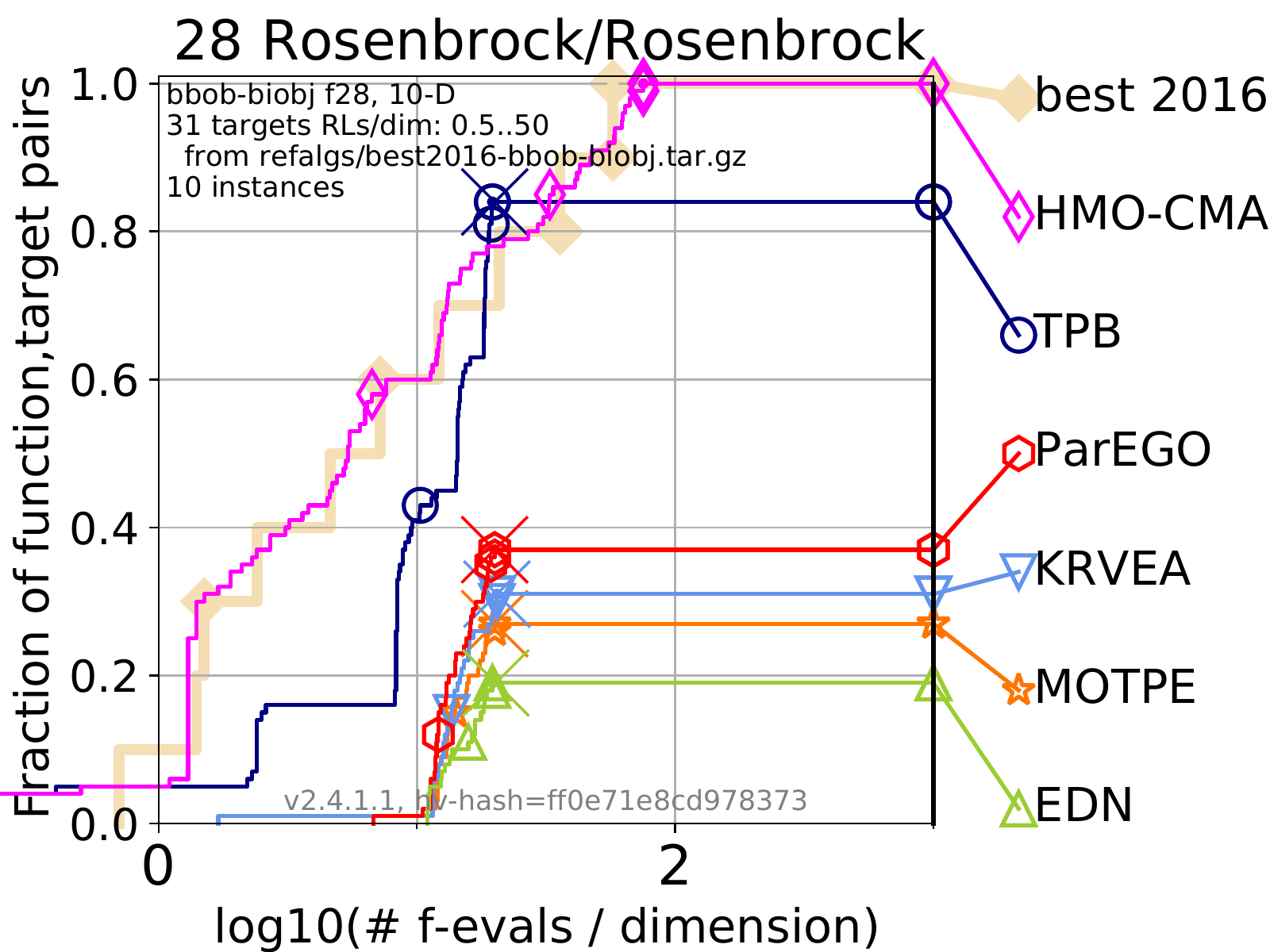}
}
\subfloat[$f_{46}$]{  
\includegraphics[width=0.24\textwidth]{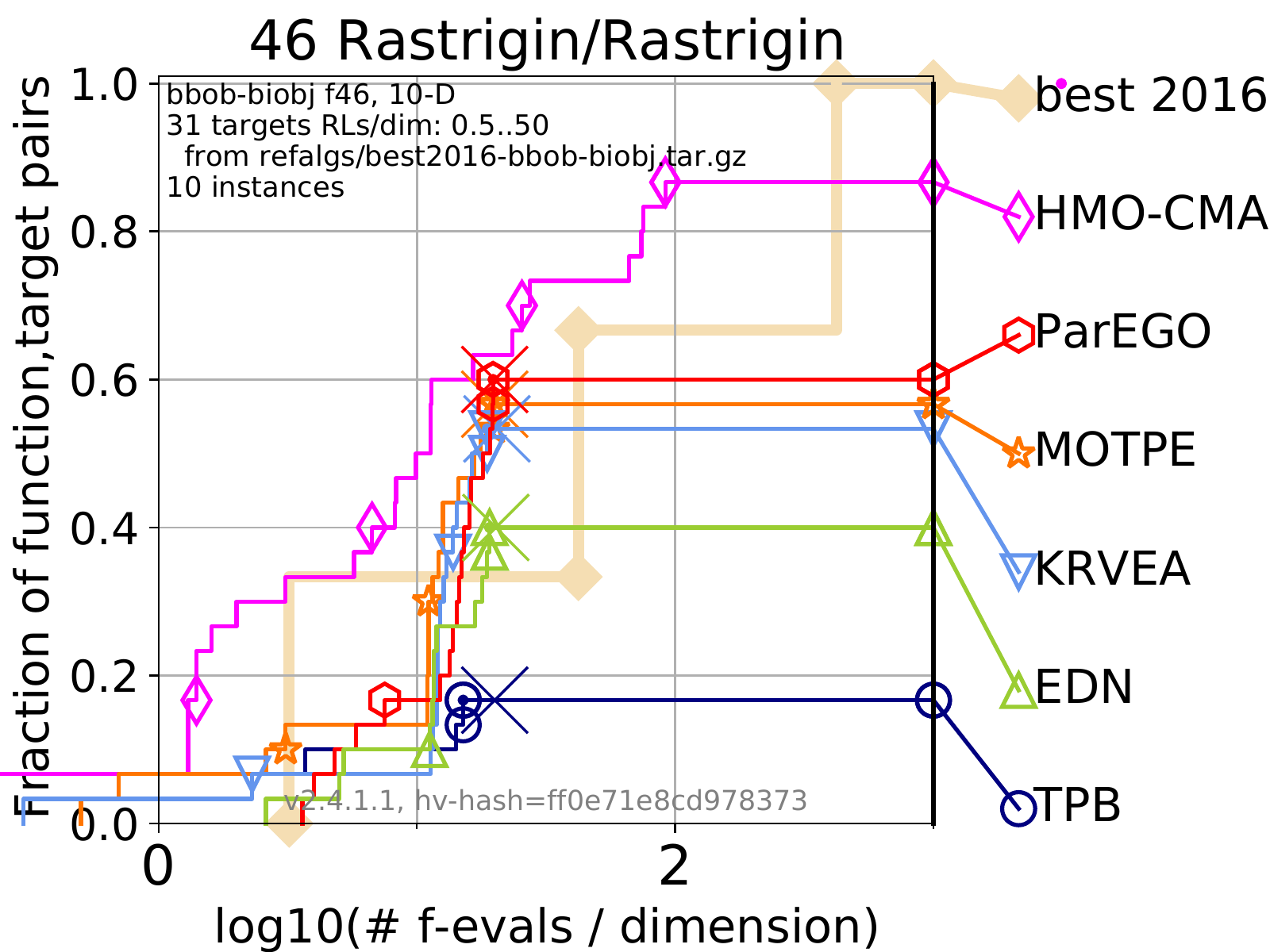}
}
\subfloat[$f_{53}$]{  
\includegraphics[width=0.24\textwidth]{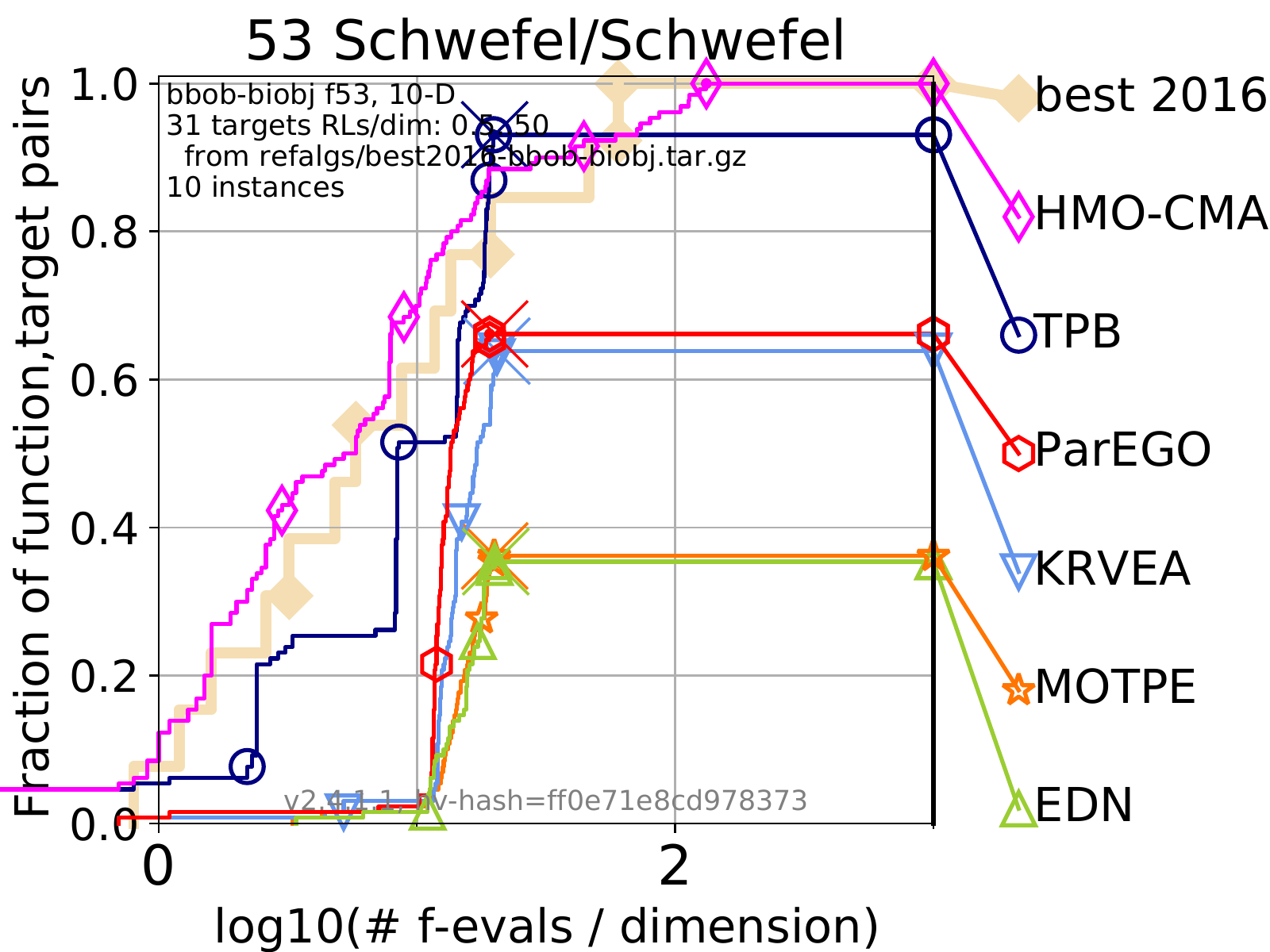}
}
\caption{Comparison with state-of-the-art optimizers on four selected problems with $N=10$ (\texttt{budget} $=20 \times N$).}
   \label{fig:vs_sota_eachf}
\end{figure*}

Figure \ref{fig:time} shows the average computation time of each optimizer over the 15 instances of $f_1$ for \texttt{budget} $=20 \times N$.
We expect that the computation time of HMO-CMA-ES is the same or less than that of TPB.
We could not measure the computation time of ParEGO, KTA2, and EDN-ARMOEA for $N=20$ in practical time due to their high computational cost.
As seen from Figure \ref{fig:time}, the computation time of TPB is lower than those of the five meta-model-based optimizers, except for the results of MOTPE for $N\leq 3$.
The computation of TPB took approximately 6.6 seconds even for $N=20$.
These results indicate that TPB is faster than meta-model-based optimizers in terms of computation time.

Figure \ref{fig:vs_sota_eachf} shows the results on $f_{1}$, $f_{28}$, $f_{46}$, and $f_{53}$, which are the multi-objective versions of the Sphere, Rosenbrock, (rotated) Rastrigin, and (rotated) Schwefel functions.
As discussed in Section \ref{sec:tpb_discussion}, the B\'{e}zier simplex model-based interpolation method assumes that a given problem is simplicial.
Although an in-depth theoretical analysis is needed, we believe that the 15 unimodal (and weakly-multimodal) bi-objective BBOB problems satisfy the assumption, including $f_{1}$ and $f_{28}$.
As shown in Figures \ref{fig:vs_sota_eachf}(a) and (b), TPB obtains a good performance on $f_{1}$ and $f_{28}$.
The results on other unimodal problems (except for $f_{11}$, $f_{12}$, and $f_{20}$) are relatively similar to Figures \ref{fig:vs_sota_eachf}(a) and (b).
In contrast, the remaining 40 multi-modal bi-objective BBOB problems do not satisfy the assumption, including $f_{46}$ and $f_{53}$.
As seen from Figure \ref{fig:vs_sota_eachf}(c), the poor performance of TPB on $f_{46}$ is consistent with our intuition.
 However, Figure \ref{fig:vs_sota_eachf}(d) shows that TPB unexpectedly performs the best on $f_{53}$.
Similar results were observed on other ten multimodal problems (e.g., $f_{9}$ and $f_{10}$).
These results suggest that the solution interpolation method can possibly perform well even when a given problem is not simplicial.
A further investigation is needed in future research.

In summary, we demonstrated the effectiveness of TPB for computationally expensive multi-objective optimization.
Our results on the bi-objective BBOB problems show that TPB performs better than HMO-CMA-ES and meta-model-based optimizers for $N \geq 10$.
We also observed that TPB is computationally cheaper than meta-model-based optimizers for $N \geq 5$.

\subsection{Importance of the two-phase mechanism}
\label{sec:inv_tp}

Here, let us consider the first phase-only TPB (TPB1) and the second phase-only TPB (TPB2).
We investigate the importance of the two-phase mechanism in TPB by comparing it with TPB1 and TPB2.
While TPB1 does not perform the second phase, TPB2 does not perform the first phase.
First, as in most meta-model-based optimizers (e.g., K-RVEA), TPB2 generates the initial solution set of size $11 \times N - 1$ by Latin hypercube sampling.
Then, TPB2 performs the second phase based on the best $K$ out of the $11 \times N - 1$ solutions.

\begin{figure}[t]
  \centering
\subfloat[$N=2$]{  
\includegraphics[width=0.22\textwidth]{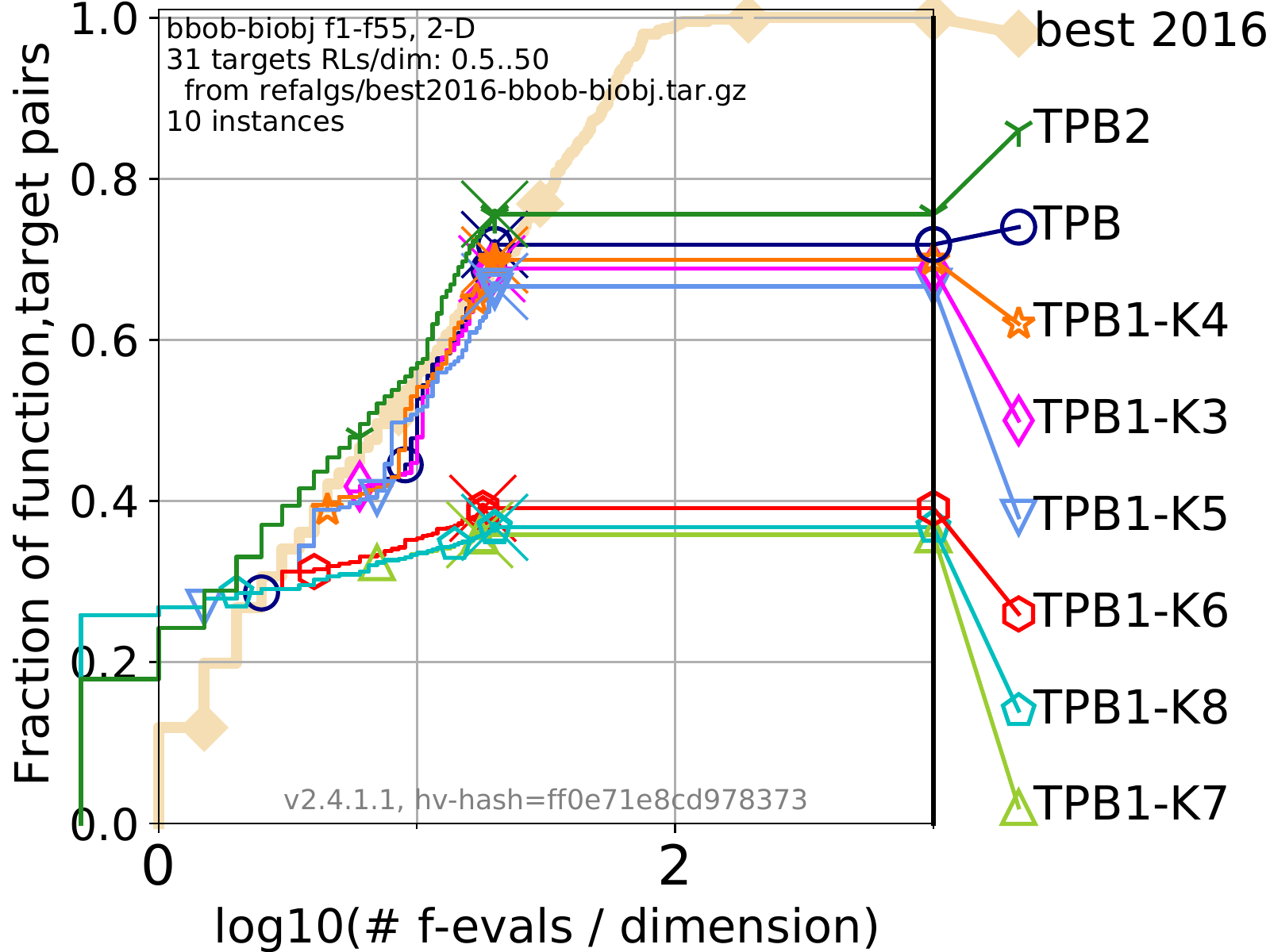}
}
\subfloat[$N=10$]{  
\includegraphics[width=0.22\textwidth]{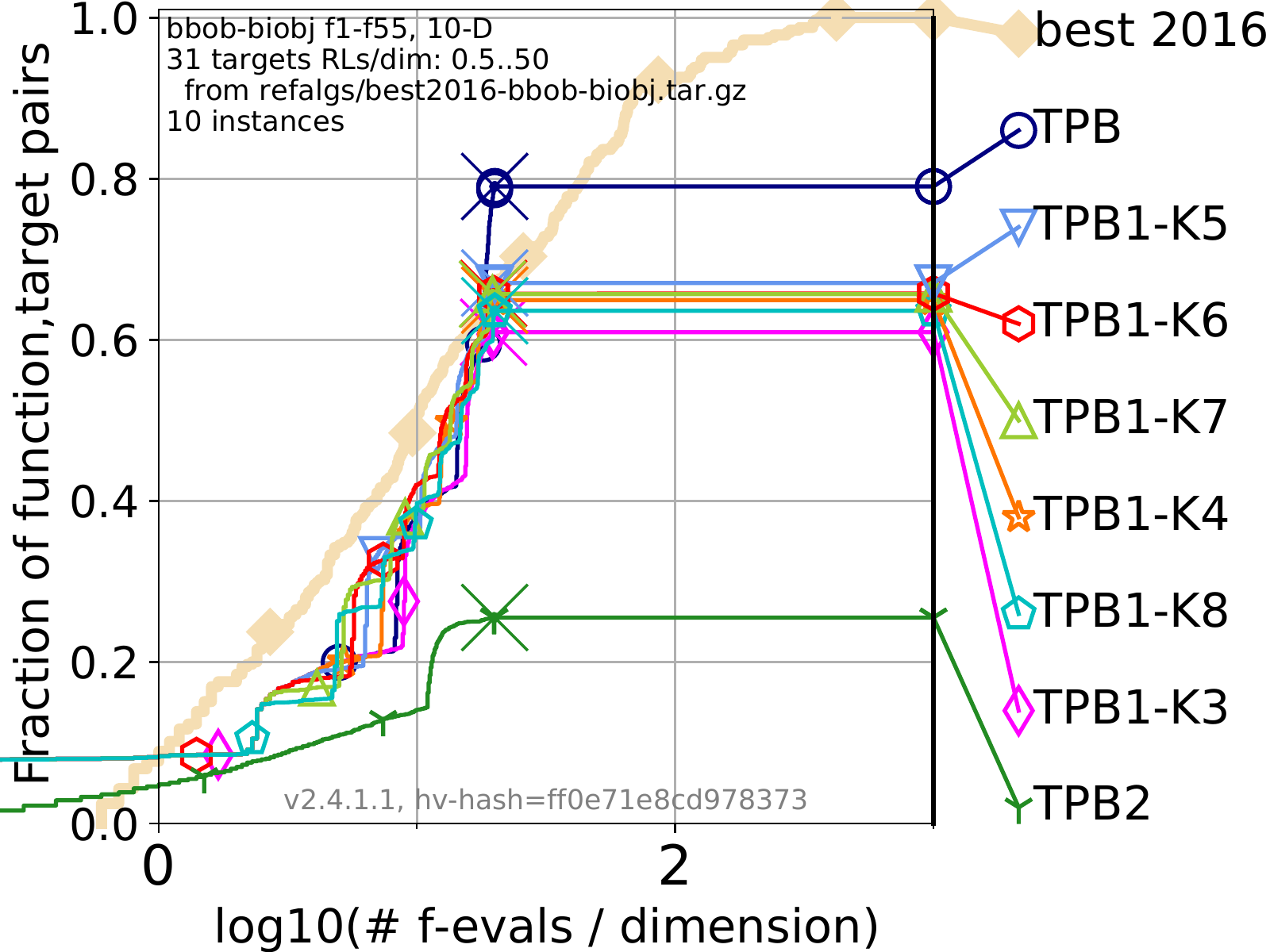}
}
\caption{Comparison of TPB, TPB1, and TPB2.}
  \label{fig:vs_tpb1_tpb2}
\end{figure}

Figure \ref{fig:vs_tpb1_tpb2} shows the comparison of TPB, TPB1, and TPB2 on the 55 bi-objective BBOB problems with $N=2$ and $10$ for \texttt{budget} $=20 \times N$.
Note that the results for $N \in \{3, 5, 20\}$ are similar to the results for $N=10$.
The results show that TPB1 performs worse than TPB at the end of the run for $N=2$ and $10$.
Interestingly, as shown in Figure \ref{fig:vs_tpb1_tpb2}(a), TPB2 outperforms TPB for $N=2$.
We observed that TPB2 performs well on multimodal problems for $N=2$.
However, as seen from Figure \ref{fig:vs_tpb1_tpb2}(b), TPB2 performs significantly poorly for $N=10$.
These results demonstrate the effectiveness of the two-phase mechanism in TPB.

\subsection{Impact of $K$ and $r^{\mathrm{1st}}$}
\label{sec:para_study}

\begin{figure}[t]
  \centering
\subfloat[\texttt{budget}$=20 \times N$]{  
\includegraphics[width=0.21\textwidth]{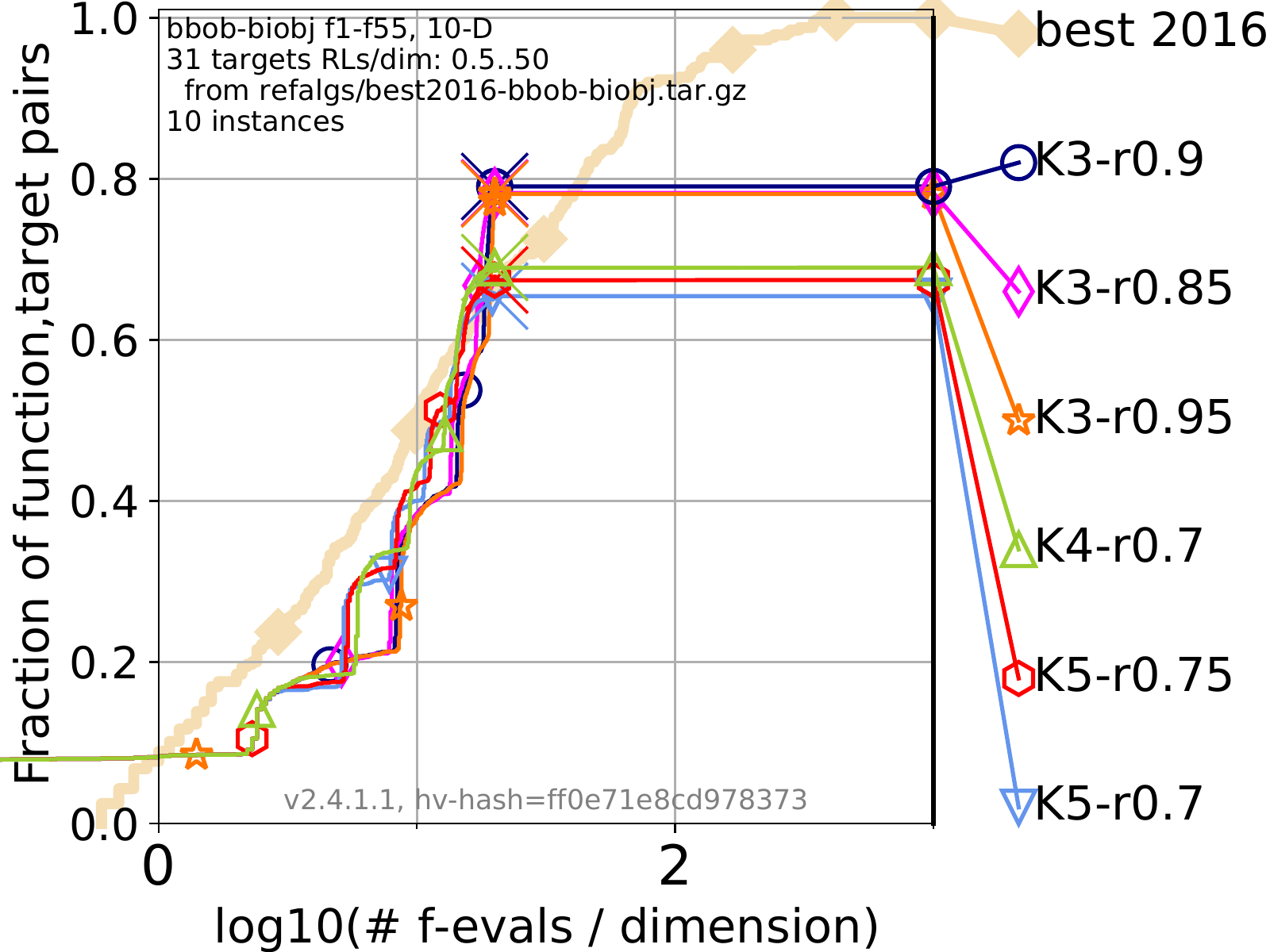}
}
\subfloat[\texttt{budget}$=40 \times N$]{  
\includegraphics[width=0.21\textwidth]{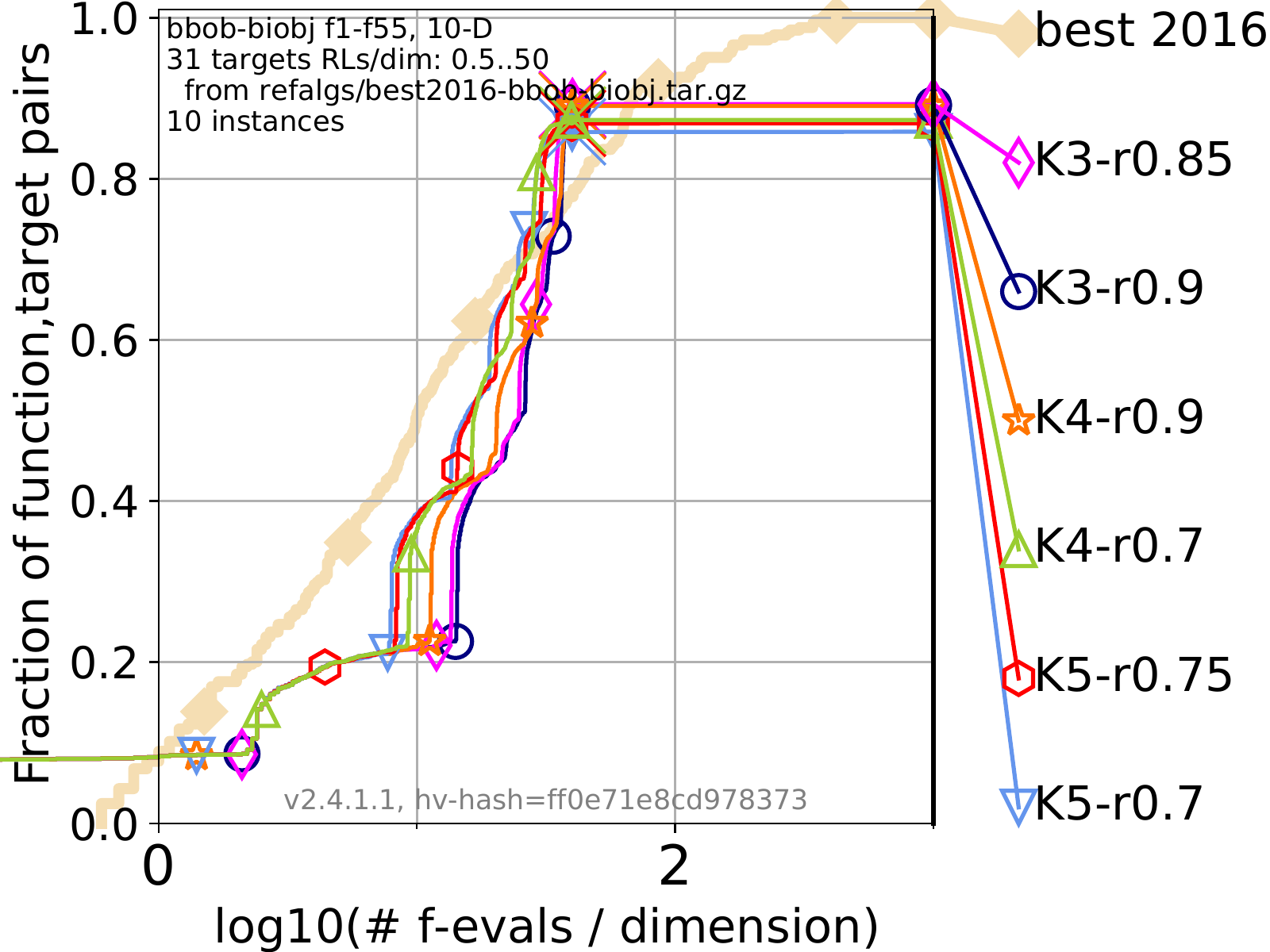}
}
\caption{Comparison of TPB with various $K$ and $r^{\mathrm{1st}}$.}
  \label{fig:param_tpb}
\end{figure}

Although Section \ref{sec:param_tpb} gave the default values of $K$ and $r^{\mathrm{1st}}$, it is important to understand their impact on the performance of TPB.
Figure \ref{fig:param_tpb} shows the results of TPB with $K \in \{3, 4, 5\}$ and $r^{\mathrm{1st}} \in \{0.7, 0.75, 0.8, 0.85, 0.9, 0.95\}$ on the 55 bi-objective BBOB problems for $N=10$, where \texttt{budget} $=20 \times N$ and $40 \times N$.
For example, ``K3-r0.9'' represents the results of TPB with $K=3$ and $r^{\mathrm{1st}}=0.9$.
For the sake of clarity, Figure \ref{fig:param_tpb} shows only the results of TPB with the three best parameter settings and the three worst parameter settings.
As seen from Figure \ref{fig:param_tpb}(a), the best performance of TPB for \texttt{budget} $=20 \times N$ is obtained when using $K=3$ and $r^{\mathrm{1st}}=0.9$.
In contrast, as shown in Figure \ref{fig:param_tpb}(b), TPB with $K=3$ and $r^{\mathrm{1st}}=0.85$ performs the best for \texttt{budget} $=40 \times N$.
Figure \ref{fig:param_tpb} shows that the gap between the best and worst performance of TPB is relatively small for \texttt{budget} $=40 \times N$.
Although we do not show detailed results here, we observed that the best setting of $K$ and $r^{\mathrm{1st}}$ depends on a problem, $N$, and \texttt{budget}.
These results suggest that the performance of TPB can be further improved by tuning the $K$ and $r^{\mathrm{1st}}$ values.
However, $K=3$ and $r^{\mathrm{1st}}=0.9$ can be a good first choice for $M=2$.

\section{Conclusion}
\label{sec:conclusion}


We have proposed TPB for computationally expensive multi-objective black-box optimization.
The first phase in TPB fully exploits an efficient derivative-free optimizer to find well-approximated solutions of $K$ scalar problems with a small budget of function evaluations, where $K=M+1$.
The second phase in TPB interpolates the $K$ solutions by the B\'{e}zier simplex model-based method that exploits the property of the Pareto optimal solution set.
Our results show that TPB performs significantly better than HMO-CMA-ES and some state-of-the-art meta-model-based multi-objective optimizers on the bi-objective BBOB problems with $N \geq 10$ when the maximum budget of function evaluations is set to $20 \times N$, $30 \times N$, and $40 \times N$.
We have also investigated the property of TPB.

We believe that TPB gives a new perspective on the field of computationally expensive multi-objective optimization.
Although the EMO community has mainly focused on meta-model-based approaches for computationally expensive optimization, TPB provides a new research direction.
It may also be interesting to extend TPB to preference-based multi-objective optimization.



\begin{acks}
This work was supported by JSPS KAKENHI Grant Number 21K17824. 
We thank Dr. Ilya Loshchilov for providing the code of HMO-CMA-ES.
\end{acks}

\bibliographystyle{ACM-Reference-Format}
\bibliography{reference} 


\begin{thebibliography}{67}


\ifx \showCODEN    \undefined \def \showCODEN     #1{\unskip}     \fi
\ifx \showDOI      \undefined \def \showDOI       #1{#1}\fi
\ifx \showISBNx    \undefined \def \showISBNx     #1{\unskip}     \fi
\ifx \showISBNxiii \undefined \def \showISBNxiii  #1{\unskip}     \fi
\ifx \showISSN     \undefined \def \showISSN      #1{\unskip}     \fi
\ifx \showLCCN     \undefined \def \showLCCN      #1{\unskip}     \fi
\ifx \shownote     \undefined \def \shownote      #1{#1}          \fi
\ifx \showarticletitle \undefined \def \showarticletitle #1{#1}   \fi
\ifx \showURL      \undefined \def \showURL       {\relax}        \fi
\providecommand\bibfield[2]{#2}
\providecommand\bibinfo[2]{#2}
\providecommand\natexlab[1]{#1}
\providecommand\showeprint[2][]{arXiv:#2}

\bibitem[\protect\citeauthoryear{Akiba, Sano, Yanase, Ohta, and Koyama}{Akiba
  et~al\mbox{.}}{2019}]%
        {AkibaSYOK19}
\bibfield{author}{\bibinfo{person}{Takuya Akiba}, \bibinfo{person}{Shotaro
  Sano}, \bibinfo{person}{Toshihiko Yanase}, \bibinfo{person}{Takeru Ohta},
  {and} \bibinfo{person}{Masanori Koyama}.} \bibinfo{year}{2019}\natexlab{}.
\newblock \showarticletitle{Optuna: {A} Next-generation Hyperparameter
  Optimization Framework}. In \bibinfo{booktitle}{\emph{Proceedings of the 25th
  {ACM} {SIGKDD} International Conference on Knowledge Discovery {\&} Data
  Mining, {KDD} 2019, Anchorage, AK, USA, August 4-8, 2019}},
  \bibfield{editor}{\bibinfo{person}{Ankur Teredesai}, \bibinfo{person}{Vipin
  Kumar}, \bibinfo{person}{Ying Li}, \bibinfo{person}{R{\'{o}}mer Rosales},
  \bibinfo{person}{Evimaria Terzi}, {and} \bibinfo{person}{George Karypis}}
  (Eds.). \bibinfo{publisher}{{ACM}}, \bibinfo{pages}{2623--2631}.
\newblock
\urldef\tempurl%
\url{https://doi.org/10.1145/3292500.3330701}
\showDOI{\tempurl}


\bibitem[\protect\citeauthoryear{Bajer, Pitra, Repick{\'{y}}, and Holena}{Bajer
  et~al\mbox{.}}{2019}]%
        {BajerPRH19}
\bibfield{author}{\bibinfo{person}{Luk{\'{a}}s Bajer}, \bibinfo{person}{Zbynek
  Pitra}, \bibinfo{person}{Jakub Repick{\'{y}}}, {and} \bibinfo{person}{Martin
  Holena}.} \bibinfo{year}{2019}\natexlab{}.
\newblock \showarticletitle{Gaussian Process Surrogate Models for the {CMA}
  Evolution Strategy}.
\newblock \bibinfo{journal}{\emph{Evol. Comput.}} \bibinfo{volume}{27},
  \bibinfo{number}{4} (\bibinfo{year}{2019}), \bibinfo{pages}{665--697}.
\newblock
\urldef\tempurl%
\url{https://doi.org/10.1162/evco\_a\_00244}
\showDOI{\tempurl}


\bibitem[\protect\citeauthoryear{Beume, Naujoks, and Emmerich}{Beume
  et~al\mbox{.}}{2007}]%
        {BeumeNE07}
\bibfield{author}{\bibinfo{person}{Nicola Beume}, \bibinfo{person}{Boris
  Naujoks}, {and} \bibinfo{person}{Michael T.~M. Emmerich}.}
  \bibinfo{year}{2007}\natexlab{}.
\newblock \showarticletitle{{SMS-EMOA:} Multiobjective selection based on
  dominated hypervolume}.
\newblock \bibinfo{journal}{\emph{Eur. J. Oper. Res.}} \bibinfo{volume}{181},
  \bibinfo{number}{3} (\bibinfo{year}{2007}), \bibinfo{pages}{1653--1669}.
\newblock
\urldef\tempurl%
\url{https://doi.org/10.1016/j.ejor.2006.08.008}
\showDOI{\tempurl}


\bibitem[\protect\citeauthoryear{Bezerra, L{\'{o}}pez{-}Ib{\'{a}}{\~{n}}ez, and
  St{\"{u}}tzle}{Bezerra et~al\mbox{.}}{2018}]%
        {BezerraLS18}
\bibfield{author}{\bibinfo{person}{Leonardo C.~T. Bezerra},
  \bibinfo{person}{Manuel L{\'{o}}pez{-}Ib{\'{a}}{\~{n}}ez}, {and}
  \bibinfo{person}{Thomas St{\"{u}}tzle}.} \bibinfo{year}{2018}\natexlab{}.
\newblock \showarticletitle{A Large-Scale Experimental Evaluation of
  High-Performing Multi- and Many-Objective Evolutionary Algorithms}.
\newblock \bibinfo{journal}{\emph{Evol. Comput.}} \bibinfo{volume}{26},
  \bibinfo{number}{4} (\bibinfo{year}{2018}).
\newblock
\urldef\tempurl%
\url{https://doi.org/10.1162/evco\_a\_00217}
\showDOI{\tempurl}


\bibitem[\protect\citeauthoryear{Bhattacharjee, Singh, and Ray}{Bhattacharjee
  et~al\mbox{.}}{2017}]%
        {BhattacharjeeSR17J}
\bibfield{author}{\bibinfo{person}{Kalyan~Shankar Bhattacharjee},
  \bibinfo{person}{Hemant~Kumar Singh}, {and} \bibinfo{person}{Tapabrata Ray}.}
  \bibinfo{year}{2017}\natexlab{}.
\newblock \showarticletitle{An approach to generate comprehensive piecewise
  linear interpolation of pareto outcomes to aid decision making}.
\newblock \bibinfo{journal}{\emph{J. Glob. Optim.}} \bibinfo{volume}{68},
  \bibinfo{number}{1} (\bibinfo{year}{2017}), \bibinfo{pages}{71--93}.
\newblock
\urldef\tempurl%
\url{https://doi.org/10.1007/s10898-016-0454-0}
\showDOI{\tempurl}


\bibitem[\protect\citeauthoryear{Borges and Pastva}{Borges and Pastva}{2002}]%
        {BorgesP2002}
\bibfield{author}{\bibinfo{person}{Carlos~F. Borges} {and} \bibinfo{person}{Tim
  Pastva}.} \bibinfo{year}{2002}\natexlab{}.
\newblock \showarticletitle{Total least squares fitting of {B\'ezier} and
  {B-spline} curves to ordered data}.
\newblock \bibinfo{journal}{\emph{Computer Aided Geometric Design}}
  \bibinfo{volume}{19}, \bibinfo{number}{4} (\bibinfo{year}{2002}),
  \bibinfo{pages}{275--289}.
\newblock
\urldef\tempurl%
\url{https://doi.org/10.1016/s0167-8396(02)00088-2}
\showDOI{\tempurl}


\bibitem[\protect\citeauthoryear{Bouzarkouna, Auger, and Ding}{Bouzarkouna
  et~al\mbox{.}}{2011}]%
        {BouzarkounaAD11}
\bibfield{author}{\bibinfo{person}{Zyed Bouzarkouna}, \bibinfo{person}{Anne
  Auger}, {and} \bibinfo{person}{Didier~Yu Ding}.}
  \bibinfo{year}{2011}\natexlab{}.
\newblock \showarticletitle{Local-meta-model {CMA-ES} for partially separable
  functions}. In \bibinfo{booktitle}{\emph{13th Annual Genetic and Evolutionary
  Computation Conference, {GECCO} 2011, Proceedings, Dublin, Ireland, July
  12-16, 2011}}, \bibfield{editor}{\bibinfo{person}{Natalio Krasnogor} {and}
  \bibinfo{person}{Pier~Luca Lanzi}} (Eds.). \bibinfo{publisher}{{ACM}},
  \bibinfo{pages}{869--876}.
\newblock
\urldef\tempurl%
\url{https://doi.org/10.1145/2001576.2001695}
\showDOI{\tempurl}


\bibitem[\protect\citeauthoryear{Brockhoff, Auger, Hansen, and Tu{\v
  s}ar}{Brockhoff et~al\mbox{.}}{ress}]%
        {BrockhoffAHT22}
\bibfield{author}{\bibinfo{person}{Dimo Brockhoff}, \bibinfo{person}{Anne
  Auger}, \bibinfo{person}{Nikolaus Hansen}, {and} \bibinfo{person}{Tea Tu{\v
  s}ar}.} \bibinfo{year}{2022 (in press)}\natexlab{}.
\newblock \showarticletitle{{ Using Well-Understood Single-Objective Functions
  in Multiobjective Black-Box Optimization Test Suites}}.
\newblock \bibinfo{journal}{\emph{Evol. Comput.}} (\bibinfo{year}{2022 (in
  press)}).
\newblock


\bibitem[\protect\citeauthoryear{Brockhoff, Plaquevent{-}Jourdain, Auger, and
  Hansen}{Brockhoff et~al\mbox{.}}{2021}]%
        {BrockhoffPAH21}
\bibfield{author}{\bibinfo{person}{Dimo Brockhoff}, \bibinfo{person}{Baptiste
  Plaquevent{-}Jourdain}, \bibinfo{person}{Anne Auger}, {and}
  \bibinfo{person}{Nikolaus Hansen}.} \bibinfo{year}{2021}\natexlab{}.
\newblock \showarticletitle{{DMS} and MultiGLODS: black-box optimization
  benchmarking of two direct search methods on the bbob-biobj test suite}. In
  \bibinfo{booktitle}{\emph{{GECCO} '21: Genetic and Evolutionary Computation
  Conference, Companion Volume, Lille, France, July 10-14, 2021}},
  \bibfield{editor}{\bibinfo{person}{Krzysztof Krawiec}} (Ed.).
  \bibinfo{publisher}{{ACM}}, \bibinfo{pages}{1251--1258}.
\newblock
\urldef\tempurl%
\url{https://doi.org/10.1145/3449726.3463207}
\showDOI{\tempurl}


\bibitem[\protect\citeauthoryear{Brockhoff, Tran, and Hansen}{Brockhoff
  et~al\mbox{.}}{2015}]%
        {BrockhoffTH15}
\bibfield{author}{\bibinfo{person}{Dimo Brockhoff}, \bibinfo{person}{Thanh{-}Do
  Tran}, {and} \bibinfo{person}{Nikolaus Hansen}.}
  \bibinfo{year}{2015}\natexlab{}.
\newblock \showarticletitle{Benchmarking Numerical Multiobjective Optimizers
  Revisited}. In \bibinfo{booktitle}{\emph{Proceedings of the Genetic and
  Evolutionary Computation Conference, {GECCO} 2015, Madrid, Spain, July 11-15,
  2015}}, \bibfield{editor}{\bibinfo{person}{Sara Silva} {and}
  \bibinfo{person}{Anna~Isabel Esparcia{-}Alc{\'{a}}zar}} (Eds.).
  \bibinfo{publisher}{{ACM}}, \bibinfo{pages}{639--646}.
\newblock
\urldef\tempurl%
\url{https://doi.org/10.1145/2739480.2754777}
\showDOI{\tempurl}


\bibitem[\protect\citeauthoryear{Brockhoff, Tusar, Tusar, Wagner, Hansen, and
  Auger}{Brockhoff et~al\mbox{.}}{2016}]%
        {BrockhoffTTWHA16}
\bibfield{author}{\bibinfo{person}{Dimo Brockhoff}, \bibinfo{person}{Tea
  Tusar}, \bibinfo{person}{Dejan Tusar}, \bibinfo{person}{Tobias Wagner},
  \bibinfo{person}{Nikolaus Hansen}, {and} \bibinfo{person}{Anne Auger}.}
  \bibinfo{year}{2016}\natexlab{}.
\newblock \showarticletitle{Biobjective Performance Assessment with the {COCO}
  Platform}.
\newblock \bibinfo{journal}{\emph{CoRR}}  \bibinfo{volume}{abs/1605.01746}
  (\bibinfo{year}{2016}).
\newblock
\showeprint[arXiv]{1605.01746}
\urldef\tempurl%
\url{http://arxiv.org/abs/1605.01746}
\showURL{%
\tempurl}


\bibitem[\protect\citeauthoryear{Cartis, Fiala, Marteau, and Roberts}{Cartis
  et~al\mbox{.}}{2019}]%
        {CartisFMR19}
\bibfield{author}{\bibinfo{person}{Coralia Cartis}, \bibinfo{person}{Jan
  Fiala}, \bibinfo{person}{Benjamin Marteau}, {and} \bibinfo{person}{Lindon
  Roberts}.} \bibinfo{year}{2019}\natexlab{}.
\newblock \showarticletitle{Improving the Flexibility and Robustness of
  Model-based Derivative-free Optimization Solvers}.
\newblock \bibinfo{journal}{\emph{{ACM} Trans. Math. Softw.}}
  \bibinfo{volume}{45}, \bibinfo{number}{3} (\bibinfo{year}{2019}),
  \bibinfo{pages}{32:1--32:41}.
\newblock
\urldef\tempurl%
\url{https://doi.org/10.1145/3338517}
\showDOI{\tempurl}


\bibitem[\protect\citeauthoryear{Chen, Ishibuchi, and Shang}{Chen
  et~al\mbox{.}}{2020}]%
        {ChenIS20}
\bibfield{author}{\bibinfo{person}{Weiyu Chen}, \bibinfo{person}{Hisao
  Ishibuchi}, {and} \bibinfo{person}{Ke Shang}.}
  \bibinfo{year}{2020}\natexlab{}.
\newblock \showarticletitle{Proposal of a Realistic Many-Objective Test Suite}.
  In \bibinfo{booktitle}{\emph{Parallel Problem Solving from Nature - {PPSN}
  {XVI} - 16th International Conference, {PPSN} 2020, Leiden, The Netherlands,
  September 5-9, 2020, Proceedings, Part {I}}} \emph{(\bibinfo{series}{Lecture
  Notes in Computer Science}, Vol.~\bibinfo{volume}{12269})},
  \bibfield{editor}{\bibinfo{person}{Thomas B{\"{a}}ck}, \bibinfo{person}{Mike
  Preuss}, \bibinfo{person}{Andr{\'{e}}~H. Deutz}, \bibinfo{person}{Hao Wang},
  \bibinfo{person}{Carola Doerr}, \bibinfo{person}{Michael T.~M. Emmerich},
  {and} \bibinfo{person}{Heike Trautmann}} (Eds.).
  \bibinfo{publisher}{Springer}, \bibinfo{pages}{201--214}.
\newblock


\bibitem[\protect\citeauthoryear{Chugh, Jin, Miettinen, Hakanen, and
  Sindhya}{Chugh et~al\mbox{.}}{2018}]%
        {ChughJMHS18}
\bibfield{author}{\bibinfo{person}{Tinkle Chugh}, \bibinfo{person}{Yaochu Jin},
  \bibinfo{person}{Kaisa Miettinen}, \bibinfo{person}{Jussi Hakanen}, {and}
  \bibinfo{person}{Karthik Sindhya}.} \bibinfo{year}{2018}\natexlab{}.
\newblock \showarticletitle{A Surrogate-Assisted Reference Vector Guided
  Evolutionary Algorithm for Computationally Expensive Many-Objective
  Optimization}.
\newblock \bibinfo{journal}{\emph{{IEEE} Trans. Evol. Comput.}}
  \bibinfo{volume}{22}, \bibinfo{number}{1} (\bibinfo{year}{2018}),
  \bibinfo{pages}{129--142}.
\newblock
\urldef\tempurl%
\url{https://doi.org/10.1109/TEVC.2016.2622301}
\showDOI{\tempurl}


\bibitem[\protect\citeauthoryear{Chugh, Sindhya, Hakanen, and Miettinen}{Chugh
  et~al\mbox{.}}{2019}]%
        {ChughSHM19}
\bibfield{author}{\bibinfo{person}{Tinkle Chugh}, \bibinfo{person}{Karthik
  Sindhya}, \bibinfo{person}{Jussi Hakanen}, {and} \bibinfo{person}{Kaisa
  Miettinen}.} \bibinfo{year}{2019}\natexlab{}.
\newblock \showarticletitle{A survey on handling computationally expensive
  multiobjective optimization problems with evolutionary algorithms}.
\newblock \bibinfo{journal}{\emph{Soft Comput.}} \bibinfo{volume}{23},
  \bibinfo{number}{9} (\bibinfo{year}{2019}), \bibinfo{pages}{3137--3166}.
\newblock
\urldef\tempurl%
\url{https://doi.org/10.1007/s00500-017-2965-0}
\showDOI{\tempurl}


\bibitem[\protect\citeauthoryear{Cust{\'{o}}dio, Madeira, Vaz, and
  Vicente}{Cust{\'{o}}dio et~al\mbox{.}}{2011}]%
        {CustodioMVV11}
\bibfield{author}{\bibinfo{person}{Ana~Lu{\'{\i}}sa Cust{\'{o}}dio},
  \bibinfo{person}{J.~F.~Aguilar Madeira}, \bibinfo{person}{A.~Ismael~F. Vaz},
  {and} \bibinfo{person}{Lu{\'{\i}}s~Nunes Vicente}.}
  \bibinfo{year}{2011}\natexlab{}.
\newblock \showarticletitle{Direct Multisearch for Multiobjective
  Optimization}.
\newblock \bibinfo{journal}{\emph{{SIAM} J. Optim.}} \bibinfo{volume}{21},
  \bibinfo{number}{3} (\bibinfo{year}{2011}), \bibinfo{pages}{1109--1140}.
\newblock
\urldef\tempurl%
\url{https://doi.org/10.1137/10079731X}
\showDOI{\tempurl}


\bibitem[\protect\citeauthoryear{Daniels, Rahat, Everson, Tabor, and
  Fieldsend}{Daniels et~al\mbox{.}}{2018}]%
        {DanielsRETF18}
\bibfield{author}{\bibinfo{person}{Steven~J. Daniels}, \bibinfo{person}{Alma
  As{-}Aad~Mohammad Rahat}, \bibinfo{person}{Richard~M. Everson},
  \bibinfo{person}{Gavin~R. Tabor}, {and} \bibinfo{person}{Jonathan~E.
  Fieldsend}.} \bibinfo{year}{2018}\natexlab{}.
\newblock \showarticletitle{A Suite of Computationally Expensive Shape
  Optimisation Problems Using Computational Fluid Dynamics}. In
  \bibinfo{booktitle}{\emph{Parallel Problem Solving from Nature - {PPSN} {XV}
  - 15th International Conference, Coimbra, Portugal, September 8-12, 2018,
  Proceedings, Part {II}}} \emph{(\bibinfo{series}{Lecture Notes in Computer
  Science}, Vol.~\bibinfo{volume}{11102})},
  \bibfield{editor}{\bibinfo{person}{Anne Auger}, \bibinfo{person}{Carlos~M.
  Fonseca}, \bibinfo{person}{Nuno Louren{\c{c}}o}, \bibinfo{person}{Penousal
  Machado}, \bibinfo{person}{Lu{\'{\i}}s Paquete}, {and}
  \bibinfo{person}{L.~Darrell Whitley}} (Eds.). \bibinfo{publisher}{Springer},
  \bibinfo{pages}{296--307}.
\newblock
\urldef\tempurl%
\url{https://doi.org/10.1007/978-3-319-99259-4\_24}
\showDOI{\tempurl}


\bibitem[\protect\citeauthoryear{Deb, Agrawal, Pratap, and Meyarivan}{Deb
  et~al\mbox{.}}{2002}]%
        {DebAPM02}
\bibfield{author}{\bibinfo{person}{Kalyanmoy Deb}, \bibinfo{person}{Samir
  Agrawal}, \bibinfo{person}{Amrit Pratap}, {and} \bibinfo{person}{T.
  Meyarivan}.} \bibinfo{year}{2002}\natexlab{}.
\newblock \showarticletitle{A fast and elitist multiobjective genetic
  algorithm: {NSGA-II}}.
\newblock \bibinfo{journal}{\emph{{IEEE} Trans. Evol. Comput.}}
  \bibinfo{volume}{6}, \bibinfo{number}{2} (\bibinfo{year}{2002}),
  \bibinfo{pages}{182--197}.
\newblock
\urldef\tempurl%
\url{https://doi.org/10.1109/4235.996017}
\showDOI{\tempurl}


\bibitem[\protect\citeauthoryear{Deb, Thiele, Laumanns, and Zitzler}{Deb
  et~al\mbox{.}}{2005}]%
        {DebTLZ05}
\bibfield{author}{\bibinfo{person}{Kalyanmoy Deb}, \bibinfo{person}{Lothar
  Thiele}, \bibinfo{person}{Marco Laumanns}, {and} \bibinfo{person}{Eckart
  Zitzler}.} \bibinfo{year}{2005}\natexlab{}.
\newblock \showarticletitle{{Scalable Test Problems for Evolutionary
  Multi-Objective Optimization}}.
\newblock In \bibinfo{booktitle}{\emph{Evolutionary Multiobjective
  Optimization. Theoretical Advances and Applications}}.
  \bibinfo{publisher}{Springer}, \bibinfo{pages}{105--145}.
\newblock


\bibitem[\protect\citeauthoryear{Dubois{-}Lacoste,
  L{\'{o}}pez{-}Ib{\'{a}}{\~{n}}ez, and St{\"{u}}tzle}{Dubois{-}Lacoste
  et~al\mbox{.}}{2011}]%
        {Dubois-LacosteLS11amai}
\bibfield{author}{\bibinfo{person}{J{\'{e}}r{\'{e}}mie Dubois{-}Lacoste},
  \bibinfo{person}{Manuel L{\'{o}}pez{-}Ib{\'{a}}{\~{n}}ez}, {and}
  \bibinfo{person}{Thomas St{\"{u}}tzle}.} \bibinfo{year}{2011}\natexlab{}.
\newblock \showarticletitle{Improving the anytime behavior of two-phase local
  search}.
\newblock \bibinfo{journal}{\emph{Ann. Math. Artif. Intell.}}
  \bibinfo{volume}{61}, \bibinfo{number}{2} (\bibinfo{year}{2011}),
  \bibinfo{pages}{125--154}.
\newblock
\urldef\tempurl%
\url{https://doi.org/10.1007/s10472-011-9235-0}
\showDOI{\tempurl}


\bibitem[\protect\citeauthoryear{Dubois{-}Lacoste,
  L{\'{o}}pez{-}Ib{\'{a}}{\~{n}}ez, and St{\"{u}}tzle}{Dubois{-}Lacoste
  et~al\mbox{.}}{2013}]%
        {Dubois-LacosteLS13}
\bibfield{author}{\bibinfo{person}{J{\'{e}}r{\'{e}}mie Dubois{-}Lacoste},
  \bibinfo{person}{Manuel L{\'{o}}pez{-}Ib{\'{a}}{\~{n}}ez}, {and}
  \bibinfo{person}{Thomas St{\"{u}}tzle}.} \bibinfo{year}{2013}\natexlab{}.
\newblock \showarticletitle{Combining Two Search Paradigms for Multi-objective
  Optimization: Two-Phase and Pareto Local Search}.
\newblock In \bibinfo{booktitle}{\emph{Hybrid Metaheuristics}},
  \bibfield{editor}{\bibinfo{person}{El{-}Ghazali Talbi}} (Ed.).
  \bibinfo{series}{Studies in Computational Intelligence},
  Vol.~\bibinfo{volume}{434}. \bibinfo{publisher}{Springer},
  \bibinfo{pages}{97--117}.
\newblock
\urldef\tempurl%
\url{https://doi.org/10.1007/978-3-642-30671-6\_3}
\showDOI{\tempurl}


\bibitem[\protect\citeauthoryear{Dymond, Kok, and Heyns}{Dymond
  et~al\mbox{.}}{2013}]%
        {DymondKH13}
\bibfield{author}{\bibinfo{person}{Antoine S.~D. Dymond},
  \bibinfo{person}{Schalk Kok}, {and} \bibinfo{person}{P.~Stephan Heyns}.}
  \bibinfo{year}{2013}\natexlab{}.
\newblock \showarticletitle{The sensitivity of multi-objective optimization
  algorithm performance to objective function evaluation budgets}. In
  \bibinfo{booktitle}{\emph{Proceedings of the {IEEE} Congress on Evolutionary
  Computation, {CEC} 2013, Cancun, Mexico, June 20-23, 2013}}.
  \bibinfo{publisher}{{IEEE}}, \bibinfo{pages}{1868--1875}.
\newblock
\urldef\tempurl%
\url{https://doi.org/10.1109/CEC.2013.6557787}
\showDOI{\tempurl}


\bibitem[\protect\citeauthoryear{Farin}{Farin}{2002}]%
        {Farin02}
\bibfield{author}{\bibinfo{person}{G.E. Farin}.}
  \bibinfo{year}{2002}\natexlab{}.
\newblock \bibinfo{booktitle}{\emph{Curves and Surfaces for {CAGD}: A Practical
  Guide}}.
\newblock \bibinfo{publisher}{Morgan Kaufmann}.
\newblock
\showISBNx{9781558607378}
\showLCCN{2001094373}
\urldef\tempurl%
\url{https://books.google.co.jp/books?id=5HYTP1dIAp4C}
\showURL{%
\tempurl}


\bibitem[\protect\citeauthoryear{Giagkiozis and Fleming}{Giagkiozis and
  Fleming}{2014}]%
        {GiagkiozisF14}
\bibfield{author}{\bibinfo{person}{Ioannis Giagkiozis} {and}
  \bibinfo{person}{Peter~J. Fleming}.} \bibinfo{year}{2014}\natexlab{}.
\newblock \showarticletitle{Pareto Front Estimation for Decision Making}.
\newblock \bibinfo{journal}{\emph{Evol. Comput.}} \bibinfo{volume}{22},
  \bibinfo{number}{4} (\bibinfo{year}{2014}), \bibinfo{pages}{651--678}.
\newblock
\urldef\tempurl%
\url{https://doi.org/10.1162/EVCO\_a\_00128}
\showDOI{\tempurl}


\bibitem[\protect\citeauthoryear{Goodfellow, Pouget{-}Abadie, Mirza, Xu,
  Warde{-}Farley, Ozair, Courville, and Bengio}{Goodfellow
  et~al\mbox{.}}{2014}]%
        {GoodfellowPMXWOCB14}
\bibfield{author}{\bibinfo{person}{Ian~J. Goodfellow}, \bibinfo{person}{Jean
  Pouget{-}Abadie}, \bibinfo{person}{Mehdi Mirza}, \bibinfo{person}{Bing Xu},
  \bibinfo{person}{David Warde{-}Farley}, \bibinfo{person}{Sherjil Ozair},
  \bibinfo{person}{Aaron~C. Courville}, {and} \bibinfo{person}{Yoshua Bengio}.}
  \bibinfo{year}{2014}\natexlab{}.
\newblock \showarticletitle{Generative Adversarial Nets}. In
  \bibinfo{booktitle}{\emph{Advances in Neural Information Processing Systems
  27: Annual Conference on Neural Information Processing Systems 2014, December
  8-13 2014, Montreal, Quebec, Canada}},
  \bibfield{editor}{\bibinfo{person}{Zoubin Ghahramani}, \bibinfo{person}{Max
  Welling}, \bibinfo{person}{Corinna Cortes}, \bibinfo{person}{Neil~D.
  Lawrence}, {and} \bibinfo{person}{Kilian~Q. Weinberger}} (Eds.).
  \bibinfo{pages}{2672--2680}.
\newblock
\urldef\tempurl%
\url{https://proceedings.neurips.cc/paper/2014/hash/5ca3e9b122f61f8f06494c97b1afccf3-Abstract.html}
\showURL{%
\tempurl}


\bibitem[\protect\citeauthoryear{Guo, Wang, Gao, Jin, Ding, and Chai}{Guo
  et~al\mbox{.}}{2022}]%
        {GuoWGJDC22}
\bibfield{author}{\bibinfo{person}{Dan Guo}, \bibinfo{person}{Xilu Wang},
  \bibinfo{person}{Kailai Gao}, \bibinfo{person}{Yaochu Jin},
  \bibinfo{person}{Jinliang Ding}, {and} \bibinfo{person}{Tianyou Chai}.}
  \bibinfo{year}{2022}\natexlab{}.
\newblock \showarticletitle{{Evolutionary Optimization of High-Dimensional
  Multiobjective and Many-Objective Expensive Problems Assisted by a Dropout
  Neural Network}}.
\newblock \bibinfo{journal}{\emph{{IEEE} Trans. Syst. Man Cybern. Syst.}}
  \bibinfo{volume}{52}, \bibinfo{number}{4} (\bibinfo{year}{2022}),
  \bibinfo{pages}{2084--2097}.
\newblock
\urldef\tempurl%
\url{https://doi.org/10.1109/TSMC.2020.3044418}
\showURL{%
\tempurl}


\bibitem[\protect\citeauthoryear{Hamada, Hayano, Ichiki, Kabata, and
  Teramoto}{Hamada et~al\mbox{.}}{2020}]%
        {HamadaHIKT2020}
\bibfield{author}{\bibinfo{person}{Naoki Hamada}, \bibinfo{person}{Kenta
  Hayano}, \bibinfo{person}{Shunsuke Ichiki}, \bibinfo{person}{Yutaro Kabata},
  {and} \bibinfo{person}{Hiroshi Teramoto}.} \bibinfo{year}{2020}\natexlab{}.
\newblock \showarticletitle{Topology of {Pareto} Sets of Strongly Convex
  Problems}.
\newblock \bibinfo{journal}{\emph{SIAM Journal on Optimization}}
  \bibinfo{volume}{30}, \bibinfo{number}{3} (\bibinfo{year}{2020}),
  \bibinfo{pages}{2659--2686}.
\newblock
\urldef\tempurl%
\url{https://doi.org/10.1137/19M1271439}
\showDOI{\tempurl}
\showeprint{https://doi.org/10.1137/19M1271439}


\bibitem[\protect\citeauthoryear{Hamada, Sakuma, Kobayashi, and Ono}{Hamada
  et~al\mbox{.}}{2008}]%
        {HamadaSKO08}
\bibfield{author}{\bibinfo{person}{Naoki Hamada}, \bibinfo{person}{Jun Sakuma},
  \bibinfo{person}{Shigenobu Kobayashi}, {and} \bibinfo{person}{Isao Ono}.}
  \bibinfo{year}{2008}\natexlab{}.
\newblock \showarticletitle{Functional-Specialization Multi-Objective
  Real-Coded Genetic Algorithm: {FS-MOGA}}. In
  \bibinfo{booktitle}{\emph{Parallel Problem Solving from Nature - {PPSN} X,
  10th International Conference Dortmund, Germany, September 13-17, 2008,
  Proceedings}} \emph{(\bibinfo{series}{Lecture Notes in Computer Science},
  Vol.~\bibinfo{volume}{5199})},
  \bibfield{editor}{\bibinfo{person}{G{\"{u}}nter Rudolph},
  \bibinfo{person}{Thomas Jansen}, \bibinfo{person}{Simon~M. Lucas},
  \bibinfo{person}{Carlo Poloni}, {and} \bibinfo{person}{Nicola Beume}} (Eds.).
  \bibinfo{publisher}{Springer}, \bibinfo{pages}{691--701}.
\newblock
\urldef\tempurl%
\url{https://doi.org/10.1007/978-3-540-87700-4\_69}
\showDOI{\tempurl}


\bibitem[\protect\citeauthoryear{Hansen}{Hansen}{2019}]%
        {Hansen19}
\bibfield{author}{\bibinfo{person}{Nikolaus Hansen}.}
  \bibinfo{year}{2019}\natexlab{}.
\newblock \showarticletitle{A global surrogate assisted {CMA-ES}}. In
  \bibinfo{booktitle}{\emph{Proceedings of the Genetic and Evolutionary
  Computation Conference, {GECCO} 2019, Prague, Czech Republic, July 13-17,
  2019}}, \bibfield{editor}{\bibinfo{person}{Anne Auger} {and}
  \bibinfo{person}{Thomas St{\"{u}}tzle}} (Eds.). \bibinfo{publisher}{{ACM}},
  \bibinfo{pages}{664--672}.
\newblock
\urldef\tempurl%
\url{https://doi.org/10.1145/3321707.3321842}
\showDOI{\tempurl}


\bibitem[\protect\citeauthoryear{Hansen, Auger, Brockhoff, Tusar, and
  Tusar}{Hansen et~al\mbox{.}}{2016}]%
        {HansenABTT16}
\bibfield{author}{\bibinfo{person}{Nikolaus Hansen}, \bibinfo{person}{Anne
  Auger}, \bibinfo{person}{Dimo Brockhoff}, \bibinfo{person}{Dejan Tusar},
  {and} \bibinfo{person}{Tea Tusar}.} \bibinfo{year}{2016}\natexlab{}.
\newblock \showarticletitle{{COCO:} Performance Assessment}.
\newblock \bibinfo{journal}{\emph{CoRR}}  \bibinfo{volume}{abs/1605.03560}
  (\bibinfo{year}{2016}).
\newblock
\showeprint[arXiv]{1605.03560}
\urldef\tempurl%
\url{http://arxiv.org/abs/1605.03560}
\showURL{%
\tempurl}


\bibitem[\protect\citeauthoryear{Hansen, Auger, Ros, Finck, and
  Pos{\'{\i}}k}{Hansen et~al\mbox{.}}{2010}]%
        {HansenARFP10}
\bibfield{author}{\bibinfo{person}{Nikolaus Hansen}, \bibinfo{person}{Anne
  Auger}, \bibinfo{person}{Raymond Ros}, \bibinfo{person}{Steffen Finck}, {and}
  \bibinfo{person}{Petr Pos{\'{\i}}k}.} \bibinfo{year}{2010}\natexlab{}.
\newblock \showarticletitle{Comparing results of 31 algorithms from the
  black-box optimization benchmarking {BBOB-2009}}. In
  \bibinfo{booktitle}{\emph{Genetic and Evolutionary Computation Conference,
  {GECCO} 2010, Proceedings, Portland, Oregon, USA, July 7-11, 2010, Companion
  Material}}, \bibfield{editor}{\bibinfo{person}{Martin Pelikan} {and}
  \bibinfo{person}{J{\"{u}}rgen Branke}} (Eds.). \bibinfo{publisher}{{ACM}},
  \bibinfo{pages}{1689--1696}.
\newblock
\urldef\tempurl%
\url{https://doi.org/10.1145/1830761.1830790}
\showDOI{\tempurl}


\bibitem[\protect\citeauthoryear{Hansen, Auger, Ros, Mersmann, Tu{\v s}ar, and
  Brockhoff}{Hansen et~al\mbox{.}}{2021}]%
        {HansenARMTB21}
\bibfield{author}{\bibinfo{person}{Nikolaus Hansen}, \bibinfo{person}{Anne
  Auger}, \bibinfo{person}{Raymond Ros}, \bibinfo{person}{Olaf Mersmann},
  \bibinfo{person}{Tea Tu{\v s}ar}, {and} \bibinfo{person}{Dimo Brockhoff}.}
  \bibinfo{year}{2021}\natexlab{}.
\newblock \showarticletitle{{{COCO:} a platform for comparing continuous
  optimizers in a black-box setting}}.
\newblock \bibinfo{journal}{\emph{Optim. Methods Softw.}} \bibinfo{volume}{36},
  \bibinfo{number}{1} (\bibinfo{year}{2021}), \bibinfo{pages}{114--144}.
\newblock


\bibitem[\protect\citeauthoryear{Hansen, Finck, Ros, and Auger}{Hansen
  et~al\mbox{.}}{2009}]%
        {HansenFRA09bbob}
\bibfield{author}{\bibinfo{person}{Nikolaus Hansen}, \bibinfo{person}{Steffen
  Finck}, \bibinfo{person}{Raymond Ros}, {and} \bibinfo{person}{Anne Auger}.}
  \bibinfo{year}{2009}\natexlab{}.
\newblock \bibinfo{booktitle}{\emph{{Real-Parameter Black-Box Optimization
  Benchmarking 2009: Noiseless Functions Definitions}}}.
\newblock \bibinfo{type}{{T}echnical {R}eport} RR-6829.
  \bibinfo{institution}{INRIA}.
\newblock


\bibitem[\protect\citeauthoryear{Hartikainen, Miettinen, and
  Wiecek}{Hartikainen et~al\mbox{.}}{2011}]%
        {HartikainenMW11}
\bibfield{author}{\bibinfo{person}{Markus Hartikainen}, \bibinfo{person}{Kaisa
  Miettinen}, {and} \bibinfo{person}{Margaret~M. Wiecek}.}
  \bibinfo{year}{2011}\natexlab{}.
\newblock \showarticletitle{Constructing a Pareto front approximation for
  decision making}.
\newblock \bibinfo{journal}{\emph{Math. Methods Oper. Res.}}
  \bibinfo{volume}{73}, \bibinfo{number}{2} (\bibinfo{year}{2011}),
  \bibinfo{pages}{209--234}.
\newblock
\urldef\tempurl%
\url{https://doi.org/10.1007/s00186-010-0343-0}
\showDOI{\tempurl}


\bibitem[\protect\citeauthoryear{Hartikainen, Miettinen, and
  Wiecek}{Hartikainen et~al\mbox{.}}{2012}]%
        {HartikainenMW12}
\bibfield{author}{\bibinfo{person}{Markus Hartikainen}, \bibinfo{person}{Kaisa
  Miettinen}, {and} \bibinfo{person}{Margaret~M. Wiecek}.}
  \bibinfo{year}{2012}\natexlab{}.
\newblock \showarticletitle{{PAINT:} Pareto front interpolation for nonlinear
  multiobjective optimization}.
\newblock \bibinfo{journal}{\emph{Comput. Optim. Appl.}} \bibinfo{volume}{52},
  \bibinfo{number}{3} (\bibinfo{year}{2012}), \bibinfo{pages}{845--867}.
\newblock
\urldef\tempurl%
\url{https://doi.org/10.1007/s10589-011-9441-z}
\showDOI{\tempurl}


\bibitem[\protect\citeauthoryear{Hirano and Yoshikawa}{Hirano and
  Yoshikawa}{2013}]%
        {HiranoY13}
\bibfield{author}{\bibinfo{person}{Hiroyuki Hirano} {and}
  \bibinfo{person}{Tomohiro Yoshikawa}.} \bibinfo{year}{2013}\natexlab{}.
\newblock \showarticletitle{A study on two-step search based on {PSO} to
  improve convergence and diversity for Many-Objective Optimization Problems}.
  In \bibinfo{booktitle}{\emph{Proceedings of the {IEEE} Congress on
  Evolutionary Computation, {CEC} 2013, Cancun, Mexico, June 20-23, 2013}}.
  \bibinfo{publisher}{{IEEE}}, \bibinfo{pages}{1854--1859}.
\newblock
\urldef\tempurl%
\url{https://doi.org/10.1109/CEC.2013.6557785}
\showDOI{\tempurl}


\bibitem[\protect\citeauthoryear{Hu, Yen, and Luo}{Hu et~al\mbox{.}}{2017}]%
        {HuYL17}
\bibfield{author}{\bibinfo{person}{Wang Hu}, \bibinfo{person}{Gary~G. Yen},
  {and} \bibinfo{person}{Guangchun Luo}.} \bibinfo{year}{2017}\natexlab{}.
\newblock \showarticletitle{Many-Objective Particle Swarm Optimization Using
  Two-Stage Strategy and Parallel Cell Coordinate System}.
\newblock \bibinfo{journal}{\emph{{IEEE} Trans. Cybern.}} \bibinfo{volume}{47},
  \bibinfo{number}{6} (\bibinfo{year}{2017}), \bibinfo{pages}{1446--1459}.
\newblock
\urldef\tempurl%
\url{https://doi.org/10.1109/TCYB.2016.2548239}
\showDOI{\tempurl}


\bibitem[\protect\citeauthoryear{Huband, Hingston, Barone, and While}{Huband
  et~al\mbox{.}}{2006}]%
        {HubandHBW06}
\bibfield{author}{\bibinfo{person}{Simon Huband}, \bibinfo{person}{Philip
  Hingston}, \bibinfo{person}{Luigi Barone}, {and} \bibinfo{person}{Lyndon
  While}.} \bibinfo{year}{2006}\natexlab{}.
\newblock \showarticletitle{A review of multiobjective test problems and a
  scalable test problem toolkit}.
\newblock \bibinfo{journal}{\emph{{IEEE} Trans. Evol. Comput.}}
  \bibinfo{volume}{10}, \bibinfo{number}{5} (\bibinfo{year}{2006}),
  \bibinfo{pages}{477--506}.
\newblock
\urldef\tempurl%
\url{https://doi.org/10.1109/TEVC.2005.861417}
\showDOI{\tempurl}


\bibitem[\protect\citeauthoryear{Hutter, Hoos, and Leyton{-}Brown}{Hutter
  et~al\mbox{.}}{[n.\,d.]}]%
        {HutterHL11}
\bibfield{author}{\bibinfo{person}{Frank Hutter}, \bibinfo{person}{Holger~H.
  Hoos}, {and} \bibinfo{person}{Kevin Leyton{-}Brown}.}
  \bibinfo{year}{[n.\,d.]}\natexlab{}.
\newblock \showarticletitle{{}Sequential Model-Based Optimization for General
  Algorithm Configuration}, \bibfield{editor}{\bibinfo{person}{Carlos A.~Coello
  Coello}} (Ed.).
\newblock


\bibitem[\protect\citeauthoryear{Igel, Suttorp, and Hansen}{Igel
  et~al\mbox{.}}{2006}]%
        {IgelSH06}
\bibfield{author}{\bibinfo{person}{Christian Igel}, \bibinfo{person}{Thorsten
  Suttorp}, {and} \bibinfo{person}{Nikolaus Hansen}.}
  \bibinfo{year}{2006}\natexlab{}.
\newblock \showarticletitle{Steady-State Selection and Efficient Covariance
  Matrix Update in the Multi-objective {CMA-ES}}. In
  \bibinfo{booktitle}{\emph{Evolutionary Multi-Criterion Optimization, 4th
  International Conference, {EMO} 2007, Matsushima, Japan, March 5-8, 2007,
  Proceedings}} \emph{(\bibinfo{series}{Lecture Notes in Computer Science},
  Vol.~\bibinfo{volume}{4403})}, \bibfield{editor}{\bibinfo{person}{Shigeru
  Obayashi}, \bibinfo{person}{Kalyanmoy Deb}, \bibinfo{person}{Carlo Poloni},
  \bibinfo{person}{Tomoyuki Hiroyasu}, {and} \bibinfo{person}{Tadahiko Murata}}
  (Eds.). \bibinfo{publisher}{Springer}, \bibinfo{pages}{171--185}.
\newblock
\urldef\tempurl%
\url{https://doi.org/10.1007/978-3-540-70928-2\_16}
\showDOI{\tempurl}


\bibitem[\protect\citeauthoryear{Ishibuchi, Setoguchi, Masuda, and
  Nojima}{Ishibuchi et~al\mbox{.}}{2017}]%
        {IshibuchiSMN16}
\bibfield{author}{\bibinfo{person}{Hisao Ishibuchi}, \bibinfo{person}{Yu
  Setoguchi}, \bibinfo{person}{Hiroyuki Masuda}, {and} \bibinfo{person}{Yusuke
  Nojima}.} \bibinfo{year}{2017}\natexlab{}.
\newblock \showarticletitle{{Performance of Decomposition-Based Many-Objective
  Algorithms Strongly Depends on Pareto Front Shapes}}.
\newblock \bibinfo{journal}{\emph{{IEEE} Trans. Evol. Comput.}}
  \bibinfo{volume}{21}, \bibinfo{number}{2} (\bibinfo{year}{2017}),
  \bibinfo{pages}{169--190}.
\newblock
\urldef\tempurl%
\url{https://doi.org/10.1109/TEVC.2016.2587749}
\showDOI{\tempurl}


\bibitem[\protect\citeauthoryear{Knowles}{Knowles}{2006}]%
        {Knowles06}
\bibfield{author}{\bibinfo{person}{Joshua~D. Knowles}.}
  \bibinfo{year}{2006}\natexlab{}.
\newblock \showarticletitle{ParEGO: a hybrid algorithm with on-line landscape
  approximation for expensive multiobjective optimization problems}.
\newblock \bibinfo{journal}{\emph{{IEEE} Trans. Evol. Comput.}}
  \bibinfo{volume}{10}, \bibinfo{number}{1} (\bibinfo{year}{2006}),
  \bibinfo{pages}{50--66}.
\newblock
\urldef\tempurl%
\url{https://doi.org/10.1109/TEVC.2005.851274}
\showDOI{\tempurl}


\bibitem[\protect\citeauthoryear{Kobayashi, Hamada, Sannai, Tanaka, Bannai, and
  Sugiyama}{Kobayashi et~al\mbox{.}}{2019}]%
        {KobayashiHSTBS19}
\bibfield{author}{\bibinfo{person}{Ken Kobayashi}, \bibinfo{person}{Naoki
  Hamada}, \bibinfo{person}{Akiyoshi Sannai}, \bibinfo{person}{Akinori Tanaka},
  \bibinfo{person}{Kenichi Bannai}, {and} \bibinfo{person}{Masashi Sugiyama}.}
  \bibinfo{year}{2019}\natexlab{}.
\newblock \showarticletitle{B{\'{e}}zier Simplex Fitting: Describing Pareto
  Fronts of Simplicial Problems with Small Samples in Multi-Objective
  Optimization}. In \bibinfo{booktitle}{\emph{The Thirty-Third {AAAI}
  Conference on Artificial Intelligence, {AAAI} 2019, The Thirty-First
  Innovative Applications of Artificial Intelligence Conference, {IAAI} 2019,
  The Ninth {AAAI} Symposium on Educational Advances in Artificial
  Intelligence, {EAAI} 2019, Honolulu, Hawaii, USA, January 27 - February 1,
  2019}}. \bibinfo{publisher}{{AAAI} Press}, \bibinfo{pages}{2304--2313}.
\newblock
\urldef\tempurl%
\url{https://doi.org/10.1609/aaai.v33i01.33012304}
\showDOI{\tempurl}


\bibitem[\protect\citeauthoryear{Kraft}{Kraft}{1988}]%
        {Kraft88}
\bibfield{author}{\bibinfo{person}{Dieter Kraft}.}
  \bibinfo{year}{1988}\natexlab{}.
\newblock \bibinfo{booktitle}{\emph{{A Software Package for Sequential
  Quadratic Programming}}}.
\newblock \bibinfo{type}{{T}echnical {R}eport} DFVLR-FB 88-28.
  \bibinfo{institution}{Deutsche Forschungs- und Versuchsanstalt f\"{u}r Luft-
  und Raumfahrt}.
\newblock


\bibitem[\protect\citeauthoryear{Loshchilov and Glasmachers}{Loshchilov and
  Glasmachers}{2016}]%
        {LoshchilovG16}
\bibfield{author}{\bibinfo{person}{Ilya Loshchilov} {and}
  \bibinfo{person}{Tobias Glasmachers}.} \bibinfo{year}{2016}\natexlab{}.
\newblock \showarticletitle{Anytime Bi-Objective Optimization with a Hybrid
  Multi-Objective {CMA-ES} {(HMO-CMA-ES)}}. In
  \bibinfo{booktitle}{\emph{Genetic and Evolutionary Computation Conference,
  {GECCO} 2016, Denver, CO, USA, July 20-24, 2016, Companion Material
  Proceedings}}, \bibfield{editor}{\bibinfo{person}{Tobias Friedrich},
  \bibinfo{person}{Frank Neumann}, {and} \bibinfo{person}{Andrew~M. Sutton}}
  (Eds.). \bibinfo{publisher}{{ACM}}, \bibinfo{pages}{1169--1176}.
\newblock
\urldef\tempurl%
\url{https://doi.org/10.1145/2908961.2931698}
\showDOI{\tempurl}


\bibitem[\protect\citeauthoryear{Maree, Alderliesten, and Bosman}{Maree
  et~al\mbox{.}}{2020}]%
        {MareeAB20}
\bibfield{author}{\bibinfo{person}{Stefanus~C. Maree}, \bibinfo{person}{Tanja
  Alderliesten}, {and} \bibinfo{person}{Peter A.~N. Bosman}.}
  \bibinfo{year}{2020}\natexlab{}.
\newblock \showarticletitle{Ensuring Smoothly Navigable Approximation Sets by
  B{\'{e}}zier Curve Parameterizations in Evolutionary Bi-objective
  Optimization}. In \bibinfo{booktitle}{\emph{Parallel Problem Solving from
  Nature - {PPSN} {XVI} - 16th International Conference, {PPSN} 2020, Leiden,
  The Netherlands, September 5-9, 2020, Proceedings, Part {II}}}
  \emph{(\bibinfo{series}{Lecture Notes in Computer Science},
  Vol.~\bibinfo{volume}{12270})}, \bibfield{editor}{\bibinfo{person}{Thomas
  B{\"{a}}ck}, \bibinfo{person}{Mike Preuss}, \bibinfo{person}{Andr{\'{e}}~H.
  Deutz}, \bibinfo{person}{Hao Wang}, \bibinfo{person}{Carola Doerr},
  \bibinfo{person}{Michael T.~M. Emmerich}, {and} \bibinfo{person}{Heike
  Trautmann}} (Eds.). \bibinfo{publisher}{Springer}, \bibinfo{pages}{215--228}.
\newblock
\urldef\tempurl%
\url{https://doi.org/10.1007/978-3-030-58115-2\_15}
\showDOI{\tempurl}


\bibitem[\protect\citeauthoryear{Mastroddi and Gemma}{Mastroddi and
  Gemma}{2013}]%
        {MastroddiG2013}
\bibfield{author}{\bibinfo{person}{Franco Mastroddi} {and}
  \bibinfo{person}{Stefania Gemma}.} \bibinfo{year}{2013}\natexlab{}.
\newblock \showarticletitle{Analysis of {Pareto} frontiers for
  multidisciplinary design optimization of aircraft}.
\newblock \bibinfo{journal}{\emph{Aerospace Science and Technology}}
  \bibinfo{volume}{28}, \bibinfo{number}{1} (\bibinfo{year}{2013}),
  \bibinfo{pages}{40--55}.
\newblock
\urldef\tempurl%
\url{https://doi.org/10.1016/j.ast.2012.10.003}
\showDOI{\tempurl}


\bibitem[\protect\citeauthoryear{Miettinen}{Miettinen}{1998}]%
        {Miettinen98}
\bibfield{author}{\bibinfo{person}{Kaisa Miettinen}.}
  \bibinfo{year}{1998}\natexlab{}.
\newblock \bibinfo{booktitle}{\emph{Nonlinear Multiobjective Optimization}}.
\newblock \bibinfo{publisher}{Springer}.
\newblock


\bibitem[\protect\citeauthoryear{Ozaki, Tanigaki, Watanabe, and Onishi}{Ozaki
  et~al\mbox{.}}{2020}]%
        {OzakiTWO20}
\bibfield{author}{\bibinfo{person}{Yoshihiko Ozaki}, \bibinfo{person}{Yuki
  Tanigaki}, \bibinfo{person}{Shuhei Watanabe}, {and} \bibinfo{person}{Masaki
  Onishi}.} \bibinfo{year}{2020}\natexlab{}.
\newblock \showarticletitle{Multiobjective tree-structured parzen estimator for
  computationally expensive optimization problems}. In
  \bibinfo{booktitle}{\emph{{GECCO} '20: Genetic and Evolutionary Computation
  Conference, Canc{\'{u}}n Mexico, July 8-12, 2020}},
  \bibfield{editor}{\bibinfo{person}{Carlos Artemio~Coello Coello}} (Ed.).
  \bibinfo{publisher}{{ACM}}, \bibinfo{pages}{533--541}.
\newblock
\urldef\tempurl%
\url{https://doi.org/10.1145/3377930.3389817}
\showDOI{\tempurl}


\bibitem[\protect\citeauthoryear{Paquete and St{\"{u}}tzle}{Paquete and
  St{\"{u}}tzle}{2003}]%
        {PaqueteS03}
\bibfield{author}{\bibinfo{person}{Lu{\'{\i}}s Paquete} {and}
  \bibinfo{person}{Thomas St{\"{u}}tzle}.} \bibinfo{year}{2003}\natexlab{}.
\newblock \showarticletitle{A Two-Phase Local Search for the Biobjective
  Traveling Salesman Problem}. In \bibinfo{booktitle}{\emph{Evolutionary
  Multi-Criterion Optimization, Second International Conference, {EMO} 2003,
  Faro, Portugal, April 8-11, 2003, Proceedings}}
  \emph{(\bibinfo{series}{Lecture Notes in Computer Science},
  Vol.~\bibinfo{volume}{2632})}, \bibfield{editor}{\bibinfo{person}{Carlos~M.
  Fonseca}, \bibinfo{person}{Peter~J. Fleming}, \bibinfo{person}{Eckart
  Zitzler}, \bibinfo{person}{Kalyanmoy Deb}, {and} \bibinfo{person}{Lothar
  Thiele}} (Eds.). \bibinfo{publisher}{Springer}, \bibinfo{pages}{479--493}.
\newblock
\urldef\tempurl%
\url{https://doi.org/10.1007/3-540-36970-8\_34}
\showDOI{\tempurl}


\bibitem[\protect\citeauthoryear{Pos{\'{\i}}k and Huyer}{Pos{\'{\i}}k and
  Huyer}{2012}]%
        {PosikH12}
\bibfield{author}{\bibinfo{person}{Petr Pos{\'{\i}}k} {and}
  \bibinfo{person}{Waltraud Huyer}.} \bibinfo{year}{2012}\natexlab{}.
\newblock \showarticletitle{Restarted Local Search Algorithms for Continuous
  Black Box Optimization}.
\newblock \bibinfo{journal}{\emph{Evol. Comput.}} \bibinfo{volume}{20},
  \bibinfo{number}{4} (\bibinfo{year}{2012}), \bibinfo{pages}{575--607}.
\newblock
\urldef\tempurl%
\url{https://doi.org/10.1162/EVCO\_a\_00087}
\showDOI{\tempurl}


\bibitem[\protect\citeauthoryear{Powell}{Powell}{2008}]%
        {Powell08}
\bibfield{author}{\bibinfo{person}{M.~J.~D. Powell}.}
  \bibinfo{year}{2008}\natexlab{}.
\newblock \showarticletitle{Developments of NEWUOA for minimization without
  derivatives}.
\newblock \bibinfo{journal}{\emph{IMA J. Numer. Anal.}} \bibinfo{volume}{28},
  \bibinfo{number}{4} (\bibinfo{year}{2008}), \bibinfo{pages}{649--664}.
\newblock
\urldef\tempurl%
\url{https://doi.org/10.1093/imanum/drm047}
\showDOI{\tempurl}


\bibitem[\protect\citeauthoryear{Powell}{Powell}{2009}]%
        {Powell09}
\bibfield{author}{\bibinfo{person}{M.~J.~D. Powell}.}
  \bibinfo{year}{2009}\natexlab{}.
\newblock \bibinfo{booktitle}{\emph{{The BOBYQA algorithm for bound constrained
  optimization without derivatives}}}.
\newblock \bibinfo{type}{{T}echnical {R}eport} DAMTP 2009/NA06.
  \bibinfo{institution}{University of Cambridge}.
\newblock


\bibitem[\protect\citeauthoryear{Regis}{Regis}{2021}]%
        {Regis21}
\bibfield{author}{\bibinfo{person}{Rommel~G. Regis}.}
  \bibinfo{year}{2021}\natexlab{}.
\newblock \showarticletitle{A two-phase surrogate approach for high-dimensional
  constrained discrete multi-objective optimization}. In
  \bibinfo{booktitle}{\emph{{GECCO} '21: Genetic and Evolutionary Computation
  Conference, Companion Volume, Lille, France, July 10-14, 2021}},
  \bibfield{editor}{\bibinfo{person}{Krzysztof Krawiec}} (Ed.).
  \bibinfo{publisher}{{ACM}}, \bibinfo{pages}{1870--1878}.
\newblock
\urldef\tempurl%
\url{https://doi.org/10.1145/3449726.3463204}
\showDOI{\tempurl}


\bibitem[\protect\citeauthoryear{Rios and Sahinidis}{Rios and
  Sahinidis}{2013}]%
        {RiosS13}
\bibfield{author}{\bibinfo{person}{Luis~Miguel Rios} {and}
  \bibinfo{person}{Nikolaos~V. Sahinidis}.} \bibinfo{year}{2013}\natexlab{}.
\newblock \showarticletitle{Derivative-free optimization: a review of
  algorithms and comparison of software implementations}.
\newblock \bibinfo{journal}{\emph{J. Glob. Optim.}} \bibinfo{volume}{56},
  \bibinfo{number}{3} (\bibinfo{year}{2013}), \bibinfo{pages}{1247--1293}.
\newblock
\urldef\tempurl%
\url{https://doi.org/10.1007/s10898-012-9951-y}
\showDOI{\tempurl}


\bibitem[\protect\citeauthoryear{Shoval, Sheftel, Shinar, Hart, Ramote, Mayo,
  Dekel, Kavanagh, and Alon}{Shoval et~al\mbox{.}}{2012}]%
        {ShovalSSHRNDKA2012}
\bibfield{author}{\bibinfo{person}{O. Shoval}, \bibinfo{person}{H. Sheftel},
  \bibinfo{person}{G. Shinar}, \bibinfo{person}{Y. Hart}, \bibinfo{person}{O.
  Ramote}, \bibinfo{person}{A. Mayo}, \bibinfo{person}{E. Dekel},
  \bibinfo{person}{K. Kavanagh}, {and} \bibinfo{person}{U. Alon}.}
  \bibinfo{year}{2012}\natexlab{}.
\newblock \showarticletitle{Evolutionary Trade-Offs, {Pareto} Optimality, and
  the Geometry of Phenotype Space}.
\newblock \bibinfo{journal}{\emph{Science}} \bibinfo{volume}{336},
  \bibinfo{number}{6085} (\bibinfo{year}{2012}), \bibinfo{pages}{1157--1160}.
\newblock
\showISSN{0036-8075}
\urldef\tempurl%
\url{https://doi.org/10.1126/science.1217405}
\showDOI{\tempurl}
\showeprint{http://science.sciencemag.org/content/336/6085/1157.full.pdf}


\bibitem[\protect\citeauthoryear{Song, Wang, He, and Jin}{Song
  et~al\mbox{.}}{2021}]%
        {SongWHJ22}
\bibfield{author}{\bibinfo{person}{Zhenshou Song}, \bibinfo{person}{Handing
  Wang}, \bibinfo{person}{Cheng He}, {and} \bibinfo{person}{Yaochu Jin}.}
  \bibinfo{year}{2021}\natexlab{}.
\newblock \showarticletitle{A Kriging-Assisted Two-Archive Evolutionary
  Algorithm for Expensive Many-Objective Optimization}.
\newblock \bibinfo{journal}{\emph{{IEEE} Trans. Evol. Comput.}}
  \bibinfo{volume}{25}, \bibinfo{number}{6} (\bibinfo{year}{2021}),
  \bibinfo{pages}{1013--1027}.
\newblock
\urldef\tempurl%
\url{https://doi.org/10.1109/TEVC.2021.3073648}
\showDOI{\tempurl}


\bibitem[\protect\citeauthoryear{Tabatabaei, Hakanen, Hartikainen, Miettinen,
  and Sindhya}{Tabatabaei et~al\mbox{.}}{2015}]%
        {TabatabaeiHHMS15}
\bibfield{author}{\bibinfo{person}{Mohammad Tabatabaei}, \bibinfo{person}{Jussi
  Hakanen}, \bibinfo{person}{Markus Hartikainen}, \bibinfo{person}{Kaisa
  Miettinen}, {and} \bibinfo{person}{Karthik Sindhya}.}
  \bibinfo{year}{2015}\natexlab{}.
\newblock \showarticletitle{A survey on handling computationally expensive
  multiobjective optimization problems using surrogates: non-nature inspired
  methods}.
\newblock \bibinfo{journal}{\emph{Struct. Multidiscipl. Optim.}}
  \bibinfo{volume}{52} (\bibinfo{year}{2015}), \bibinfo{pages}{1--25}.
\newblock


\bibitem[\protect\citeauthoryear{Tanabe and Ishibuchi}{Tanabe and
  Ishibuchi}{2020}]%
        {TanabeI2020}
\bibfield{author}{\bibinfo{person}{Ryoji Tanabe} {and} \bibinfo{person}{Hisao
  Ishibuchi}.} \bibinfo{year}{2020}\natexlab{}.
\newblock \showarticletitle{An easy-to-use real-world multi-objective
  optimization problem suite}.
\newblock \bibinfo{journal}{\emph{Applied Soft Computing}}
  \bibinfo{volume}{89} (\bibinfo{year}{2020}), \bibinfo{pages}{106078}.
\newblock
\urldef\tempurl%
\url{https://doi.org/10.1016/j.asoc.2020.106078}
\showDOI{\tempurl}


\bibitem[\protect\citeauthoryear{Tian, Cheng, Zhang, and Jin}{Tian
  et~al\mbox{.}}{2017}]%
        {TianCZJ17}
\bibfield{author}{\bibinfo{person}{Ye Tian}, \bibinfo{person}{Ran Cheng},
  \bibinfo{person}{Xingyi Zhang}, {and} \bibinfo{person}{Yaochu Jin}.}
  \bibinfo{year}{2017}\natexlab{}.
\newblock \showarticletitle{PlatEMO: {A} {MATLAB} Platform for Evolutionary
  Multi-Objective Optimization [Educational Forum]}.
\newblock \bibinfo{journal}{\emph{{IEEE} Comput. Intell. Mag.}}
  \bibinfo{volume}{12}, \bibinfo{number}{4} (\bibinfo{year}{2017}),
  \bibinfo{pages}{73--87}.
\newblock
\urldef\tempurl%
\url{https://doi.org/10.1109/MCI.2017.2742868}
\showDOI{\tempurl}


\bibitem[\protect\citeauthoryear{Tour{\'{e}}, Hansen, Auger, and
  Brockhoff}{Tour{\'{e}} et~al\mbox{.}}{2019}]%
        {ToureHAB19}
\bibfield{author}{\bibinfo{person}{Cheikh Tour{\'{e}}},
  \bibinfo{person}{Nikolaus Hansen}, \bibinfo{person}{Anne Auger}, {and}
  \bibinfo{person}{Dimo Brockhoff}.} \bibinfo{year}{2019}\natexlab{}.
\newblock \showarticletitle{Uncrowded hypervolume improvement: {COMO-CMA-ES}
  and the sofomore framework}. In \bibinfo{booktitle}{\emph{Proceedings of the
  Genetic and Evolutionary Computation Conference, {GECCO} 2019, Prague, Czech
  Republic, July 13-17, 2019}}, \bibfield{editor}{\bibinfo{person}{Anne Auger}
  {and} \bibinfo{person}{Thomas St{\"{u}}tzle}} (Eds.).
  \bibinfo{publisher}{{ACM}}, \bibinfo{pages}{638--646}.
\newblock
\urldef\tempurl%
\url{https://doi.org/10.1145/3321707.3321852}
\showDOI{\tempurl}


\bibitem[\protect\citeauthoryear{Vrugt, Gupta, Bastidas, Bouten, and
  Sorooshian}{Vrugt et~al\mbox{.}}{2003}]%
        {VrugtGBBS2003}
\bibfield{author}{\bibinfo{person}{Jasper~A. Vrugt}, \bibinfo{person}{Hoshin~V.
  Gupta}, \bibinfo{person}{Luis~A. Bastidas}, \bibinfo{person}{Willem Bouten},
  {and} \bibinfo{person}{Soroosh Sorooshian}.} \bibinfo{year}{2003}\natexlab{}.
\newblock \showarticletitle{Effective and Efficient Algorithm for
  Multiobjective Optimization of Hydrologic Models}.
\newblock \bibinfo{journal}{\emph{Water Resources Research}}
  \bibinfo{volume}{39}, \bibinfo{number}{8} (\bibinfo{year}{2003}),
  \bibinfo{pages}{1214--1232}.
\newblock
\urldef\tempurl%
\url{https://doi.org/10.1029/2002WR001746}
\showDOI{\tempurl}


\bibitem[\protect\citeauthoryear{Wang, Hong, Ye, Zhang, Jiang, and Tan}{Wang
  et~al\mbox{.}}{ress}]%
        {WangHYZJT22}
\bibfield{author}{\bibinfo{person}{Zhenzhong Wang}, \bibinfo{person}{Haokai
  Hong}, \bibinfo{person}{Kai Ye}, \bibinfo{person}{Guang-En Zhang},
  \bibinfo{person}{Min Jiang}, {and} \bibinfo{person}{Kay~Chen Tan}.}
  \bibinfo{year}{2022 (in press)}\natexlab{}.
\newblock \showarticletitle{Manifold Interpolation for Large-Scale
  Multiobjective Optimization via Generative Adversarial Networks}.
\newblock \bibinfo{journal}{\emph{{IEEE} Trans. Neural Networks Learn. Syst.}}
  (\bibinfo{year}{2022 (in press)}).
\newblock


\bibitem[\protect\citeauthoryear{Yang, Palar, Emmerich, Shimoyama, and
  B{\"{a}}ck}{Yang et~al\mbox{.}}{2019}]%
        {YangPESB19}
\bibfield{author}{\bibinfo{person}{Kaifeng Yang},
  \bibinfo{person}{Pramudita~Satria Palar}, \bibinfo{person}{Michael Emmerich},
  \bibinfo{person}{Koji Shimoyama}, {and} \bibinfo{person}{Thomas B{\"{a}}ck}.}
  \bibinfo{year}{2019}\natexlab{}.
\newblock \showarticletitle{A multi-point mechanism of expected hypervolume
  improvement for parallel multi-objective bayesian global optimization}. In
  \bibinfo{booktitle}{\emph{Proceedings of the Genetic and Evolutionary
  Computation Conference, {GECCO} 2019, Prague, Czech Republic, July 13-17,
  2019}}, \bibfield{editor}{\bibinfo{person}{Anne Auger} {and}
  \bibinfo{person}{Thomas St{\"{u}}tzle}} (Eds.). \bibinfo{publisher}{{ACM}},
  \bibinfo{pages}{656--663}.
\newblock
\urldef\tempurl%
\url{https://doi.org/10.1145/3321707.3321784}
\showDOI{\tempurl}


\bibitem[\protect\citeauthoryear{Zhang and Li}{Zhang and Li}{2007}]%
        {ZhangL07}
\bibfield{author}{\bibinfo{person}{Qingfu Zhang} {and} \bibinfo{person}{Hui
  Li}.} \bibinfo{year}{2007}\natexlab{}.
\newblock \showarticletitle{{MOEA/D:} {A} Multiobjective Evolutionary Algorithm
  Based on Decomposition}.
\newblock \bibinfo{journal}{\emph{{IEEE} Trans. Evol. Comput.}}
  \bibinfo{volume}{11}, \bibinfo{number}{6} (\bibinfo{year}{2007}),
  \bibinfo{pages}{712--731}.
\newblock
\urldef\tempurl%
\url{https://doi.org/10.1109/TEVC.2007.892759}
\showDOI{\tempurl}


\bibitem[\protect\citeauthoryear{Zhang, Zhou, and Jin}{Zhang
  et~al\mbox{.}}{2008}]%
        {ZhangZJ08}
\bibfield{author}{\bibinfo{person}{Qingfu Zhang}, \bibinfo{person}{Aimin Zhou},
  {and} \bibinfo{person}{Yaochu Jin}.} \bibinfo{year}{2008}\natexlab{}.
\newblock \showarticletitle{{RM-MEDA:} {A} Regularity Model-Based
  Multiobjective Estimation of Distribution Algorithm}.
\newblock \bibinfo{journal}{\emph{{IEEE} Trans. Evol. Comput.}}
  \bibinfo{volume}{12}, \bibinfo{number}{1} (\bibinfo{year}{2008}),
  \bibinfo{pages}{41--63}.
\newblock
\urldef\tempurl%
\url{https://doi.org/10.1109/TEVC.2007.894202}
\showDOI{\tempurl}


\bibitem[\protect\citeauthoryear{Zitzler and Thiele}{Zitzler and
  Thiele}{1998}]%
        {ZitzlerT98}
\bibfield{author}{\bibinfo{person}{Eckart Zitzler} {and}
  \bibinfo{person}{Lothar Thiele}.} \bibinfo{year}{1998}\natexlab{}.
\newblock \showarticletitle{Multiobjective Optimization Using Evolutionary
  Algorithms - {A} Comparative Case Study}. In
  \bibinfo{booktitle}{\emph{Parallel Problem Solving from Nature - {PPSN} V,
  5th International Conference, Amsterdam, The Netherlands, September 27-30,
  1998, Proceedings}} \emph{(\bibinfo{series}{Lecture Notes in Computer
  Science}, Vol.~\bibinfo{volume}{1498})},
  \bibfield{editor}{\bibinfo{person}{A.~E. Eiben}, \bibinfo{person}{Thomas
  B{\"{a}}ck}, \bibinfo{person}{Marc Schoenauer}, {and}
  \bibinfo{person}{Hans{-}Paul Schwefel}} (Eds.).
  \bibinfo{publisher}{Springer}, \bibinfo{pages}{292--304}.
\newblock
\urldef\tempurl%
\url{https://doi.org/10.1007/BFb0056872}
\showDOI{\tempurl}


\end{thebibliography}










\end{document}